\input amstex
\documentstyle{amsppt}
\magnification=\magstep1
\define\cc{\Bbb C}
 \define\z{\Bbb Z}
\define\r{\Bbb R}
\define\k{\Bbb K}
\define\N{\Bbb N}
\define\jj{\Bbb J}

\define\Q{\Bbb Q}
\define\A{\Cal A}
\define\h{\Cal D}
\define\E{\Cal E}
\define\m{\Cal M}
\define\T{\Cal T}
\define\f{\Cal S}

\define \W{\Cal W}

\define\p{\Cal P}
\define\apf{almost periodic }

\define\la{\lambda}
\define\om{\omega}

\define\e{\varepsilon}
\define\va{\varphi }
\define\CB#1{\Cal C_b(#1)}
\define\st{\subset }
\define\al{\alpha  }
\topmatter
  \title
 Generalized  vector valued almost periodic and ergodic distributions
 \endtitle
  \endkeywords
  \author
  Bolis Basit and Hans G\"unzler
\endauthor
 \abstract
{For $\Cal A\subset L^1_{loc}(\Bbb J,X)$ let $\Cal M\Cal A$
consist of all $f\in L^1_{loc}$ with $ M_h f (\cdot):=\frac
{1}{h}\int _{0}^{h}f(\cdot +s)\,ds  \in \Cal A$ for all $h>0$.
Here $X$ is a Banach space, $\Bbb J= (\alpha ,\infty),  [\alpha
,\infty)$ or  $\Bbb R$. The class $\widetilde { \Cal M} \Cal A
\subset \Cal D'(\Bbb J, X)$ the space of distributions is
similarly defined. Usually $\Cal A\subset\Cal M\Cal A\subset \Cal
M^2\Cal A\subset \cdots$. One has
 $  \cup_{n=0}^{\infty} \widetilde { \Cal M}^n\Cal {A}= \Cal D'_{\Cal A} (\Bbb {R},X) :=\{T\in   \Cal {D}'(\Bbb {R},X):  T*\varphi\in \Cal A  \text { \, for\,\,all \,}\varphi \in \Cal {D}(\Bbb {R},\Bbb C)  \}$,
 $ \Cal D'_{\Cal A} (\Bbb {R},X)\cap  L^1_{loc} (\Bbb {R},X)  =\cup_{n=0}^{\infty} \Cal {M}^n\Cal A$.
 The map  $ \Cal A \to   \Cal {D}'_{\Cal A}$
is iteration complete, that is $  \Cal {D}'_{ \Cal {D}'_{\Cal A}}=
\Cal {D}'_{\Cal A}$. Under suitable assumptions  $\widetilde {
\Cal M}^n \Cal {A}= \Cal A + \{T^{(n)} : T \in \Cal A\}$, and
similarly for $\Cal {M}^n \Cal A$.  Almost periodic $X$-valued
distributions $\h'_{\A}$ with $\A = $ almost
 periodic (ap) functions are characterized in several ways.
 Various generalizations of the
Bohl-Bohr-Kadets theorem on the almost periodicity of the
indefinite  integral of an ap or almost automorphic function are
obtained.
 On $ \Cal {D}'_{\Cal E} $, $ \Cal E $ the class of  ergodic functions, a mean can be
constructed which gives Fourier series. Special cases of  $\Cal A
$ are the  Bohr  ap, Stepanoff ap, almost automorphic,
asymptotically ap, Eberlein weakly ap, pseudo ap and (totally)
ergodic functions $(\T)\E$.
 Then always $\Cal {M}^n \Cal A$  is strictly contained in  $ \Cal {M}^{n+1} \Cal
A$.   The relations between $\m^n \E$, $\m^n\T\E$ and subclasses
are discussed.  For many of the above results a new
$(\Delta)$-condition is needed,
       we show that it holds for most of the $\A$ needed in applications.
 Also,   we obtain    new tauberian theorems
  for $f\in L^1_{loc}(\Bbb J,X)$ to belong to a class $\A$ which are decisive
in describing the    asymptotic behavior of unbounded solutions of many  abstract differential-integral equations.
This generalizes various recent results}.
  \endabstract
 \endtopmatter
\rightheadtext{ Generalized almost periodic distributions}
\leftheadtext{Bolis Basit and Hans G\"unzler}
 \TagsOnRight
\document
\pageno=1
 \baselineskip=15pt


\newpage

\qquad\qquad\qquad\qquad\qquad \qquad\qquad { \text {Contents}}

\bigskip

\S
0.Introduction...................................................................................................3

\S 1.  Notation,  Definitions and
Preliminaries........................................................6

\S 2. New classes of  vector valued
distributions....................................................17

\S 3.   Mean classes.............................................
.................................................26

\S 4.   The $(\Delta)$
Condition......................................................................................30

 \S 5.  Mean classes, derivatives and uniform continuity
........................................39

\S 6. Almost periodic distributions  and indefinite integrals
.................................42

\S 7.  Ergodic
classes............................................................................................49

\S 8. Fourier
analysis...........................................................................................54

\S 9. Applications to the study of asymptotic behavior of
solutions of

\qquad differential
equations.................................................................................59

\S 10 Differential-difference  equations
................................................................64

\S 11 open
questions............................................................................................62

References.........................................................................................................67
\newpage

 \head{\bf\S 0. Introduction}\endhead

The primary objective of these notes is to present a
self-contained study of the theory of mean type classes introduced
during our study of the asymptotic behavior of solutions of
differential-difference systems, [13, (1.1)]. This theory does
what the Stepanoff extension of Bohr's almost periodic functions
accomplished (see below), now for quite general function classes,
for example  asymptotic almost periodic, almost automorphic or
Eberlein weakly almost periodic functions:
 For a
  given
 $\A \st L^1_{loc}(\jj,X)$, $X$ Banach space, $\jj$
infinite interval,

 $\m\A= \{ f\in L^1_{loc}:$ the mean  $M_h f(\cdot) := (1/h) \int_{0}^{h}
f(\cdot+s)\,ds \in \A$ for all $h>0 \}$.

\noindent This extension $\m\A$ of  $\A$  can also be defined for
$U $  a subset of Schwartz's vector valued distributions
$\h'(\jj,X)$, giving $\widetilde{\m} U$.

 These extensions proved
useful in [13].
 There they made
possible a considerable weakening of the necessary assumptions,
and simpler and more flexible formulations. For example bounded
solutions of differential-difference systems may not belong to the
class of Banach  space valued Bohr almost periodic (ap) functions
$AP=AP(\r,X)$, but belong to $\m AP$.  This means that the study
of mean classes makes an important contribution to the
understanding of asymptotic behavior of solutions of many
differential-difference systems. Also, this
 gives rise to an
extension of many known spaces. For example for   the class $AP$,
one has $AP\st S^p AP =$
 $\{$ Stepanoff  almost periodic functions$\}$  $\st \m AP$.

In the following we give a short account of how these classes
evolved starting with Bohr's  almost periodic functions.

In his study of the asymptotic behavior of solutions of
differential equations in perturbation theory of astronomy, Bohl
[26] 1906 introduced quasi-periodic functions as (essentially)
sums of finitely many complex-valued periodic functions on the
real line $\r$ with arbitrary, but  usually rationally independent
periods. The study of his differential equations could be reduced
to the question: When is the indefinite integral $F$ of a
quasi-periodic $f$ again quasi-periodic? As a necessary and
sufficient condition he obtained the boundedness of $F$ on $\r$.

The class of quasi-periodic functions has been extended by Bohr
[28] 1925 to his almost periodic (ap) functions, which can be
characterized as uniform limits on $\r$ of trigonometric
polynomials $\sum_1 ^n  c_j e^{i \omega_j \,t}$, with $c_j$
complex and $\omega_j$ arbitrary real. Again the indefinite
integral of an ap function is ap provided it is  bounded (theorem
of Bohl-Bohr), with corresponding applications to linear
differential equations.

To include also discontinuous only locally integrable periodic
functions, Stepanoff [73] in 1926 introduced what are now called
Stepanoff-ap functions, limits of trigonometric polynomials with
respect to the norm  $||f||_{S^1}  : =$  sup $\{ \int _t  ^{t+1}
|f(s)|\ ds : t$ real $\}$. Bochner noticed that a bounded
indefinite integral $F$ of $f$ is already Bohr-ap if f is
Stepanoff-ap. Doss [37] 1961 extended this  in the following way :
If $ F : \r \to  \cc$ is bounded and if the differences $\Delta_h
F  = F_h - F$  are Bohr-ap for all $h>0$, with $F_h (t) : =
F(t+h)$, then $F$ is ap. In the case where $F$  is the indefinite
integral $Pf$ of $f$,  this gives: If all means $M_h f$ are ap for
$h >0$ and if $Pf$ is bounded, $Pf$ is ap.

\noindent   For
 $\m AP(\r,\cc)$, one gets:
$f \in   \m AP$ and $Pf$ bounded implies  $Pf \in AP$.

\noindent  This is an extension of the above mentioned Bohl-Bohr
theorem, since one has $AP \st \{$ Stepanoff-ap functions$ \} \st
\m AP$, and both inclusions are strict. Such an extension can now
be carried over to solutions of linear differential-difference
systems ([13, e.g. Theorem 4.1]), again instead of the right side
$f$ being Bohr-ap it suffices to assume only $f  \in  \m AP$.

Non-periodic ap functions appear not only in the above mentioned
perturbation problem.  For example all bounded (non-resonant)
solutions of the boundary value problem for the wave equation are
ap in the time variable ([2, Chapter 5]).

To treat partial differential equations, one has to consider ap
functions with values in a Banach space $X$. Furthermore there are
other types of interesting asymptotic behavior, for example
Bochner's almost automorphic functions [24], Fr\'{e}chet's
asymptotically ap  functions ($ f+g$, with $f$  ap  and  $ g(t)
\to  0$  as  $|t| \to \infty$)[43], Eberlein ap functions [40] or
the recently introduced pseudo-ap functions of Zhang [77] 1994.

For a discussion of such function classes and especially for
harmonic analysis and synthesis it is essential to have a mean
$lim _{T \to \infty} (1/2T) \int _{-T} ^T  f (s)\, ds$ , and
indeed  all these functions are ergodic in some sense.

With all the above and especially applications to differential
equations in mind, starting with an arbitrary class $\A$ of
functions $f : \jj \to X $ which are only locally integrable,
$\jj$ any infinite interval from $\r$ containing  some
$[n,\infty)$, $X$ a Banach space, we introduce as above the
classes $\m\A$, $\m^2\A = \m(\m \A$),..... Usually (but not
always)  one has $ \A \st \m \A$ and then $\A \st \m\A \st \m^2
\cdots $; for $\A = AP$, or a similar class, or $\A = \{$ ergodic
functions $\}$ all these inclusions are strict, giving a whole
hierarchy of generalized $\A$-functions. Since the Bohr-mean can
be extended to ergodic distributions (\S 8), harmonic analysis
becomes possible for all the $\m^n\A$ and $\widetilde {\m}^n \A$,
as soon as it is possible for $\A$.

Consider  the extensions $\widetilde{\m} U$ of a subset $U$ of
Schwartz's vector valued distributions $ \h'(\jj,X)$. The space
$\h'_U(\r,X)$ consists of those $T  \in  \h'(\r,X)$ for which
convolution $T*\varphi \in U$ for all test functions $\varphi \in
\h'(\r,\cc)$. One gets  $\cup_0^{\infty} \widetilde{\m}^n U =
\h'_U(\r,X)$. With $U=AP$ and $X=\cc$ Schwartz's [71, p. 206]
${\Cal {B}}'_{pp}(\r,\cc)= \h'_{AP}(\r,\cc)$.

If a function $f  \in \A$ has a derivative $f'$, then almost by
definition  $f'\in \m\A$. Under some additional assumptions one
can show the converse: $\m\A = \{ f + g' : f, g  \in  \A, g'$
"exists" $\}$, similarly for the distribution class
$\widetilde{\m} U$. Here and for several other results the
essential assumption is that $\A$ satisfies the
"$(\Delta)$-condition", which means that if for a locally
integrable  $f : \jj \to X$  the differences $ \Delta_h f   \in
\A$ for all $h>0$ then also $ f - (1/s) \int_t^{t+s} f (v)\, dv
\in \A$ for all $s>0$. Fortunately practically all the important
function spaces satisfy this
 condition (\S 4).

The introduction of these mean classes $\m^n\A$ not only makes
possible the generalization and sharper formulation  of results on
the asymptotic behavior of solutions of linear
differential-difference systems ([13]), it also collects and
unifies diverse results in the theory of ap functions and
analogues, as exemplified in Proposition 5.6 below, furthermore
some new results on the global behavior of derivatives can be
obtained thereby (\S 10). Finally with these concepts an extension
of a tauberian criterion  essential in the study of differential
equations to also unbounded functions is possible (\S 9).

We give now a detailed description of the contents of the
following sections:

In \S 1 we introduce the various function spaces treated,
especially ap functions, their analogues and generalizations, and
ergodic functions. Also we give here some characterizations of
functions with (uniformly) continuous differences  needed later.


 In section 2 we study the extension
$\A\st\m\A\st \m^2\A\st \cdots$  systematically, as well as a
corresponding  one for
 $\widetilde { \m} \A $ for classes $\A$ of vector-valued  distributions from
$ \h'(\jj,X)$, completing and generalizing results from [13,
Section 3] and Schwartz [71, pp. 206-207, 247-248]. It turns out
that under weak assumptions all inclusions are strict (Proposition
3.8) and
 $ \m^n \A $, $\widetilde { \m}^n \A $  are subsumed by ``$\A$-distributions'' $ \h'_{\A}$   consisting of   all  $T  \in \h' (\r,X)$ with $  T*\varphi\in \A $, for all test functions  $\varphi \in \h(\r,\cc)$:
 $  \h'_{\A}= \cup_{n=0}^{\infty} \widetilde { \m}^n\A,\,$
 $ \h'_{\A} (\r,X)\cap  L^1_{loc} (\r,X)  =\cup_{n=0}^{\infty} \m^n\A,\,$
and $  \h'_{ \h'_{\A}}= \h'_{\A}$.
 We give various characterizations of  $\h'_{\A}$    and similar ones needed for
$ \f'(\r,X)$ (Theorems 2.11, 2.16).

  In Sections 3, 4 we give conditions under which the various assumptions needed for Section 2 are fulfilled,
especially $A \st { \m} \A $ and the central condition  $(\Delta)$
of Definition 1.4 below. For practically all  the spaces of
generalized almost periodic or ergodic functions $\E$ they hold
(Examples 3.4-3.7, 4.7,  4.15, 4.16). In Propositions 4.9-4.10 we
show that both $L^p_w (\jj,X)$ and   $O(w) (\jj,X)$ satisfy
$(\Delta)$. In Proposition 3.8, Examples 3.9 and Corollary 3.10 we
demonstrate that the inclusions $\m^k A \st { \m^{k+1}} \A $,
$k\in \N_0$ are strict for various classes $\A$. See also
Proposition 7.1 and Remark 7.2.

In Section 5 we show   $\widetilde {\m}^n \A = \A + \{ T^{(n)} : T
\in \A\}$ and similarly for
 $ { \m^n} \A $ (Corollaries 5.2- 5.4).
With this any $S^p$-ap $\phi$ is of the form $f+g'$, $f,g$ ap
(Corollary 5.4). Results of Stepanoff and Bochner are united and
generalized to $ \m^n \A \cap C_{u } \st  \A$ (Proposition 5.6).
Furthermore, it is shown that  $\m AP(\r,\r)|\, [0,\infty)$ is
strictly contained in $ \m [AP([0,\infty),\r)]$ (Example 5.5)).

In Theorem 6.1  we characterize the distributions from $\h'_{AP}=
\h'_{S^p AP}$ by Bohr's, Bochner's  and von Neumann's definitions.
Further relations between  $\A= $ Bohr-ap  functions $AP$,
asymptotic ap functions $AAP$, Eberlein ap functions $EAP$,  Weyl
ap functions $W^p AP$, Besicovitch  ap functions  $B^p AP$ and
their distribution spaces $\h'_{\A}$ can be found above
Proposition 6.5.

In various places we need and derive results of the form `` if all
the differences $\Delta_h f$
 are in $U$, then $f$ is in $V$" (Propositions 1.5, 1.7, Lemmas 4.4, 4.19),
extending several earlier results  [11], [30], [48], [8]. As
applications we get several extensions (Proposition 6.5, Corollary
6.6 and after Corollary 6.6) of the Bohl-Bohr-Kadets theorem: if
for ap $f:\r\to X$ the indefinite integral $Pf$ is bounded, it is
ap too, provided $c_0\not \st X$ [54].
 These results are not covered by the theorems of [13].

 Further  we obtained several
new properties  of
 the classes of ergodic functions $\E$, $\h'_{\E}$,  needed in  applications
 (Proposition 7.1, Remark 7.2). Especially, we prove $(\Delta)$
 for $\E$, $\E_0 =\{f\in \E$ with mean $ (f)= 0 \}$, $Av$, $Av_0$.

 In Section 8 the ergodic mean is extended to  $\widetilde {m} : \h'_{\E}\to X$.
If  $\A$ is  invariant with respect to multiplication by
characters $e^{i\omega t}$, so  also are  $ { \m^n} \A $,
$\h'_{\A}$. Hence Fourier series  and Bochner spectra may be
defined for functions and distributions in these classes. For
distributions from   $\h'_{\A}$ where   ${\A}$  is  Bohr ap,
asymptotically ap
 or Eberlein  weakly ap  functions  the classical results carry over,
 for  $ { \m}^n \A $  with a suitable topology (see after (8.7)).

  In Section 9  we introduce
the notion of spectrum of a function $\phi\in \h'_{
L^{\infty}}(\r,X)$  with respect to a class $\A (\jj,X)$ and
obtain  tauberian criteria (Proposition 9.5, Theorem 9.7,
Corollary 9.8, Remark 9.11) which play a decisive role in
describing the asymptotic behavior of also unbounded solutions of
many differential equations, and with the Loomis condition $(L_U)$
a unified treatment of such situations.
 For the case $\phi\in  L^{\infty}(\r,X)$  see [3-4], [16], [9, Sections 4-5], [10, Sections 4-6], [67, Sec. 4-5], [5-6],
  [17], [34] and references therein.

\head{\bf\S 1.  Notation,  Definitions and Preliminaries}\endhead

 In the following $\jj$ will always  be an interval  of the form  $\r,  (\al, \infty),
 [\al, \infty)$ for some  $\alpha \in \r$,
$\r_+ =[0, \infty)$, $\r^+=(0, \infty)$, $\N= \{1, 2,\cdots \}$
and $\N_0=\{0 \} \cup \N$.
 Denote by $X$ a real or  complex Banach space, with scalar field $\k=\k (X)$, $\k=\r$ or $\cc$,   and by $L(X)$ the Banach
space of all linear bounded operators from $X$ to $X$.

 If $f$ is a function defined on $\jj\to X$,
then   $f_s$, $\Delta_s f$ will stand for the functions defined on $\jj$  by $f_{s}(t) = f(t+s)$,
 $\Delta_s f(t) = f_s (t) - f (t)$ for all $s\in \r $ with $s+\jj \st
 \jj$,
  $|f|$ will denote the function $|f|(t):= ||f(t)||$  for all $t\in \jj$ and $||f|| _{\infty} := \text {  sup}_{x\in \jj} ||f(x)||$.

For $U, V \st X$ or $\st X^{\jj}$, $U-V:= \{u-v :\, u\in U, \,
v\in V\}$; similarly for $U+V$.

 If $f \in L_{loc}^1 (\jj, X)$, then $Pf$ will denote the
indefinite integral defined by $Pf(t) = \int_{\alpha_0}^{t}
f(s)\,ds$ (where   $\al_0 =\al +1$  respectively  $0$ if $\jj=\r$,
all integrals are Lebesgue-Bochner integrals (see [53, p. 79],
[39], [50]), similarly for measurable).

The  space of functions with  continuous   derivatives of order up
to $k$ defined on $\jj$ with values in  a Banach space $X$ will be
denoted by $\, C^k (\jj, X)$,  $\, C(\jj, X):=\, C^0 (\jj, X)$.
The spaces of all  constants,  bounded continuous, continuous with
relatively compact range, continuous  with relatively weakly
compact range, uniformly continuous, bounded uniformly continuous,
uniformly continuous with relatively compact  range, uniformly
continuous with relatively weakly compact range, uniformly
continuous and vanishing at infinity, and continuous with compact
support in $\jj$ functions will respectively be denoted by $\, X$,
 $\, C_{b}(\jj,X)$, $\, C_{rc}(\jj,X)$, $ C_{wrc}(\jj,X)$, $C_{u}(\jj,X)$,
 $ C_{ub}(\jj,X)$, $ C_{urc}(\jj,X)$, $ C_{uwrc}(\jj,X)$, $C_0 (\jj, X)$ and $ C_c(\jj,X)$.

 $C^{\infty}(\r,X)$ will stand for the space of infinitely differentiable
 functions $f :\r\to X$.

\smallskip

   The character  $\gamma_{\om}$ is defined by $\gamma_{\om}(t):= e^{i\om\,t}$, $\om,\,t\in\r$. The Fourier transform of $f\in L^1(\r,X)$ with complex $X$ is
$\hat {f}(\om)= \int_{-\infty}^{\infty} \gamma_{\om} (-t)
f(t)\,dt$ and the reflection of $f$ is defined   by $\check
{f}(t)=f(-t)$.

The Sobolev spaces $W^{1,n}_{loc} (\jj,X)$ are defined below in
(2.5).

 $\h (\jj,\k)$ denotes  the Schwartz test functions (infinitely differentiable $\k$-valued functions
with compact support in $\jj$) [71, pp. 21, 24], [41, pp.
299-302].

 $\h'(\jj,X)$
denotes the set of linear continuous $T: \h(\jj,\k) \to X$ as in
[71, pp. 24, 30] or [70, p. 49].

Here $\jj$ in  $\h (\jj,\k)$, $\h' (\jj,X)$ is always open.

 Similarly, $\f(\r,\k)$ will stand
for the  Schwartz space of all rapidly  decreasing infinitely
differentiable  $\k$-valued functions defined on $\r$ (see [75, p.
146]) and

  $\f'(\r,X)$ is the space of Banach valued tempered
distributions of linear continuous $T: \f(\r,\k) \to X$ (see [71,
p. 234]).

 $ (\h_{L^p})' (\r,X):=\{ T:  \h_{L^p} (\r,\k)\to X:  T\text { linear,\,\, continuous}\}$, $1\le p\le
 \infty$  where $\h_{L^p} (\r,\k)$ contains all $\va  \in  C^{\infty}
(\r,\k)$ with $\va^{(j)}  \in  L^p (\r,\k)$ for $ 0\le j< \infty$,
and $\va_m \to 0$ in $ \h_{L^p}$ meaning
                          $\va_m^{(j)} \to 0$  in the $L^p(\r,\k)$-norm for
$ 0\le j< \infty$ see [71 , p. 199]. The topology on $(\h_{L^p})'$
is given  by the  seminorms
 $|| T||_{V}:= sup \{ ||T(\va)||: \va\in V\}$, $ V$  bounded set in
$\h_{L^p}(\r,\k)$; $V$ bounded means here  sup $\{|| \varphi
^{(j)} ||_{L^p}: \varphi \in V\} < \infty$ for each $j\in\N_0$.

 $ T\in\h' (\r,X)$, $T\in  (\h_{L^1})' (\r,X)$ means there is a (unique)
$\widetilde {T}\in  (\h_{L^1})' (\r,X)$    with
 $\widetilde {T}\, | \h (\r,\k) =T$.

 Translates $T_a$ for distributions $T$ are defined in accordance with the above definition of translates $f_a$ for functions  $f$
by $T_a(\va):= T(\va_{-a})$ (contrary to the definition in [71,
(II, 5; 2), p. 55]).

\smallskip

\smallskip

Let  $\A\subset L^1_{loc}(\jj, X)$ or $\A\subset \h'(\jj, X)$.   We use  the following assumptions  for $\A$.

\smallskip

\noindent $Cone$ :     $s F+ tG\in \A $ if $ F, G \in \A $  and $ s, t \in  [0, \infty)$.

\smallskip

\noindent  $\Q$-$convex$ : $-\sum _{k=1}^{m}  F_k/m \in \A$ for
all $F_k \in \A$, $ m\in \N$.

\smallskip

\noindent $Real$-$linear$:   $s F+ tG\in \A $ if $ F, G \in \A $  and $ s, t \in  \r$.

 \smallskip

\noindent $Positive$-$invariant$: translate $ F_a\in \A $ if $ F\in \A $  and $  a \in  [0, \infty)$.

\smallskip

\noindent $C_{ub}$-$invariant$: For   $ f\in C_{ub}(\r,X)$  with $ f|\,\jj  \in \A$,  $ f_a|\,\jj \in \A$ for all $a\in\r$.

\smallskip

\noindent $Invariant$:  $F_a \in \A$ if $F\in\A$ for all $a\in\r$.

\smallskip

\smallskip

\noindent $Uniformly\,\, closed$  :
 $(\phi_n) \st  \A  $  and  $\phi_n \to  \phi $ uniformly on $\jj$
implies  $\phi \in  \A $.

\smallskip

\noindent $C^{\infty}$-$uniformly\,\, closed$  :  Uniformly closed, but only for
 $(\phi_n) \st  \A \cap C^{\infty}(\jj,X)$.

\smallskip

\noindent $(\Gamma)$:   $ \,\,\,\,\gamma_{\om}\phi \in\A$ for all
characters $\gamma_{\om} (t)= e^{i\om \, t}$, $\om\in\r$ and all
$\phi\in \A$.

\smallskip

\noindent $(L_ U)$:  For  $\A, U  \subset L^1_{loc}(\jj, X)$ :
$\phi \in U$, $\Delta_h \phi\in    \A $ for all $h > 0$ implies
$\phi \in \A$.

\smallskip

\noindent $(L'_ U)$:  For  $\A, U  \subset \h'(\jj, X)$ :  $T \in
U$, $\Delta_h T\in    \A $ for all $h > 0$ implies $T \in \A$.

\smallskip

 For  $U= C_b (\jj, X)$, $C_u (\jj, X)$, $C_{ub} (\jj, X)$, $C_{uwrc} (\jj, X)$ we write

\smallskip

 $ (L_b)$, $(L_u)$, $(L_{ub})$,  $(L_{uwrc})$.

\smallskip

  There in  $(L'_U)$     the ``$\in$'' for $T=f\in  L^1_{loc}(\jj, X)$,  $U\subset L^1_{loc}(\jj, X)$ means
    for example there is $g\in\A$
with   $f =g $ almost  everywhere  on $\jj$; in $(L_U)$ however
$\phi =g$ everywhere on $\jj$. (See special case  to Lemma 4.4
below).

\smallskip

For $\A$ with $(L_ U)$ respectively $(L'_ U)$ we say that $\A$ is
a $U$-Loomis class (see Loomis [60, p. 365], Caracosta-Doss [31],
[13, pp. 117, 120]).

\smallskip

\smallskip

\noindent A function $\phi\in L^1_{loc}(\jj, X)$ is called
$ergodic$ if there is a constant $m (\phi) \in X$ such that

\smallskip

 sup $ _{x\in \jj} || \frac{1}{T}\int_{0}^{T}\phi (x+s)\,ds - m (\phi)||\to 0$  as $T\to \infty$.

\smallskip

\noindent The limit $m(\phi)$  (clearly unique) is called the mean
of $\phi$.

\smallskip

\noindent  A function $\phi\in L^1_{loc}(\jj, X)$ is called
$totally \,$ $ergodic $ if $\gamma_{\om} \, \phi $  is ergodic for
all characters $\gamma_{\om} (t)= e^{i\om \, t}$, $\om\in\r$.

\smallskip

\smallskip

  In this paper we adopt the following notation (note the difference from [9], [10], [13] (there $\E, \T\E$ stands for $\E_{ub}, \T\E_{ub}$ defined below)):

 The space of all ergodic (totally ergodic) functions from  $L^1 _{loc}(\jj, X)$ will be denoted
by

$\E(\jj,X),\qquad$ ($\T\E(\jj,X)$, then $\k =\cc$ ).

\smallskip

 We set

\smallskip

$ \E_{0} (\jj,X):= \{\phi\in \E (\jj,X): m(\phi)=0 \}$,
   $\,\,\, \E_{ub} (\jj,X):= {\E} (\jj,X) \cap C_ {ub} (\jj,X)$,

  $ \E_{0,ub}
   (\jj,X)= \E_{0} (\jj,X) \cap C_{ub} (\jj,X) $,

$ \T\E_{0} (\jj,X):= \{\phi\in \T\E (\jj,X):  \gamma_{\om}\phi \in
\E_0 (\jj,X),\,\, \om\in\r\}$,

    $\T \E_{ub} (\jj,X):= \T {\E} (\jj,X)\cap  C_ {ub} (\jj,X)$,
$\,\,\T \E_{0,ub} (\jj,X)=  \T\E_{0} (\jj,X) \cap C_{ub} (\jj,X)
$.

\smallskip

\noindent If $Y$ is a locally convex  complex space, let $ \Pi
(\r,Y)$ stand for the set of all trigonometric polynomials

$\pi (t):= \sum _{j=1}^{m} a_j e^{i\, \om_{j}\, t}$, where  $a_j
\in Y$, $\om_{j}\in \r$, $j=1, \cdots, m \in \N$,

\smallskip

$P_{\tau} (\jj,Y): = \{\phi \in C(\jj,Y): \phi (t+\tau)= \phi (t)$
for all $ t\in \jj\}$, $\,\,\,\tau$ fixed  $>0$.

\smallskip

 \noindent  A function  $\phi : \r\to Y$
is called $almost\,$ $\, periodic$ (ap) if and only if to each
neighborhood $V$ of $0$ of $Y$ there exists  $\psi \in \Pi (\r,Y)$
with  $\psi-\phi \in V$ for all $t\in \r$;

\noindent for $Y=$ Banach space $X$ this means the existence of
  a sequence $(\pi_n) \st \Pi (\r,X)$ such that $||\pi_n
-\phi||_{\infty} \to 0$ as $n\to \infty$.

\smallskip

\noindent The space of all almost periodic functions will be
denoted by $AP(\r, Y)$, always

$AP(\r, Y) \st C_{ub}(\r, Y)$.
\smallskip

\noindent
 Since $\Pi (\r,X) \st \T \E_{ub} (\r,X)$, one gets also $AP (\r,X) \st \T\E_{ub}
 (\r,X)$, $X$ Banach space. See [2, p. 21], [33, Lemma 1].

 \noindent  A function  $\phi \in C_b
(\r,X)$ is  $almost\,$ $\, periodic$ if and only if  for each $\e
>0$ the set  of $\e$-periods
$\,\,\T (\phi, \e ):= \{ \tau \in \r: ||\phi (t+\tau)-\phi
(t)||\le \e$ for all $t\in\r\}$

\noindent is relatively dense in $\r$ (Bohr's Definition);

\noindent here $M\st \r$ relatively dense means the existence of
an $l > 0$ such that

$M\cap [t, t+l]\not = \emptyset$ for each $t\in \r$.

\smallskip

\smallskip

\noindent $ \phi \in C_b (\r,X)$ is almost periodic  if and only
if $\{\phi_s : s\in \r\}$ is relatively compact in $C_b (\r,X)$
(Bochner's definition). See [2, p. 7], Proposition 6.2,  Lemma 6.4
and $(vii)$ after Lemma 6.4.

\smallskip

\noindent A function $\phi \in C (\r,X)$ is called   Bochner
almost automorphic (aa)  if to each sequence $(s_m)$ from $\r$
there exists a subsequence $(s_{m_n})$ such that for each   $t
\in \r $
 the $\lim_{n\to \infty} \phi(t+s_{m_n})$ exists, $= : g(t)$, and furthermore $g(t -
s_{m_n})\to \phi(t)$ as  $n \to \infty$ for each  $ t \in\r$.
 (See [24, p 2041, Definition 2], [46, pp. V, 11] and the references there).

\smallskip

$UAA (\r,X)$ will denote the set of all uniformly continuous
Bochner aa $\phi : \r \to X$. For characterizations via
$\e$-periods "for $|t| \le M$ " see [13, p. 119] and the
references there.

\smallskip

$LAP_{ub}(\r,X)$ will denote the set of Levitan ap  $=\, N$-ap
functions which are in $C_{ub}(\r,X)$; for various equivalent
definitions with "local $\e$-periods" see again [13, p. 119] and
the references there, also [59, section 4], [9, Definition 2.1.3
(iii), Theorem 2.1.4].

\smallskip

 \noindent  Let $\phi \in L^p_{loc}(\jj,X)$, $1\le p < \infty$. Denote  by

$||\phi ||_{S^p_l} := $ sup$_{x\in\jj}\, [\frac{1}{l}\,
\int_{x}^{x+l} ||\phi(t)||^p\, dt]^{1/p}$, $\qquad ||\phi ||_{S^p}
:= ||\phi ||_{S^p_1}$,

\smallskip

$||\phi ||_{W^p} :=\lim _{l\to \infty}\, ||\phi ||_{S^p_l} $,

\smallskip

$||\phi ||_{B^p} :=\overline{\lim} _{T\to \infty}\,[\frac{1}{2T}\,
\int_{-T}^{T} ||\phi(t)||^p\, dt]^{1/p}$.

\noindent The limit $\lim _{l\to \infty}\, ||\phi ||_{S^p_l} $
always exists (see [21, pp. 72-73], also valid for $X$-valued
$\phi$).
\smallskip

$S^p_b (\jj,X) := \qquad \{\phi \in L^p_{loc} (\jj,X): ||\phi
||_{S^p <\infty}\}$, $1\le p< \infty$.
\smallskip

 \noindent  A function  $\phi \in L^p_{loc} (\jj,X)$ is called
respectively

$Stepanoff$ $S^p$-$almost$ $periodic$,

  $ Weyl$ $W^p$-$almost$
$periodic$,

 $Besicovitch$ $B^p$-$almost$ $periodic$

\noindent if and only if there exists a sequence $(\pi_n) \st \Pi
(\jj,X)$ such that respectively

$||\pi_n -\phi||_{S^p} \to 0$ as $n\to \infty$,

 $||\pi_n- \phi||_{W^p} \to 0$ as $n\to \infty$,

$||\pi_n- \phi||_{B^p} \to 0$ as $n\to \infty$.

\smallskip

 \noindent   Stepanoff $S^p$-almost periodic  functions, Weyl $W^p$-almost
periodic  functions,
 Besicovitch $B^p$-almost periodic  functions
 for all $1\le p <\infty$  will be denoted respectively by

 $S^p AP(\jj,X)$,  $\qquad W^p AP(\jj,X)$,  $\qquad B^p AP(\jj,X)$

 \noindent ( see [2, p. 76],  [21, pp. 71-78], [27, pp. 34-36]).

\smallskip

\smallskip

\noindent  Here for $\A= AP$,  $S^pAP$,  $W^pAP$,  $B^pAP$, $UAA$,
$LAP_{ub}$ for simplicity reasons we define for  $\jj\not= \r $

\smallskip

(1.1) $\,\,\, \qquad\,  \A(\jj,X)\, =\, \A(\r,X) |\,\jj$.

\smallskip

\noindent  For most of these spaces one can give definitions for
$\jj\not = \r$ such that (1.1) becomes a theorem; for $AP(\jj,X)$
see [47, p. 96, footnote 11, Satz], but  also  Example 5.5.

\smallskip

 \noindent
 A function $\phi\in C_b (\jj, X)$ is called
$asymptotically$  $almost$  $periodic$  (respectively $Eberlein$
$almost$ $periodic$)  if  $\jj$ is a semigroup,  $H(\phi)
=\{\phi_s : s\in \jj\}$ is relatively compact   (respectively
weakly relatively compact) in $C_b (\jj, X)$; here for
asymptotically \apf  only $\jj\not =\r$ is allowed. A function
$\phi\in C_b (\r, X)$ is called $asymptotically$  $almost$
$periodic$  if $f=g+\xi$  with $g\in AP(\r,X)$, $\xi\in C_0
(\r,X)$ (see [43, footnote (1) p. 521], [68, pp. 14-15]).
\smallskip

$AAP(\jj, X)$, $ EAP(\jj, X)$, $ EAP_{rc}(\jj, X) $

\noindent will respectively  stand for
 the spaces of
asymptotically almost periodic functions,  Eberlein  almost
periodic functions, $EAP$-functions with relatively compact
 range.

\smallskip

By the  Eberlein-$\check{S}$mulian  theorem [39, p. 430],  for
$\phi\in C_b (\jj,X)$, $\phi$ Eberlein ap  is equivalent with: To
each sequence $(t_n)\subset\jj $ there is a subsequence
$(t_{n_m})_{m\in \N} $ and $\psi\in C_b (\jj,X)$ such that
$\phi_{t_{n_m}}\to \psi$ weakly in $ C_b (\jj,X)$.

\noindent For $EAP$ one has

 $ EAP(\r,X)\, |\,\jj \st\, EAP(\jj,X)$.

\noindent For closed $\jj$ see \S 11 question 19.

\smallskip

\noindent $EAP (\jj,X)$ is defined only for semigroups $\jj$, that
is $\jj=\r$, $(\al,\infty)$ or $[\al,\infty)$ with $\alpha \ge 0$
([20, p. 138], $X=\cc$, [45], [62], [68, p. 15], [69, p. 424], [9,
pp. 12-14]).
\smallskip

\noindent  Furthermore for example

\smallskip

(1.2a)\qquad   $EAP (\overline {\jj}, X) |\,\jj$   strictly  $\st
EAP (\jj,X)$ for $\jj= (\al,\infty)$ with $\al >0 $:

\smallskip

\noindent  There exists  $\phi \in C_b ((\al,\infty),X)$ vanishing
on $[2\al,\infty)$ which cannot be extended continuously to
$[\al,\infty)$; such a    $\phi \in EAP( (\al, \infty),X)$. Then
$\phi$ gives the ``strictly''.

\noindent $\st$ follows from the definitions with $y(f) : =
z(f|(\al,\infty))$ for
       $z  \in  [C_b((\al,\infty),X)]^*$  (see also [20, p. 143
       Theorem 4.2.10] for $X = \cc$).
\smallskip

 \noindent Whereas one can show

\smallskip

(1.2b)\qquad $EAP (\overline {\jj}, X) |\,\jj = EAP (\jj,X)$ if
$\jj=
 (0,\infty)$:

\smallskip

\noindent Let $\phi \in EAP ((0,\infty),X)$ and $(t_n) \st
(0,\infty)$ with $t_n \to 0$ as $n\to \infty$. Then $(\phi
_{t_n})\st C_b ([0,\infty),X)$ and $\phi _{t_{n_k}} \to \phi $
weakly in $C_b ((0,\infty),X)$. By Mazur's theorem [75, p. 130] a
sequence  $(\psi_m)$ of convex linear combinations of $(\phi
_{t_{n_k}})$ converges to $\phi$ in $C_b ((0,\infty),X)$. This
implies $\lim_{m\to \infty} \psi_m \in C_b ([0,\infty),X)$ showing
that $\phi$ can be extended  by continuity  to $[0,\infty)$. By
the Hahn-Banach and Mazur theorems this extension $\in EAP
([0,\infty),X)$.

\smallskip

\noindent Also, one has, for any $X$ and $\,\,0\le \al <\beta<
\infty$

\smallskip

 (1.2c)\,\,\,\,\,\, $EAP ( [\beta, \infty), X) = EAP ([\al, \infty),X)| [\beta,
\infty) = EAP ((\al, \infty),X)| [\beta, \infty) $:

\smallskip

$ A_1 : = EAP([\al,\infty),X)|[\beta,\infty)  \st
EAP((\al,\infty),X)|[\beta,\infty)= : A_2$:

   \noindent  This follows with " $\st$ " from (1.2a).

   $ A_2  \st  EAP([\beta,\infty),X) = : A_3$ :

\noindent Follows with $y(\phi) :
           = z(\phi|[\beta,\infty))$ for $z  \in  [C_b([\beta,\infty),X)]^*$.

  $ A_3  \st  A_1$ :

\noindent     If   $f \in EAP ( [\beta, \infty), X)$, then
$f_{2\beta}  \in  EAP(\r_+,X)$ : If $(r_m) \st \r_+$,  $t_m : =
\beta + r_m  \in  [\beta, \infty)$, so there is a subsequence
$(t_{m_n})$ and a $g  \in  W : = C_b([\beta,\infty),X)$ with
$f_{{\beta}+r_{m_n}}  \to g$ weakly in $W$. If now $y \in U^*$, $U
: = C_b (\r_+ , X)$, $z(w) : = y(w_{\beta})$ defines a $z \in
W^*$, so $y((f_{2\beta})_{r_{m_n}} ) = y((f_{\beta})_{\beta+
r_{m_n}}) = z(f_{\beta+r_{m_n}})\to z(g)$. Since $f_{2\beta}$,
$f_{\beta} \in U$;  with $ G : = g_{-\beta}$ on $[2\beta,\infty)$,
$G : = g(\beta)$ on $[0,2\beta)$ one has $G \in U$ and $G_{\beta}
= g$ on $[\beta,\infty)$, so $y( (f_{2\beta})_{r_{m_n}}) \to y
(G)$, showing $f_{2\beta} \in EAP(\r_+,X)$. Then $f_{2\beta} \in
C_{ub}(\r_+,X)$ by [69, Proposition 2.1] or [9, Theorem 2.3.4];
this implies $F  \in C_{ub}([\al, \infty),X)$ if $F : = f$ on
$[\beta, \infty)$, $: = f(\beta)$ on $[\al,\beta)$. If now $(s_m)
\st [\al,\infty)$ has a subsequence $s_{m_n} \to s_0  \in
[\al,\infty)$, $F_{s_{m_n}} \to F_{s_0}$ even uniformly on
$[\al,\infty)$. Else  $s_m \to \infty$, so there is a subsequence
$(t_n)$ with $t_n >2\beta$ for $n  \in \N$. As above, to $y  \in
V^*$, $V : = C_b([\al,\infty),X)$, $z(w) : = y(w_{\beta - \al})$
defines a $z \in W^*$, one gets $y(F_{t_n}) = y(f_{t_n})  = y(
(f_{t_n + \al - \beta})_{\beta - \al}) = z( f_{t_n + \al -
\beta})$ . Since $t_n + \al - \beta \ge ß\beta$, there are a
subsequence $(t'_n)$ of $(t_n)$  and a  $g \in  W$, both
independent of $y$, with $y(F_{t'_n}) \to z(g) = y(g_{\beta -
\al})$ as above with $g_{\beta - \al} \in  V$, i.e. $F \in
EAP([\al,\infty),X)$.

\smallskip

 \noindent (1.2c), (1.2b) and [9, Theorem 2.3.4] give for any $X$

 (1.2d)\qquad    $EAP(\jj,X)  \st  C_{ub}(\jj,X)$ for closed $\jj$ or $\jj = (0,\infty)$.

\smallskip

\smallskip

\noindent  For any admissible  $\jj$ and arbitrarily $X$ one has

\smallskip

$\qquad AP(\jj,X)  \st UAA(\jj,X)   \st  LAP_{ub}(\jj,X)\qquad$
((1.1), [13, p. 119]).

\smallskip

(1.2e)\qquad $EAP (\jj,X) \cap  LAP_{ub}(\jj,X) = AP(\jj,X)$, for
$\jj= \r$ or $\r_+$, any $X$:

\smallskip

\noindent This follows with $\phi \in EAP (\jj,X)$ implies $\phi =
\psi +\xi$, where $\psi \in AP(\jj,X)$, $\xi \in EAP_0 (\jj,X)$
 ([68, p. 18] for $\jj=\r$, [69, Theorem 2.4] for $\jj=\r_+$) and
$EAP_0 (\jj,X)\cap LAP_{ub}(\jj,X)=\{0\}$ (see [12?, p. 1142,
Proposition 2.4]).

\smallskip

\noindent Also, one has (see [62, Example 1])

 $ AP(\jj,X)  \st EAP_{rc}(\jj,X) \st  EAP(\jj,X)$,
any $X$, admissible $\jj$,

\smallskip


\smallskip

(1.2f)$\qquad EAP(\jj,X)\st\T\E(\jj,X)\,\,\,$   for any $X$,
admissible $\jj$:

\smallskip

\noindent If $\jj = (\al,\infty)$ or $[\al,\infty)$ and $f  \in
EAP$, then, with $\jj' : = [\al + 1, \infty)$,   $f|\jj'  \in
EAP(\jj',X)$ by (1.2c), and, again with (1.2c),  $f|\jj'  =
F|\jj'$ with $F \in EAP(\r_+,X)$,  $ \st \E(\r_+,X)$ by  [69, p.
425, Theorem 2.4 and Proposition 2.3]. With $f \in  C_b(\jj,X)$
one gets $f\in \E(\jj,X)$; since $EAP(\jj,X)$ satisfies $(\Gamma)$
almost by definition, $f \in \T\E(\jj,X)$ follows. If $\jj = \r$,
$f|\r_+  \in EAP(\r_+,X)$ as in (1.2c), $\st \E(\r_+,X)$; since
then also $\check {f} \in
 EAP(\r,X)$, $\check {f} | \r_+  \in \E(\r_+,X)$; together $f  \in \E(\r,X)$, then $f
\in \T\E(\r,X)$.

\smallskip

\noindent With (1.2d) one gets

(1.2$f'$)\qquad $ EAP(\jj,X) \st \T\E_{ub}(\jj,X)$ for $X$, $\jj$
as in (1.2d).

\bigskip

\noindent   If $V(\jj,X)  \st L^1_{loc}(\jj,X)$, we  define

\smallskip

$\,\,        VAP(\jj,X)  : =$

\smallskip

              $\,\,   V (\jj,X) +  AP(\jj,X) : =$
$ \{ v+f : \,\, v \in V (\jj,X),\,\, f \in AP(\jj,X) \}$;

\smallskip

\noindent $ V = C_0(\jj,X)$ gives $ AAP (\jj,X)$,

\smallskip

\noindent $V = EAP_0  (\jj,X)=$ $EAP$-null functions $:=$

 $ \{f \in  EAP(\jj,X) :  0 $ in  weak closure of
$\{f_h : h  \in  \jj\}$ in $C_b(\jj,X) \}$

\smallskip

\noindent gives $EAP (\jj,X)$
                 (see [69, p. 424 ], [68, p. 34], [64, Theorem]).

\noindent Define

$Av_0 (\jj,X) : = \{f \in  L^1_{loc}(\jj,X):  lim_{T \to \infty}
(1/T) \int^{\al_0 +T} _{\alpha_0}  f(t)\, dt\,\,$ exists, $= 0
\}$,

 $ Av_n (\jj,X):= \{f  \in  L^1_{loc} (\jj,X) : lim_{T
\to \infty} (1/T) \int^{\al_0 +T} _{\alpha_0}  ||f(t)||\, dt$
exists, $= 0 \}$.

 \noindent Then $V = Av_n (\jj,X) \cap C_b (\jj,X)$ gives ($\jj \not= \r$) Zhang's
  pseudo-almost periodic functions  $PAP_0 (\jj,X) $ ([78, p. 168], [77], [14]).

\smallskip

\noindent  Other examples would be (ergodic null functions)

$V =  \E_n (\jj,X): = \{f  \in \E(\jj,X) : |f|  \in \E_0
(\jj,\r)\}$,

  $V
= \E_{n,ub}(\jj,X): = \{f \in \E_{ub} (\jj,X) : |f| \in \E_0
(\jj,\r)\}$, $= \E_{ub}(\jj,X) \cap \E_n (\jj,X)$.

\bigskip

\noindent If an $\A =  \A(\jj,\k)  \subset L^1_{loc}(\jj,\k)$ is
given, the corresponding $weak $ class is defined  by (see [2, p.
38], $\A= AP(\r,\k)$)

\smallskip

(1.3) $\qquad \,\, \W\A (\jj,X) := \{ f\in L^1_{loc}(\jj,X):\,\,\,
y\circ  f \in\A (\jj,\k) \text { \,\,\, for\,\,all \,}  y\in X^*
\}$.

\smallskip

\noindent $All\,\, the\,\, above\,\, \A\,\, except\,\,\, C_0
(\jj,X),\, \,\,C_c (\jj,X)\,\,\, are\,\, linear,\,\,
\,positive\,\,\, invariant\,\,\, and$

\noindent $uniformly
 \,closed,\,$ for $VAP$ see Proposition 7.13.

 This follows mostly from the definitions; uniformly closed for
$C_{rc}$        is clear with relatively compact $=$ totally
bounded in $X$, for $C_{wrc}$
        follows as in [13, p. 119 below, for $C_{uwc}$ there, continuity is not necessary].

\smallskip

In the sequel if there is no danger of confusion we omit $\jj, X $
in referring to any of the above classes; for example we write
$AP$ instead of  $AP(\jj,X)$.
\smallskip

 \proclaim{Definition 1.1} Let $\phi\in L^1_{loc}(\jj,X)$, $ \A \subset L^1_{loc}(\jj,X)$  or   $\st \h'(\jj,X)$ and $ 0\not =h \in \r$.
Define

\smallskip

(1.4) $\,\,\,\qquad M_h\phi (t) := \frac{1}{h}\int_{0}^{h}
\phi(t+s)\, ds,\qquad t\in \jj$.

\smallskip

(1.5) $\,\,\,\qquad \m\A := \{\phi\in L^1_{loc}(\jj,X):\qquad
M_h\phi \in\A \text { \qquad for\,\,all \,}h\,>0 \}$.

\smallskip

\noindent  Recursively   if $k \in \N$

\smallskip

 (1.6) $ \,\,\,\qquad \m^k\A=\m (\m^{k-1}\A)$, $\, \m^1\A := \m\A $,   $\, \m ^0 \A := \A \cap L^1_{loc}  = : \A_{Loc}$.

\endproclaim

 Similarly, we  define these  means  for  distributions:

 \proclaim{Definition 1.2} Let   $T\in \h'(\jj,X)$ and  $\A\st \h'(\jj,X)$.
Define

\smallskip

(1.7) $\,\,\,\qquad\widetilde { M}_h T (\varphi):=  T(M_{-h}
\va)$, $ \, \varphi \in   \h (\jj ,\k)$  with $\varphi:= 0$   on $
\r \setminus\jj$    if $\jj \not=\r$.

\smallskip

 Set

\smallskip

(1.8) $\,\,\,\qquad\widetilde { \m} \A=\{T\in
\h'(\jj,X):\widetilde { M}_h T\in \A  \text { \qquad for\,\,all
\,}h\,>0 \}$.

\smallskip

\noindent Recursively  if $k\in \N$

\smallskip

(1.9) $\,\,\,\qquad\widetilde{ \m}^k\A :=\widetilde{
\m}(\widetilde{\m}^{k-1}\A) \,$,  $\,\widetilde { \m}^{1}\A
=\widetilde { \m}\A\,$ and $\,\widetilde { \m}^{0}\A =\A$.

\smallskip

\endproclaim

\noindent (1.7)    indeed  defines a distribution from
$\h'(\jj,X)$ if   $ T\in \h'(\jj,X)$. The operator $\widetilde{
M}_h $ is linear, continuous   and coincides with  $M_h$  of
Definition 1.1  on $ L^1_{loc}(\jj,X)$ since

\smallskip

  $ \,\,\,\qquad\, \int_{\jj} M_h (\phi)(s)\va (s) ds= \int_{\jj} \phi(s)(M_{-h} \va)(s) ds $

\smallskip

\noindent for all $\phi\in L^1_{loc}(\jj,X)$ and  $\va \in
\h(\jj,\k)$ with $\varphi$  extended by $0$ on $ \r \setminus\jj$
if $\jj \not=\r$.

\noindent For $ T\in \h'(\jj,X)$    and $a > 0$ also
 $ T_a\in \h'(\jj,X)$,  and $\widetilde {M}_h$
                          commutes with translation,

  $\,\,\,\qquad\widetilde {M}_h (T_a) = (\widetilde {M}_h T)_a$  if $h >0$ (respectively $a\in \r$ if $\jj=\r$).

 \proclaim{Definition 1.3}    For  $\A\st \h'(\r,X)$
define

\smallskip

(1.10) $\,\,\,\qquad\h'_{ \A} (\r,X)=\{T\in   \h'(\r,X):
T*\varphi\in \A  \text { \,\,\, for\,\,all \,}\varphi \in
\h(\r,\k) \}$.

\endproclaim

 \noindent Here  $(T*\va)(x):= T((\check {\va})_{-x})\,\,$ for $x\in \r$

 \noindent (see [70, (I, 3; 12), p. 72] or [75,  p. 156, (2)]).
\smallskip

 We use  $\, s_{h}:= (1/h)\chi_{(-h,0)}$, where $\chi_{(-h,0)}$ is the characteristic function of the interval $(-h, 0)$, $h> 0$; $\, s_{h}:= (1/(-h))\chi_{(0,-h)}$ if $h< 0$.

\noindent   For  $ \phi\in L^1_{loc}(\jj,X)$ the convolution $\phi*s_h$ is defined on  $\jj$, $\in C(\jj, X)$   and in this sense one has

\smallskip

 (1.11) $\,\,\,\qquad    M_{h_1}\cdots  M_{h_k}\phi= (\phi * s_{h_k})*\cdots *s_{h_1}$ on $\jj$,  for all $h_1>0, \cdots, h_k>0$.

\smallskip

For closed $\jj$ and $\tilde {\phi} =\phi$ on $\jj$ and $:=0$ on
$\r\setminus \jj$, and substituting     $\tilde {\phi}$ for $\phi$
in (1.11), we can use there the usual convolution, defined and
continuous on $\r$.
\smallskip
 \noindent Moreover, if  $\A \st { \m}\A$ respectively  $\A \st\widetilde{ \m}\A$ (see Proposition 3.2 ), then

\smallskip

 (1.12) $\,\,\,\qquad { \m}^{k-1} \A \st { \m}^k\A \,$  respectively
 $\, \,\widetilde{ \m}^{k-1} \A \st\widetilde{ \m}^k\A \,$ for all $ \, k\in\N$.

\smallskip

For any      $\A \st  L^1_{loc}(\jj,X)$ or $ \h' (\jj,X)$ and any $n\in \N$ one has, almost by definition

\smallskip

 (1.13) $\,\,\,\qquad$ If $ f\in  \m^{n} \A $,  $g = f $  almost everywhere on $\jj$, then   $ g\in  \m^{n} \A $.

\smallskip

Here and later we need the following Doss-condition (see
Caracosta-Doss [31], Doss [37] and Lemma 2.3, Propositions 2.4,
 4.2, 5.1, 8.1,  Corollary 2.15 below ):

\proclaim {Definition 1.4} We say that $\A \st  L^1_{loc}(\jj,X)$ or $\st  \h' (\jj,X)$  satisfies $(\Delta)$  if for any  $\phi\in L^ 1 _{loc}(\jj,X)$ for which  the differences
 $\Delta_s \phi \in \A$, $ 0< s \in\r$,  one has   $(\phi-  {M_h }\phi)\in \A$ for all $h >0$;

\smallskip

\noindent  $\A$ satisfies $(\Delta_1)$
 if the conclusion holds for $h=1$.

\smallskip

\noindent  $\A$  satisfies $(\Delta')$ respectively
$(\Delta'_1)$  if for any  $T \in  \h'(\jj,X)$ for which  the
differences
 $\Delta_s T\in \A$  for all $s > 0$ one has   $(T- \widetilde {M_h }T)\in \A$  for all $0< h\in \r$ respectively $h=1$
 (here both ``$\in$" are meant in the distribution sense).

\endproclaim

\noindent $(\Delta')$ for $\A\st C(\jj,X)$ implies $(\Delta)$;

for $\A\st L^1_{loc}(\jj,X)$  one only gets  $(\phi- \widetilde
{M_h }\phi) \,  ``\in \A$ a.e''.

\noindent For a converse see Proposition 4.3.

\proclaim {Proposition 1.5}  Let $I\st \r$ be an  arbitrary interval, $\e_0 > 0$, $k\in \N_0 $.
If then   for $\phi\in L^1_{loc} (I, X)$ the difference
 $\Delta_h\phi \in C^k (I^{-h},X)$  for  all $h\in (0, \e_0]$, then  $\phi \in C^k (I, X)$ (see Remark 1.6 (i)).
\endproclaim

\noindent Here   $ I^{-h}= (\alpha, \beta - h)$ respectively  $(\alpha, \beta - h]$  if $I =(\alpha, \beta )$ respectively   $(\alpha, \beta ]$,
similarly for  $[\alpha, \beta )$, $[\alpha, \beta )$, $-\infty \le  \alpha < \beta   \le \infty$.

\demo {Proof} Case $k=0$: For  $X=\k$ one gets continuity of $\phi$ on     $(\alpha, \beta )$ by
  [49, Lemma 13, p. 224 ]  which reads
          : Assume  $I = (\alpha,\beta)$ open interval  $\st \r$, $k  \in  \N_0$, $1\le  p\le
          \infty$, $f : I \to \k \, $ Lebesgue measurable, $\, \Omega  \st \r$
such that for each $ \e > 0$ the set  $\Omega  \cap (-\e,\e)$ is
Lebesgue  measurable with positive measure. If then the
differences
    $\Delta_h f  \in  \A (I^{-h},\k)$ for each $h  \in  \Omega$,
then $f  \in  \A (I,\k)$. Here $\A$ can be $W^{p,k}_{loc} $ or
$C^k$ (among others), $I^{-h} : = (\alpha + |h|, \beta - |h|)$
 (see also [30, p. 197, theorem 1.3] or [10]).

 If  $\alpha$ or $ \beta \in I$  continuity  there follows from the continuity
 of the $\Delta_h\phi$ in  $\alpha$ respectively $ \beta $.

For general $X$,
             $\phi$ is weakly continuous on $I$ by the above and
            therefore locally norm bounded. Choose a fixed $t_0$  such that
            $I_0 : = [t_0 -\e, t_0 + \e]  \st  I$ with positive
$\e$, then  define   $\psi : \r \to X$ by  $\psi = \phi$ on $I_0$
and $\psi = \phi(t_0 -\e)$
            left of $I_0$ and  $ \psi = \phi(t_0 +\e)$ to the right; then
            $\psi$
is boundded and weakly continuous on $\r$. To $\psi$ one can apply
a result of Gelfand [44, p. 237, Satz 1] which implies: If X is
separable Banach space and $f: \r \to X$ is weakly continuous,
then the set of discontinuity of $f$ is at most of first category
(see [75, p. 12, Baire's theorem 2, valid for $X$-valued $x_n$]),
 one gets thus norm-continuity of $\psi$, $\phi$ at $t_n \to t_0$
and then  at $t_{0}$  with the continuity of $\Delta_{t_0-t_n}
\phi$  at $t_n$, so on the interior of $I$ (see also [10, Theorem
2.1 (b)]). Continuity at eventual end points follows as above.

Case $k=1$: By  [49, Lemma 13]  for each $ x^*\in X^*$,
$(x^*\circ\phi)'$  exists and is continuous first on $I^{\circ}$,
then on $I$. Define $G(t)(x^*)  :=(x^*\circ\phi)'(t) $, $t\in I$,
$ x^*\in X^*$. By the uniform boundedness theorem [39, p. 55],  $
G(t)\in X^{**}$. The definition of $G$ gives      $ \Delta_h G
=(\Delta_h\phi)'  $, $    \in C (I^{-h},X)$ even. We show that $G$
is weakly continuous on $I$. Indeed, take $y\in X^{***}$, then $y
|\, X =: x^* \in X^*$. One has
 $ y (G (t+h) -G(t))= x^*(\Delta_h G(t)) = \Delta_h G (t) (x^*), $ $\to 0$  as $ h\to  0$,
  since  $G$  is w$^*$-continuous by its  definition above.
This gives  weak continuity of $G$ on $I^o$, then on $I$ with
$\Delta_h\phi \in C^1 (I^{-h},X)$. This implies   $G \in
L^{\infty}_{loc} (I, X^{**}) $.  The case $k=0$ gives  $G : I\to
X^{**}$ is  continuous. Now $ h G (t)=    \int_{0}^{h} G (t+s)\,
ds - \int_{0}^{h} \Delta_s G(t) \,ds $ \noindent  (with  the
$X^{**}$-Bochner integral, existing because of continuity),  for
all $0 < h <\e_0$, $t\in I^{-h}$. Since  $\Delta_s G (t) \in X$,
by a theorem for the Bochner integral (see [50, p. 62, \S2,
Aufgabe 29]), one gets that the second integral is in $X =$ closed
linear subspace of $X^{**}$. With Hille's theorem for the Bochner
integral ([75, p. 134, Corollary 2]) one gets
 $\int_{0}^{h} G (t+s)\, ds (x^*)= \int_{0}^{h} G (t+s)(x^*)\, ds=  \int_{0}^{h}  (x ^*\circ\phi)'(s+t)\, ds = \Delta_h (x^*\circ \phi)(t) = x^* (\Delta_h  \phi)(t)$, $x^* \in X^*$.

This gives  $\int_{0}^{h} G (t+s)\, ds= \Delta_h  \phi (t)$, $ \in X$, $t\in I^{-h}$.
The first equation then gives $G(t) \in X$ for $t\in I^o $, then $t \in I$ with
$\Delta_h G (t) \in X$, one gets $G(I) \st X$.

\noindent Similarly as in the calculations above one shows  $\Delta_h  \phi (t) =  \int_{t}^{t+ h}  G(s)\, ds $.
Since     $G \in C(I, X)$, differentiation with respect to $h$ at $h=0$ is possible and gives:  the norm derivative $\phi' (t)$  exists $= G(t)$, first from $t\in I^o$,
then from $I$ with $\Delta_h \phi \in C^1 (I,X)$ and also $\phi' =G$ there. Since $G$ is continuous by the above, $\phi \in C^1 (I,X)$.

Case $k\Rightarrow  k+1 $: Obvious.\P
\enddemo

\proclaim {Remarks 1.6} (i) Proposition 1.5 holds also if $\phi:
\jj\to X $ is only Lebesgue-measurable instead of locally
integrable (see Lemma 4.19).

 (ii) Lemma 13 [49, p. 224] (see
also [30, p. 200, Lemma 3.1]) holds also for $X$-valued $f$, $\A
:= C^k $, $ M^k_{loc}$, $W^{p,k}_{loc}$.
\endproclaim

\smallskip

\noindent Let  $f \in C_u(\jj,X)$  and $w(t)= 1+|t|$. Set

 (1.14) $\,\,\,\qquad  ||f||_u := || f /w||_{\infty} +$ sup $_{ t\in [0,1]}$ $||f-f_t||_{\infty}$.

\proclaim {Proposition 1.7}

\noindent (i) If  $\phi\in L^1_{loc} (\jj, X)$ and  $\Delta_h\phi \in C_{ub} (\jj, X)$  for
all $ 0 <h \le $  some positive $\e_0$, then $\phi$ is uniformly  continuous.

\noindent (ii)  If  $\phi\in L^1_{loc} (\jj, X)$
and  $\Delta_h\phi \in C_{u} (\jj, X)$  for
all $h> 0$, then       $ ||\Delta_h\phi)/w ||_{\infty} \to 0$ as $h\to 0$, $w(t)= 1+|t|$.

\noindent (iii)
   $C_u (\jj,X)$ endowed with $ ||\cdot||_u $ defined by (1.14) is a Banach space.
\endproclaim

\demo {Proof}  (i) By Proposition 1.5 $\phi \in C (\jj, X)$. For
the $\jj=\r$, the assumption gives   $\Delta_h\phi \in C_{ub}
(\jj, X)$  for all $h\in \r$, then    $\phi \in C_{u} (\jj, X)$ by
[11, Corollary 5.5]. If  $\jj= (\alpha, \infty) \not =\r$,
$\Delta_h\phi$ and therefore $\phi$  can  be continuously extended
to $ [\alpha, \infty)$ and then to $\r$ by $\phi (t)= \phi(\al)$
for $t\in \r\setminus \jj$, implying  $\Delta_h\phi \in C_{ub}
(\r, X)$  for all $ 0 <h \le \e_0$ for  the extension. By the
above        $\phi \in C_u (\jj, X)$.

   (ii) Let  $\phi\in L^1_{loc} (\jj, X)$
and  $\Delta_h\phi \in C_{u} (\jj, X)$  for all $h> 0$. Then
$\Delta_s (\Delta_h\phi) \in C_{ub} (\jj, X)$ for all $ h,\, s
>0$. It follows that $\Delta_s (\phi-M_h \phi )= \Delta_s (\phi)
-M_h \Delta_s (\phi)= \lim _{n\to \infty} (1/n)\sum _{k=1}^{n}
(\Delta_s \phi- (\Delta_s \phi)_{s_k}), \in C_{ub} (\jj, X)$ for
all $  h,\, s>0$, here $s_k = (hk)/n$, $n\in \N$. One gets
$\phi-M_h \phi \in C_{u} (\jj, X)$  for all $0< h \in \jj$, by
(i). Hence $\phi = M_1 \phi +g$ with $g \in C_{u} (\jj, X)$. Since
even $||\Delta_ h g||_{\infty} \to 0$ as $h\to 0$, one has to
consider only $u:= M_1 \phi$. Since $u\in C^1 (\jj, X)$ with $u'
=\Delta_1 \phi \in C_{u} (\jj, X)$, $\Delta_h u =\int _0^{h} u'
(t+s)\, ds$, $||\Delta_s u' ||_{\infty}\le \e$ for all $0< s<
\delta(\e )$. Since $u'/w$ is bounded, one gets $||(\Delta_h u)/w
||_{\infty}\le  h || u'/w ||_{\infty} + \e ||h/w||_{\infty}$ for
all for all $0< h < \delta (\e )$. This proves (ii).

 (iii) That $||\cdot||_ u$ is a
norm is obvious. Using  (ii) one can  show that     $(C_u,
||\cdot||_u)$  is complete. (The strictly stronger topology of
uniform convergence  for $ C_u$ is not linear and  so is not
normable; see  also Kolmogoroff's theorem  [39, p. 91]).    \P

\enddemo

\head{\bf \S 2 New classes of  vector valued
distributions}\endhead

In the following we use freely Schwartz's  theory for Banach space
valued distributions; as remarked in [71, p. 30], most
non-topological results carry over to $\h'(\jj,X)$. Especially,
the convolution $S*T$ is well defined and belongs to  $\h'(\r,X)$,
if $S\in \h'(\r,\k)$  and $T\in\h'(\r,X)$ or vice versa, and
(independently) $S$ or $T$  has compact support. This convolution
is commutative, bilinear and $T*\va\in C^{\infty}(\r,X)$ if $ \va
\in \h (\r,\k)$. Moreover, one has

\smallskip

(2.1) $\,\,\,\qquad\, (S*T)'=S'*T=S*T'$;

\smallskip

(2.2)  $\,\,\,\qquad\,(S*T)_a =S_a*T=S*T_a $, $a\in\r$;

\smallskip

(2.3) $\,\,\,\qquad\,(U*V)*W=U*(V*W)$,

\smallskip

\noindent if one of $U, V, W \in \h'$ is $X$-valued and the  other two
 are $\k$-valued and furthermore (independently ) if two of the $U, V, W$ have
compact support.

\noindent This can be reduced to the scalar case as follows:

\smallskip

\noindent If $S \in \h' (\r,X)$, $T \in \h' (\r,\k)$ with compact
support, then with $\check {\va} (t) = {\va} (-t)$,

\smallskip

(a) $(S*T)(\va):= S(T*\check {\va})$, $\qquad\varphi \in
\h(\r,\k)$;

\smallskip

\noindent If $T$ does not have compact support, then $S$ has,
choose $\rho \in \h(\r,\k)$ with $\rho \equiv 1$ on an open
neighborhood of supp $S$, and define

\smallskip

(b) $(S*T)(\va):= S(\rho\cdot T*\check {\va})$, $\qquad \varphi
\in \h(\r,\k)$;

$S*T:=T*S,\,\quad$ if $\,\,T \in \h' (\r,X)$, $S \in \h' (\r,\k)$.

\smallskip

\noindent  For $X=\k$ this gives the scalar convolution [75, p.
62, Definition 1, p. 156, (2), p. 158, (8) and Theorem 3].

\noindent  In (a) respectively (b)

$T*\check {\va} \in \h(\r,\k)$ respectively $C^{\infty}(\r,\k)$
with $(T*\check {\va})(t)= T(\va_{-t})$ for $t\in \r$ [75, p. 156,
Proposition 1].

\noindent In (b) supp $S$ is defined as in the scalar case, as
there one has $S(\va)=0$ if supp $\va \,\cap$ supp $S =\emptyset$,
$\varphi \in \h(\r,\k)$ [75, p. 62, Definition 1, Theorem 1]; this
immediately shows that the definition in (b) is independent of the
choice of $\rho$.

\noindent So $S*T: \h(\r,\k)\to X$ is  always well defined, linear
and commutative.

By Proposition 2.7 below (with $(T*\check {\va})^{(j)}=T*(\check
{\va})^{(j)}$),

if $\va_n \to 0$ in $\h(\r,\k)$, then

$T*\check {\va}_n$  respectively $\rho\cdot(T*\check {\va}_n) \to
0$ in $\h(\r,\k)$,

\noindent  so  $S*T \in \h'(\r,X)$.

\noindent If $y\in $ dual $X^*$, by definition one has $y\circ
(S*T)= (y\circ S)*T$ respectively   $S*(y\circ T)$ with $y\circ S$
 respectively $y\circ T$ $\in \h'(\r,\k)$.

\noindent (2.1)-(2.3) then follow from the scalar case ([75, p.
 159, (9), (10)]).

\smallskip

\smallskip

 For $L^1_{loc}$-functions distribution convolution is given by
the usual Bochner integral.  With (2.3) and $s_h$ as after (1.10)
one can show (see (1.11))

\smallskip

(2.4) $\,\,\,\qquad\, \widetilde {M_h} T=T*s_h$,  $\qquad T\in
\h'(\r,X)$, $\qquad  0\not = h\in \r$.

\smallskip

If we cite a Theorem of Schwartz (for the scalar case), it is understood that the proofs work also for
the  Banach space case and open  $\jj\not = \r$.

\smallskip

\smallskip

In what follows $\A$ will always be a fixed subset of $L^1_{loc} (\jj,X)$ or $ \h'(\jj,X)$.

\proclaim{Lemma 2.1} If  $T\in \h'_{\A}(\r,X)$, then  $T'\in \h'_{\A}(\r,X)$.
\endproclaim
\demo  {Proof} By (2.1), $T'*\va =T*\va'$ for all  $ \va \in \h
(\r,\k)$; since  $ \va' \in \h (\r,\k)$ the statement follows from
(1.10).
        \P
\enddemo

\proclaim{Lemma 2.2} If   $\A$ is  real-linear and
positive-invariant, $\st L^1_{loc} (\jj,X)$ or $ \h'(\jj,X)$ $
n\in\N $, one has

(a) if  $T\in\A$ with $\jj$ open, then   $T^{(n)}\in \widetilde
{M}^n \A$;

(b)  if $T\in\A$ and  $\jj$ open with distribution derivative
$T^{(n)}\in L^1_{loc}(\jj,X)$, then
   $T^{(n)}\in \m^n \A$  and   $T\in  W^{1,n}_{loc} (\jj,X)$;

(c) if $F  \in  \A  \cap  W^{1,n}_{loc} (\jj,X)$, then $F^{(n)}
\in M^n \A$ (any $\jj$; see
     (1.13)).

\endproclaim

\noindent Here for every interval $I$ of positive length  and $n\in \N$

\smallskip

(2.5) $\qquad  W^{1,n}_{loc} (I,X):=\{f\in C^{n-1}(I,X):
f^{(n-1)}$ locally absolutely continuous

  $\,\,\,\qquad\, \qquad \,  $   on $I$ and $ (f^{(n-1)})'$ exists a.e. in  $ I  \}$,

  $ \,\,\,\qquad\,\,\,\,\qquad\,= \{ T\in \h' (I,X): T^{(n)}\in L^1_{loc} (I,X) \}$ for $I$ open,

$ \,\,\,\qquad\,\ \,\,\,\qquad\,\,\,\,\qquad\W^{1,0 }_{loc}
(I,X):=L^1_{loc} (I,X)$,

\smallskip

\noindent with $f^{(n)}:=0$ where $(f^{(n-1)})'$ does not exist,
then $f^{(n)}\in L^1_{loc}(I,X)$.

\demo  {Proof}
(a) $n=1$:
  $ \widetilde {M_h} T'(\va)=T'(M_{-h} \va) =-T((M_{-h} \va)')=T(\va_{-h}- \va)/h =(1/h)(T_h-T)(\va)$.

\noindent Since with $ \A$ also  $\widetilde {\m} \A$ is linear respectively positive-invariant,
one can complete the proof by induction.

(b)   $ \{T\in  \h'(\jj,X)$:    $T^{(n)}\in L^1_{loc}(\jj,X)\} =
W^{1,n}(\jj,X) $ for open $\jj$

\noindent follows using indefinite integrals  $P^i g$, $g=
T^{(n)}$.

If $n=1$, $T'=g=f'$ as after (2.5), then  $f - T =$ constant, so
$\Delta_h f = \Delta_h T \in  \A$, and $\Delta_h f = hM_h(f')$
with the fundamental theorem of  calculus [53, p. 88, Theorem
3.8.6], yielding $f'  \in \m\A$. The general case follows by
induction with  $f - T  =$ polynomial with coefficients from $X$,
if $f  \in W^{1,n}_{loc}(\jj,X)$ with $f^{(n)} = T^{(n)} $.

(c) follows as in (a) and (b), since  $F'  \in
W^{1,n-1}_{loc}(\jj,X)$.
        \P
\enddemo

\proclaim{Lemma 2.3} If  $\A \st L^1_{loc} (\jj,X)$ or $ \h'(\jj,X)$  satisfies ($\Delta_1$)  and $n\in\N$,  then  $ {\m}^n  \A$ satisfies $(\Delta_1)$; if additionally
 $\jj=\r$, also  $\h'_{\A}(\r,X)$ satisfies $(\Delta'_1)$. The same holds for  $(\Delta)$ instead of $(\Delta_1)$,
$  \h'_{\A}(\jj,X)$ then   satisfies $(\Delta')$.  If $\A \st \h'(\jj,X)$ satisfies $(\Delta'_1)$  (respectively $(\Delta')$),   also   $\widetilde {\m}^n  \A$ satisfies $(\Delta'_1)$ (respectively $(\Delta')$).
\endproclaim

\demo  {Proof}
Case $\m\A$ : If $0< h \in\r$,  $f\in L^1_{loc}(\jj,X)$  with  $\Delta_s f\in \m\A$
for $s > 0$, one gets  $\Delta_s M_h f= M_h\Delta_s f\in \A$. By assumption then
 $ (M_h f-M_a M_h )f\in \A$  for $a=1$ (respectively $ a> 0$).
Since with Fubini-Tonelli;   $ M_h M_a f=M_a M_h f$ on $\jj$ for $a, h > 0$, one
get  $ M_h (f-M_a f )\in \A$ or  $ (f-M_a f )\in \m\A$. With induction,   $\m^n\A$
has ($\Delta_1$)) (respectively ($\Delta$)).

The case  $ \widetilde{\m}^n\A$ is treated similarly with ($\Delta'_1$) (respectively ($\Delta'$)), since  if $s, h, k > 0$

\smallskip

(2.6) $\,\,\,\qquad\,   \Delta _s\widetilde {M_h}  T=   \widetilde
{M}_h\Delta_s  T\,\,$ and    $\,\, \widetilde {M} _h\widetilde
{M}_k T=   \widetilde {M}_k\widetilde  {M}_h T$ for   $ \,\, T\in
\h'(\jj,X)$.

\smallskip

Case $  \h'_{\A}(\r,X)$: Assume $T \in  \h'(\r,X) $    with $\Delta_h T \in  \h'_ {\A}(\r,X)$,
$h > 0$, i.e.   $\Delta_h (T*\va) =(\Delta_h T)*\va \in \A$ for  $\va \in  \h (\r,\k)$,
with $T*\va \in C^{\infty}(\r,X)$. By assumption
 $ T*\va- M_a (T*\va)\in \A$. Now by (1.7), (2.3) and (2.4)

\smallskip

(2.7) $\,\,\,\qquad\,  M_a(T*\va) =  (T*\va)* s_a  =(\widetilde
{M}_a T)*\va$.

\noindent Together one gets $ (T-\widetilde {M}_a T)*\va   \in\A$
for
 $\va \in  \h (\r,\k)$ or  $(T-\widetilde {M}_a T)\in    \h'_{\A}
 (\r,X)$. \P
 \enddemo
\proclaim{Proposition 2.4}
 If $\A \st  \h'(\r,X)$,

\noindent(a)  $\,\,\,\qquad   \h'_{\A}(\r,X)\st \widetilde {\m}
\h'_{\A}(\r,X)$.

\noindent (b) If furthermore  $\A$  is a cone  satisfying
$(\Delta_1)$,  then

  $\qquad\qquad  \widetilde {\m} \h'_{ \A}(\r,X)= \h'_{\A} (\r,X) $.
\endproclaim

              \demo  {Proof} (a)
  $T\in  \widetilde {\m} \h'_{\A}(\r,X)$ means  $  \widetilde {M_h}T\in  \h'_{\A}(\r,X)$   or with (2.3) and (2.4),
$(T*s_h)*\va=T*(s_h*\va) \in \A $, $h >0, \va \in  \h(\r,\k)$; since $s_h*\va  \in  \h(\r,\k)$ and  $T\in  \h'_{\A}(\r,X)$,
this gives the proof of the inclusion.

(b) ``$\st$":   Let  $T \in  \h'(\r,X)$ with $ \widetilde {M_h}
T\in    \h'_{\A} (\r,X)$, $h >0$. By [71, p. 51, Th\'{e}or\`{e}me
I] there exists $ S\in    \h' (\r,X)$ with $S'=T$ on $\r.$ Since
(any $\jj$)

\smallskip

 (2.8) $\,\,\,\qquad\, \widetilde {M}_h (S')= (1/h)\Delta_h S$, $h> 0$

\smallskip

\noindent and since $  \h'_{\A}(\jj,X)$  is also a cone, $\Delta_h S  \in  \h'_{\A}(\jj,X)$ for $h > 0$.
With $(\Delta_1)$, Lemma 2.3 and Lemma 2.1,

 $T-(\widetilde {M}_1 S)' = (S-\widetilde {M}_1 S)'=: U \in  \h'_{\A}(\r,X)$.
Now ( also for open $\jj \not =\r$)

\smallskip

 (2.9) $\,\,\,\qquad\, (\widetilde {M}_a S)' = \widetilde {M}_a( S')$.

\smallskip

\noindent Hence    $T=(\widetilde {M}_1 S)'+U = \widetilde {M}_1 (S')+ U=  \widetilde {M}_1 T+ U \in  \h'_{\A}(\r,X)$.

``$\supset$" follows by (a).
        \P
\enddemo

\noindent The monotonicity of the $ \widetilde {M_h}$-operator gives

\proclaim{Corollary 2.5} If  $\A $  is as in Proposition  2.4 (b) with  $\A * \h(\r,\k)  \st \A$, one has

\smallskip

(2.10) $\,\,\,\qquad\,   {\m}^n\A \st \widetilde {\m}^n \A  \st
\h'_{ \A}(\r,X),\qquad\,n \in \N_0$.

\smallskip

\endproclaim
\proclaim{Corollary 2.6} For  $\A $   as in Proposition  2.4 (b), one has

     $ \h'_{{\m}^n\A} (\r,X) =   \h'_ {\widetilde {\m}^n  \A}(\r,X)= \widetilde {\m}^n \h'_{ \A}(\r,X) =  \h'_ { \A}(\r,X) $,$\qquad\, $
              $n \in \N_0$.
\endproclaim

\demo {Proof}  The definitions, (2.4) and (2.3), since $\widetilde {\m}^n\A \cap L^1_{loc} (\r,X)  ={\m}^n\A (\r,X)$.
        \P
\enddemo

 In the following (see [71, p. 72, Th\'eor\`eme IX])

\smallskip

  $U \st   \h'(\jj,X) $  $bounded$  means  sup$\{ ||T(\va)||: T\in U\}< \infty$  for each $ \va \in  \h (\jj,\k)$,

\smallskip

$ \h_{K}:=  \h_{K}(\jj,\k) := \{ \va\in  \h(\jj,\k): \text {supp}\,\,\va  \st K\}$,

\smallskip

 $ ||\va||_{\infty, m}:= \sum _{j=0}^{m} \text {sup}_{x\in \jj} |\va^{(j)}(x) |$.

\proclaim{Lemma 2.7} If    $U \st   \h'(\jj,X) $  is bounded and $K\st \jj$ is compact, then  there exists $m\in\N$  with

\smallskip

(2.11) $\,\,\,\qquad\,    ||T(\va)||\le m\, ||\va||_{\infty
,m}\qquad$ for all $T\in U$,  $\qquad\va\in \h _{K}(\jj,\k)$.

\smallskip

\endproclaim

\demo {Proof}    It is enough to find  some ball $K_{\psi,m,\e}=
\{ \va\in  \h_{K}:  ||\va -\psi||_{\infty, m}\le \e \}$  with sup
$\{  || T(\va)||: T\in U, \va  \in\h_ { K_{\psi,m,\e}} \} <
\infty$. Assuming this  not to be the case, for $K= K_{0} =
K_{0,0,2}$ there exists $T_1 \in U$,   $\psi_1\in K_{0,0,1}$ with
$||T_1(\psi_1)|| > 1+1$. By  [71, Th\'{e}or\`{e}me XX p. 82] there
is a ball  $K_{1} = K_{\psi_1,m_1,\e_1}\st K_0$ with $||T_1(\va)||
> 1$ on $K_1$. Recursively one finds $T_n\in U$,
  $\psi_n\in \h_{K_n}$, $m_n$, $\e_n > 0$ with $m_n\to \infty$, $\e_n\to 0$,
  $K_n= K_{\psi_n,m_n,\e_n} \st K_{n-1}$  and  $||T_n(\va)|| \ge n$ for $\va \in K_n$. The completeness of
$\h_K$  gives $\psi\in \h_K$  with $\psi_n\to\psi$, then $\psi \in \cap_{n=1}^{\infty} K_n$.
But then   $||T_n(\psi)|| \ge n$, contradicting the boundedness of $U$.
        \P
\enddemo

\proclaim{Lemma 2.8} If  $F_n (t):= t^{n-1}/((n-1)!)$  for  $t>0$,
$:=0$ for $t \le 0, n \in \N$, then $F_n^{(n)}=\delta$,
$F'_{n+1}=F_n$  in $\h'(\r,\r)$ and for  $T \in   \h'(\r,X) $
there exist $\gamma, \zeta_n \in\h(\r,\r)$ with

\smallskip

(2.12) $\,\,\,\qquad\,   T=  ((\gamma F_n)* T)^{(n)}+ \zeta_n*T$,
$\qquad \,n\in\N$.

\smallskip

\endproclaim
 \demo{Proof  of (2.12)} (See [71, (VI.6.22), p. 191]) : Choose $  \gamma \in \h(\r,\r) $ with $\gamma =1$ in
 a neighborhood of $0$, $\zeta_n := \delta- (\gamma F_n)^{(n)}$. Then $\zeta_1 =
 \delta- \gamma F'_1 - \gamma' F_1=-\gamma' F_1  \in \h(\r,\r)$. Inductively $\zeta_{n+1} :=
 \delta- (\gamma F'_{n+1}+ \gamma' F_{n+1})^{(n)} =\delta- (\gamma F_{n}+ \gamma' F_{n+1})^{(n)}
 =\zeta_n- ( \gamma' F_{n+1})^{(n)}, \in\h(\r,\r) $.  \P
\enddemo
\smallskip

\proclaim{Proposition 2.9}
  $\h'_{L^{\infty}} (\r,X)=  \h'_{C_{b}} (\r,X)  = \h'_{C_{ub}} (\r,X) =  \h'_{(\h_{L^{\infty} })} (\r,X)$

$\qquad  \qquad\,\,\,\qquad \qquad\,\,\,\qquad \qquad\,\,\,\qquad
= (\h_{L^1})' (\r,X) \,| \h (\r,\k)$.
\endproclaim

 \proclaim{Remark}
Our  $\h'_{L^{\infty}} (\r,X)$ of Definition 1.3  is  a priori
different from  the

\noindent $(\h'_ {L^{\infty}})  (\r,X)  = (\h_{L^{1}})' (\r,X)$ of
[71, p. 200].
\endproclaim

\demo {Proof}  Last  ``$\supset$":
 If  $ T :   \h_{L^1}(\r,\k) \to X$  is linear and sequentially continuous, the standard
contradiction argument gives $m_0\in \N$ so that, with
$||\va||_{1,m_0} :=  \sum _{j=0}^{m_0}  ||\va^{(j)}||_{L^1(\r,\k)}$,

\smallskip

(2.13) $ \qquad  \qquad ||T(\va)|| \le m_0 ||\va||_{1,m_0}$, for all    $ \va \in   \h_{L^1}(\r,\k) $.

\smallskip

\noindent With this, $T* \va(x)= T((\check {\va})_{-x})$ is bounded in $x \in \r$
for fixed  $ \va \in   \h(\r,\k) $, so  $ T \in   \h' _{L^{\infty}}(\r,X) $.

\noindent Since
 $T* \va \in  C^{\infty}(\r,X) $   if  $ \va \in   \h(\r,\k)$,
 $(T* \varphi)^{(n)}=T* (\varphi^{(n)})$ and
  $\h_{L^{\infty}} (\r,X)\st  {C_{ub}} (\r,X)  \st {C_{b}} (\r,X) \st {L^{\infty}} (\r,X)$,
 the first three equalities follow.

\noindent For  $ \h'_{L^{\infty}} (\r,X)\st (\h_{L^1})' (\r,X)$ we show, for  $\A\st \h'_{L^{\infty}} (\r,X)$
 (see [71, p. 201-202]):

If $\A$ is $C^{\infty}$-uniformly closed,  $ \h'_{\A}
\st\h'_{L^{\infty}} (\r,X)$, then to
 $T\in  \h'_{\A} (\r,X)$ there exists $F,G \in \A \cap  C(\r,X)$ and $n\in \N$ with

\smallskip

(2.14) $ \qquad  \qquad         T  = F + G^{(n)}$.

\smallskip

\noindent With (2.12) one can take $F=\zeta_n*T \in\A \cap C^{\infty}(\r,X) $, it remains to show that $ (\gamma F_n)*T\in\A$ for suitable $n$ ($\gamma$  is independent of $n$):

$U:=\{T_x: x\in\r\}$ is bounded in   $\h' (\r,X) $, since $T_x(\check{\va})=T((\check{\va})_{-x})=(T*\va)(x) $

\noindent  is continuous in $x$,  $ \st L^{ \infty}(\r,X)$. If supp $\gamma \st [-a,a]$,  to $U$ and $K:=[-a-1, a+1]$ there is $m$ with (2.11) by Lemma 2.7.
Choose $n:=m+2$, then $F_n\in C^m(\r,\k)$, supp $\,\gamma F_n \st [-a,a]$.
Therefore there exists  $(\va_j) \st \h_{K} (\r,\k)$  with  $||\va_j-\gamma F_n||_{\infty,m}\to 0$  as $j\to \infty$. (2.11) shows that $(T*\va_j)_{j\in\N}$
is Cauchy with respect to uniform convergence on $\r$, so by the assumption on $\A$
one has    $T*\va_j\to  G\in \A$, uniformly on $\r$, $G\in  C_b(\r,X)$.

\noindent One has  $(\gamma F_n)*T=G$: For this    $((\gamma F_n)*T)*\va=G*\va$ for
$\va \in\h (\r,\k)$ is enough. Now
 $(T*\varphi_j)*\va\to G*\va$ (uniformly) on $\r$,  $\varphi_j\to \gamma F_n$ uniformly  on $\r$ implies
$\psi_j:=\va_j*\va \to (\gamma F_n)*\va$ in $\h (\r,\k)$, so
$((T*\va_j)*\va)(x)= (T*\psi_j)(x)=T((\check {\psi}_j)_{-x})\to $

$ T( ((\gamma F_n *\va\check {)\,} ) _{-x}) = (T*(\gamma F_n)*\va))(x)= (((\gamma F_n)*T)*\va)(x) $  for each $x\in\r$.

\noindent This gives the desired result  $G*\va=((\gamma F_n)*T)*\va$, and so (2.14).

\noindent Since     $   L^{\infty} (\r,X)\st (\h_{L^1 })' (\r,X)$ and the latter
is closed with respect to differentiation and addition, (2.14) for   $\A = L^{\infty} (\r,X)$ gives

  $\h'_{L^{\infty}} (\r,X) \st (\h_{L^1})' (\r,X) | \,\h (\r,\k)$.
        \P
\enddemo

\noindent The extension process $ \A\to   \h'_{\A} (\r,X)$     is iteration complete:

\proclaim{Theorem 2.10}      If  $\A\st \h'(\r,X)$ is closed under addition, then
  $ \h'_{\h'_{\A}} (\r,X)=\h'_{ \A} (\r,X) $.
\endproclaim
 \demo {Proof}
 ``$\supset $":  Follows from the definitions, associativity  (2.3) and

 $\h (\r,\k) *\h (\r,\k)\st \h (\r,\k) $.

``$\st$":
If   $ T\in \h'_{\h'_{\A}} (\r,X)$,   $(T* \va)*\psi= T* (\va*\psi)\in \A$  for all
  $ \va,\psi \in\h (\r,\k)$. By [66, Theorem 3, p. 554] (see also  [42, pp. 584, 587], [36]) one has

\smallskip

 (2.15) $ \qquad  \qquad    \h (\r,\cc)= \h (\r,\cc)* \h (\r,\cc) +  \h (\r,\cc)* \h (\r,\cc)$.

\smallskip

\noindent This implies

\smallskip

 (2.16) $ \qquad  \qquad    \h (\r,\r)= \sum _{1}^{4}\h (\r,\r)* \h (\r,\r)$.

\smallskip

\noindent By the above and the assumption on $\A$, $T*\va\in \A$ for all   $ \va\in    \h (\r,\k)$, i. e.  $T\in  \h'_{\A} (\r,X)$.
        \P
\enddemo

\proclaim
 {Remark} For $\A$  as in the following Theorem 2.11, Theorem 2.10 can be shown
directly without [66], using an analogue to Lemma 2.7 for $T*\va*
\psi$,   $ T\in    \h'_{\h'_{\A}} (\r,X)$.
\endproclaim

\proclaim{Theorem 2.11} If  $\A \st\h'_{ L^{ \infty}}(\r,X)$ is $C^{\infty}$-uniformly closed, closed under addition and   $\A * \h (\r,\k)  \st  \A$, then for   $T  \in \h'( \r,X)$ the following statements are equivalent:

\noindent (a)  $T \in  \h'_{ \A} (\r,X)$, i.e.  $T*\va \in \A $ for each  $\va \in  \h( \r,\k)$ ;

\noindent  (b) there exist $F , G \in  \A \cap C(\r,X)$ and a nonnegative integer $n\in \N_0$, such that

         $T  =  F  +  G^{(n)}$       $\qquad$  (distribution derivative);

\noindent  (c) there exist  $m\in \N$,  $k_j\in \N_0$,  $F_j \in  \A$
with         $T  = F_1^{(k_1)}+\cdots   + F_m^{(k_m)}$;

\noindent  (d)  $T \in  (\h_{L^1})' (\r,X)$  and there exists
$(\phi_n)_{n\in \N} \st  \A \cap C^{\infty}(\r,X) $
    with  $ \phi_n \to T$ in       $(\h_{L^1})' (\r,X)$;

\noindent  (e)  there exists   $(\phi_n) \st  \A$ :

      $ \phi_n \to T$ uniformly on each $U$ with  $ U\st \h (\r,\k) $, $U$  bounded in    $\h_{L^1} (\r,X)$;

\noindent  (f) $T \in $  closure of $\A$ in  $(\h_{L^1})' (\r,X)$.

\endproclaim
\proclaim{Remark 2.12 } Theorem 2.11  becomes false if  $\A$ is not $C^{\infty}$-uniformly closed:  $\h'_{ C_c}(\r,X)\not = (\h_{ L^1})'(\r,X) $-closure of $C_c (\r,X)$, $=\h'_{ C_0}(\r,X)$.
However one can show directly

  $\h'_{ \h}(\r,X)= \h'_{ C_c}(\r,X) = \{ T\in \h'(\r,X):$ supp $T$ compact  $\} =\frak  {E}' (\r,X)$ (see [75, p. 62] for the case $X=\cc$).
\endproclaim

\demo {Proof of Theorem 2.11} $(a)\Rightarrow (b)$:
This has been shown in the proof of (2.14) above; the assumption  $\h'_{\A} (\r,X)\st  \h'_{L^{\infty}} (\r,X)$ needed
is now fulfilled, since   $\A\st  \h'_{L^{\infty}} (\r,X)$  implies
 $\h'_{\A} (\r,X)\st \h'_{ \h'_{L^{\infty}}}(\r,X)= \h'_{L^{\infty}} (\r,X)$ by
Theorem 2.10.

 $(b)\Rightarrow (c)$: trivial.

 $(c)\Rightarrow (d)$: Since $\A$ is closed under addition, one can assume
 $T  =  F^{(k)}$, $F\in\A$. Then  $\phi_n  = T*\rho_n = F*\rho_n^{(k)} $ again $\in\A$,  $0\le \rho_n \in \h (\r,\k)$,  $\int_{\r} \rho_n (x)\,dx=1$,
supp $ \rho_n \st [-(1/n), 1/n] $. With $\A\st \h'_{L^{\infty}} (\r,X) $,

Lemma 2.1 and Proposition 2.9 one can assume  $ T \in  ( \h_{L^1})' (\r,X) $.
Assuming $U$ bounded  $ \st  \h_{L^1} (\r,\k)$, one has to show uniformly in $\va \in U$,

 $ \phi_n (\va)=\phi_n * \check {\va} (0) = T* (\rho_n * \check {\va})(0)= T ( \check {\rho}_n * \va) \to  T(\va)$.

\noindent since $T$ satisfies (2.13) by the proof of Proposition 2.9, it is enough to show
 $||\check {\rho}_n * \va) -\va||_{1,m} \to 0$ uniformly in $\va \in U$. But this
follows with $\psi =\va^{(j)}$  and $a=1/n$ from

$ || \rho * \psi -\psi||_{L^1} \le 2a ||\rho||_{L^1} ||\psi'||_{L^1}$  if supp $\rho \st [-a,a]$,

\noindent   $\psi \in C^1 (\r,\k)$ with   $\psi, \psi' \in L^1 (\r,X)$.   The  $\psi$-inequality one gets with

 $\psi(x)- \psi(x-t)=  \int_{x-t}^{x} \psi'(s) \,ds= \int_{-t}^{0} \psi'(x+s) \, ds $ and twice Fubini-Tonelli.

 $(d)\Rightarrow (e)$: obvious.

 $(e)\Rightarrow (a)$:  Since
 $ U_{\va}=\{ (\check{ \va})_{-x} : x \in \r \}$ is bounded in   $  \h_{L^1} (\r,\k)$
for fixed  $ \va\in \h (\r,\k)$,  $ \phi_n*\va \to T* \va $ uniformly on $\r $, so
$ T*\va\in \A$ by  the assumption on  $\A$.

 $(d)\Rightarrow (f)$: obvious. Here $ T\in \overline{\A}$ means that there exists
an extension $ \widetilde {T}\in   (\h_{L^1})' (\r,\k)$ with     $ \widetilde {T}\in  \overline{\A}$.

 $(f)\Rightarrow (a)$: By the definition of the topology of  $  (\h_{L^1})'$ above, for  $  \va\in \h (\r,\k)$,
to  $ U_{\va}=\{ (\check{ \va})_{-x}  \}$, $ \e_n = 1/n$  there exists $ F_{n,\va}\in\A$  with
$| F_{n,\va}*\va -T*\va |  \le \e$  on $ \r$, $ F_{n,\va}* \va \in \A$.
        \P
\enddemo
\proclaim{Corollary 2.13}
If
  $\A$   is as in Theorem 2.11, one has

  $\,\,\,\qquad\,      \h'_{\A} (\r,X) =  \overline{\A}$ (closure of  $\A$ in   $(\h_{L^1})' (\r,\k)
=$sequential closure of $\A$

 \noindent (here ``sequential closure''  means (d) of Theorem 2.11)).

\endproclaim

\proclaim{Corollary 2.14}
If
  $\A$   is as in Theorem 2.11,  and
  $      \A_* \st U \st \h'_{\A} (\r,X) $ then  $ \h'_{U} (\r,X)= \h'_{\A} (\r,X) $.

\endproclaim

Here $\A_{0}=\{\sum_{j=1}^{m} \phi_j* \va_j: m\in\N, \phi_j \in \A,  \va_j\in \h (\r,\k) \}$, $\st \A\cap C^{\infty} (\r,X)$.

\demo{Proof} If $ T\in \h'_{\A} (\r,X)$,   one has $T= F+G^{(m)}$
with   Theorem 2.11 (b), $ F, G\in\A$, so  $T*\va=
F*\va+G^{(m)}*\va =  F*\va+G*\va^{(m)}$,  $ \in\A_{0}$, $ \va\in
\h (\r,\k) $. This means  $  \h'_{\A} (\r,X) \st \h'_{\A_{0}}
(\r,X)$, $ \st \h'_{\A} (\r,X)$ since $ \A_* \st \A$, so $
\h'_{\A} (\r,X)= \h'_{\A_*} (\r,X)$. Then $ \A_{0}\st U\st
\h'_{\A} (\r,X)$ gives

$  \h'_{\A} (\r,X)=  \h'_{\A_*} (\r,X) \st \h'_{U} (\r,X)\st\h'_{
\h'_{\A}} (\r,X)$
 $= \h'_{\A} (\r,X)$

\noindent  with Theorem 2.10.
        \P
\enddemo
\noindent The special cases $U = \m^n \A$ or  $U= \widetilde {\m}^n \A$ follow already with Corollary 2.6 for more general  $\A$.

\proclaim{Corollary 2.15}
If
  $\,\, \A\st \h' _{ L^{\infty}} (\r,X)\,\, $ is\,\, real \, linear,\,\, positive-invariant,
  $\,\,C^{\infty}$-
  uniformly closed,
 $ \A*  \h (\r,\k)\st   \A$, $ \A\st \m  \A$  and if $(\Delta_1)$)  holds for $\A$, then

\smallskip

 (2.17) $\,\,\,\qquad\,      \h'_{\A} (\r,X)  =\cup_{n=0}^{\infty} \widetilde {M}^n\A$.

\smallskip

\endproclaim
\noindent This implies, for $\A$ as in Corollary 2.15, since

   (2.18) $\qquad \,\,\widetilde {\m}^n  \A   \cap       L^1_{loc} (\jj,X)=  {\m}^n  [\A   \cap       L^1_{loc} (\jj,X)] = \m^n \A \,$ for $n\in\N$,

\smallskip

 (2.19) $\qquad \,\,      \h'_{\A} (\r,X)\cap  L^1_{loc} (\r,X)  =\cup_{n=0}^{\infty} M^n\A$ with
   $M^0\A:=    \A\cap  L^1_{loc} (\r,X)$.

\smallskip

 \demo {Proof}
If $ T\in \h'_{\A} (\r,X)$, with Theorem 2.11 (b) one has $T= F+G^{(m)}$  with  $ F, G\in\A$.
 $ \A\st \m  \A \st \widetilde {\m}  \A$ implies  $ \A\st \widetilde {\m}^{m}  \A$, then  $ F\in \widetilde {\m}^{m} \A$;
 since also  $ G^{(m)}\in \widetilde {\m}  \A$ by Lemma 2.2, and  $  \widetilde {\m}^m  \A$ is with $\A$ real linear, one gets
 $ T\in \widetilde {\m}^m  \A$.

Conversely, $ \A*  \h (\r,\k)\st \A$ implies $ \A  \st  \h'_{\A}$,
 $ \widetilde {M}^n\A \st \widetilde {M}^n \h'_{\A} $, $= \h'_{\A}$ by Proposition 2.4, with $(\Delta_1)$.
        \P
\enddemo

Special cases (for $S^p _b$ see  \S 1 or [13, p. 132,  (3.3)],
$L^{\infty} \st S^p_b$):

\smallskip

 (2.20)  $\qquad  \A\st  S^p_{b}(\r,X)$ with some   $1 <p < \infty$, $\A$ $C^{\infty}$-uniformly closed,

$\qquad\qquad\qquad\qquad\qquad \A*\h (\r,\k)\st \A$ implies $ \A\st \m\A$.

\smallskip

(2.21)  $\qquad \A\st  S^1_{b}(\r,X)$, $\A$ linear, invariant,  uniformly closed,

$\qquad\qquad\qquad\qquad\qquad \A\st \m\A$ implies
 $\qquad   \A*\h (\r,\k)\st \A$.

\smallskip

\smallskip

\noindent For a Theorem 2.11-analogue for $\f'(\r,X)$ we need

 \smallskip

(2.22) $\qquad\qquad   \p(\jj,X) :=\{  f\in C (\jj,X): \, \, f= O(|x|^n )$ as $ |x| \to \infty $  for some  $n\in \N \}$.

\smallskip

\noindent With $w_k (x)= (1+x^2)^{k/2}$ for $x, k\in \r$ one has

\smallskip

(2.23) $\qquad\qquad   \p(\r,X) :=\{ gw_{k} : g\in C_b(\r,X), k\in \N \}$.

\smallskip

\noindent The topology  of $\f(\r,X)$ is given by (see [71, p.
 234, (VII,3;2)]) the seminorms

\smallskip

(2.24) $\qquad\qquad  ||\psi||_{w_k , m} :=\sum _{j=0}^{m} ||\psi^{(j)} w_k||_{\infty}, \qquad m,k \in \N_0 $.

\smallskip

\proclaim{Theorem 2.16} If     $T  \in \h'( \r,X)$ the following statements are equivalent:

\noindent (a)  $T \in  \h'_{  \p} (\r,X) $;

\noindent  (b)     $T*\va \in  \f' (\r,X)$ for each  $\va \in  \h( \r,\k)$;

\noindent  (c) there exist  $F\in  \p( \r, X)$,  $m\in \N_0$ with  $T= F^{(m)} $ (distribution derivative );

\noindent  (d) there exists $m_o\in\N$ with $|| T(\va)|| \le m_o ||\va || _{w_{m_o} , m_o}$  for all $ \va \in \h (\r,\k)$;

\noindent  (e)    $T \in  \f' (\r,X) |  \h( \r,\k) $;

\noindent  (f) there exist  $(S_n)_{n\in\N}  \st     \f' (\r,X) $ with $S_n(\va)\to T(\va)$ for $ \va\in  \h( \r,\k)$

$ \qquad\qquad$ and $ (S_n(\psi))$ is Cauchy in $X$
for each  $\, \psi\in \f(\r,\k)$.
\endproclaim

\proclaim{ Remark} If (f) holds, then  $S_n \to \widetilde {T}$ in the  $\f' (\r,X)$- topology (i.e. uniformly on $U$, $U$ bounded   $\st \f( \r,\k)$)
with $ \widetilde {T} \in \f'( \r,X)= $ unique extension of $T$; even $\h( \r,X)$ is sequentially dense in $\f' ( \r,X)$.

\endproclaim
\demo {Proof}   (For $X=\cc$ part of this theorem is shown by
Schwartz [71, Th\'{e}or\`{e}mes IV, VI, pp .238-241, 244-245]):
 $(c)\Rightarrow (d)\Rightarrow (e)\Rightarrow  ``(d)$  for   $ \va \in\f ( \r,\k)"$ $ \Rightarrow  (a)\Rightarrow (b)$
 follows as in [71, pp. 238-241] respectively the proof of Theorem 2.11 with (for  $(d)\Rightarrow (a)$)

\smallskip

(2.25) $\qquad\qquad w_k (x+y)\le  2 ^{k/2} w_k(x) w_k (y), \qquad x,y\in  \r, k>0 $.

\smallskip

 $(b)\Rightarrow (a)$ follows from  $(e)\Rightarrow (a)$ and Theorem 2.10 with $\A= \p$.

 $(d)\Rightarrow (f)$:  $S_n := T*\rho_n \in  \p (\r,X)\st \f'(\r,X)$ by   $(d)\Rightarrow (a)$, $\rho_n$ as in the proof of theorem 2.11.
  $\rho_n*\va \to\va$ in $\h(\r,\k)$ gives  $S_n  \to T$ on  $\h(\r,\k)$.
$(S_n(\psi))$ Cauchy:  By  $(d)\Rightarrow (e)$ there is $S\in  \f'(\r,X)$ with $S | \h(\r,\k) =T$, one has for $\psi \in  \f (\r,\k)$, $\rho_n =\check{\rho}_n $,

 $ S_n(\psi)= (S *\rho_n)* \check {\psi}(0)= S *(\rho_n* \check {\psi})(0)= S(\rho_n*  {\psi})$ by  associativity (2.3) for
 $ U\in \f'(\r,X)$,  $ V, W\in \f(\r,\k)$ ( for $X=\cc$ :  [71, Th\'{e}or\`{e}me XI, p. 247-248]).   Since $\rho_n *\psi \to \psi$ in  $\f( \r,\k)$    by (2.26) below,
 $S_n (\psi) \to  S (\psi)$.

\smallskip

(2.26) $\,\,\,\qquad w_k (x) ||\psi (x)-\psi (y)||\le 3 |x-y|
\cdot ||\psi||_{w_k,2}$,

 $\qquad\qquad\qquad\qquad\qquad\qquad x, y \in \r$,  $|x-y|\le 1$, $\psi \in  C^2 ( \r, X)$, $k\ge 0$.
\smallskip

 $(f)\Rightarrow (d)$:  Lemma 2.7 and its proof hold also for $U$ bounded $ \st \f'( \r,X)$,
with $||(\cdot)||_{\infty, m}$ replaced by      $||(\cdot)||_{w_m , m}$ and
$\h_{K}$ by $\f( \r,\k)$, $U$ bounded meaning sup $\{ ||S(\psi)||: S\in U \} < \infty$ for each  $ \psi\in \f ( \r,\k)$.
This gives (d) for all $S_n$  instead of $ T$, $\va \in \f( \r,\k)$ and so (d) for $T$.

 $(a)\Rightarrow (c)$:  This is the main part of Theorem 2.16, we give a proof  different  from Schwartz
  [71, p. 240-241, case $X=\cc$] and shorter:

For fixed  compact  $ K\st \r$ and $ n\in \N$, define  $ V_n := \{
\va\in \h_{K}  (\r,\k) : ||(T*\va)(x)|| \le n w_n (x), \, x\in
\r\} $. (a) and (2.23) give $\h_{K} : = \h_{ K} ( \r,\k)=\cup
_{n=1}^{\infty} V_n$. The $V_n$ are closed in $ \h( \r,\k) $,
since $T \in \h'( \r,\k)$. $ \h_{K}$  being complete metric, by
the Baire category theorem there are $n_o$,  $\e_o >0$, $\va_o \in
\h_{K}$, $m_o$  with $\{ \va\in  \h_{K} :  ||\va-\va_o|| _{\infty,
m_o}\le  \e_o \} \st V_{n_o}$.  With suitable $c$ this gives $
||T*\va|| \le  c w_{n_o} (x) ||\va||_{  \infty,m_o}$, $\va \in
\h_{K}$. Choosing $n= m_o +2$, $\gamma$, $K$ and $\va_j \in
\h_{K}$  as in the proof of (2.14), one gets $(T*\va_j)/ cw_{n_o}
\to \text {some } g \in C_b( \r,X)$, uniformly on $\r$, or
$(T*\va_j)\to  c\,w_{n_o} g = : G  \in \p ( \r,X)$, locally
uniformly  on $\r$. As in the proof of (2.14), one shows
$T*(\gamma F_n)= G$. Lemma 2.8  gives $T= F+ G^{(n)}$ with $F, G
\in \p(\r,X)$. Now with $ f \in \p (\r,X)$    also the indefinite
integral $Pf  \in \p (\r,X)$, so $H : = P^n F  \in \p (\r,X)$.
Hence $T=H^{(n)}+ G^{(n)}$. This gives (c).   \P
\enddemo

\proclaim{Corollary 2.17}
If
  $\A$   is as in Theorem 2.11  and
  $      T \in \h' (\r,X) $,   $      T \in \h'_{\A} (\r,X) $ is also equivalent with

(g)  there is  $      \widetilde {T} \in  \f' (\r,X) $ with   $\widetilde {T}| \h (\r,\k) = T $ and  $ \widetilde {T}*\psi \in \A$ for  all $\psi \in  \f (\r,\k) $.

\endproclaim
\demo {Proof}
  $ \h'_{\A} (\r,X)\st \h'_{ \h'_{L^{\infty}}} (\r,X)=  \h'_{L^{\infty}} (\r,X)\st \f' (\r,X) |\, \h (\r,\k)$ by Theorems 2.10 and 2.16.
Since   $\h (\r,\k)$ is dense in  $\f (\r,\k)$ and $ \A$ is
$C^{\infty}$-uniformly closed, it is enough to show that
  $ ||T*\va_n ||_{\infty }  \to 0 $  if  $ (\va_n)\st \h (\r,\k) $, $\va_n \to 0$ in $ \f (\r,\k)$.
This follows from  Theorem 2.11 (b) for $\A= L^{\infty}$ since  $||\va^{(m)} ||_{L^1} \le (\int_{-\infty}^{\infty} \frac {dt}{1+t^2} ) \,|| \va||_{w_2, m}  $     for $\va \in \h (\r,\k)$.
        \P
\enddemo

\head{\bf \S 3   Mean classes }\endhead

 In this section we discuss mainly the relation  $\A  \st  \m\A$
and give various examples.
 Concerning the definition of the mean
class $ \m\A$, if $\A$ is real-linear and positive-invariant $ \st
L^1_{loc}(\jj,X)$, for $f  \in   \m\A$  one of the following is
sufficient :

\smallskip

 (a)   $H(f,\A) : = \{h > 0 : M_h(f)  \in  \A \}$ has an
interior point.

\smallskip

(b)  $H(f,\A)$ contains two rationally independent $a, b$,
provided additionally
           $\A$ is uniformly closed and  $Pf  \in  C_u(\jj,X)$, e.g. $f  \in
           S^p_b(\jj,X)$ with         $ 1<p\le \infty$.

\smallskip

  (c) $ M_h(f)  \in  \A$ for just one positive $h_0$, provided $\A$ is a
$\lambda$-class with e.g. $(L_{ub})$ and $Pf$ uniformly continuous
(see (b));
 a special case
          is $\A = AP(\r,X)$ with $c_0 \not \in X$.

\smallskip

\noindent      This is a kind of Tauberian theorem: namely,
$M_h(f) = f*s_h  \in  \A$ for
         some $h = h_0$ and $f  \in S^p_b $ (e.g. $ L^{\infty}$) implies  $f*s_h  \in
         \A$ for all $h > 0$ (see e.g. [19, p. 130, Wiener Tauberian Theorem]).

\noindent A $\lambda$-$class$ means here $\A$ is a closed linear
subspace of $C_{ub}(\jj,X)$ which is $C_{ub}$-invariant under
translations and multiplication by characters
        $\gamma_{\omega}$  and contains the constants (see [13, pp. 117-119]).

\smallskip

\noindent   To see this, we first note, for convex and
positive-invariant $\A  \st L^1_{loc}
     (\jj,X)$, that

\smallskip

  (3.1)     $\,\,\,\qquad H(f,\A)$  is an additive semigroup.

\smallskip

  This follows
from

\smallskip

 (3.2)   $\,\,\,\qquad M_{h+k} (f)  = (hM_h(f) + k (M_k(f))_h)/(h+k)$, $h, k >
 0$.

\smallskip

\noindent The assumptions of (a) and (3.2) give  $(0,\e)  \st
H(f,\A)$ with positive
   $ \e$, then $H(f,\A) = \r^+$  by (3.1).

\smallskip

  To get (b), by the assumptions it is enough to show that $H(f,\A)$ is
dense in $\r^+$. Now by the Kronecker approximation theorem ([52,
p. 436, (iv)]) to $t>0$ and $0< \e <t$ there are integers $m$,
    $n$ with  $|ma+nb-t| < \e$, so $ma+nb >0$. One can assume $m>0$; if
$n \ge 0$, $ ma+nb  \in  H (f,\A)$ by (3.1); if $n = -l < 0$,
(3.2) with  $h = ma-lb$ and $k= lb$ gives $h  \in  H (f,\A)$,
$H(f,\A)$ is dense in $\r^+$.

\smallskip

The uniform continuity of $Pf$ for $f \in S^p_b$ with $p > 1$
follows with
    H\"{o}lder's inequality.

\smallskip

For (c) one can use [13, Theorem 4.1] :  $M_{h_0}(f) = : F \in \A$
implies
   $F  \in  W^{1,1}_{loc}$  and $ \Delta_{h_0} f = F'$ a.e.,  $\in \m\A$ by Lemma 2.2(c).
   this is an equation of the form (1.1) of  [13] with $n=0$ (admissible) and
   $m=2$, the determinant condition in Corollary 2.6 of [13] is here
  $-1 + \gamma_{h_0}$ $\not \equiv 0$, obviously true.

Theorem 4.1  of [13] can be applied with  $k=1$: the needed  $f
\in \m C_{ub}$ is equivalent with $Pu  \in  C_u $ by Proposition
1.7 (i);
   right side $= F'  \in \m^{1+1} \A$ holds since $F'  \in  \m\A$ and $\A  \st  \m\A$ for
   $\lambda$-classes by Corollary 3.3 below, implying  $\m\A  \st  \m^2\A$.

\smallskip

 $ AP(\r,X)$ is obviously a  $\lambda$-class, it satisfies $(L_{ub})$  by a
result of
   Kadets if $c_0$ is not isomorphic to a subspace of $X$ (see[13, p. 120],
  $(L.1)$ there $= (L_{ub}$) now).

 Instead of $(L_{ub})$ other assumptions are possible, see the end of
section 9.

 The proof shows that in (c) already $M_{h_0} (f)  \in  \m\A$  suffices
for
    $f  \in  \m\A$.

\smallskip

\smallskip

For a given class  $\A\st  L^1_ {loc}(\jj,X)$ one always has, with
Definitions (1.4)-(1.9),
 $\m^n\A \st \widetilde {\m}^n \A$, $n\in\N$.

\noindent  If moreover  $\A\st \m \A$, then $\A\st  L^1_
{loc}(\jj,X)$ and

\smallskip

(3.3) $\,\,\,\qquad \A\st\m\A\st \m^2\A\st \cdots\st \m^n \A \st
\cdots$

 $\,\,\,\qquad \, \qquad $ $\A\st\widetilde{ \m}\A \st\widetilde{ \m}^2 \A\st\cdots \st\widetilde{\m}^{n}\A \st \cdots $.

\smallskip

The assumption $\A\st \m\A$ does not always hold:
\proclaim{Example 3.1} There exists an infinite-dimensional
invariant closed linear subspace $\A_e\st C_{b}(\r,\cc)$ (with the
sup-norm), for which  $\A_e$ is even perpendicular to $\m \A_e$,
i.e.
  $\A_e
\cap\m\A_e= \{0\} =\A_e \cap C_u (\r,\cc) $:
\endproclaim

\noindent   If $ \A_e= $ smallest such subspace containing
$g(t)=e^{it^2}$, one has

(3.4) $\,\,\,\qquad \A_e = APg =\{ \phi g : \phi \in
AP(\r,\cc)\}$.

(3.5) $\,\,\,\qquad 0\not =f\in\A_e $, $h > 0$ implies $M_h f
(t)\to 0$ as $|t|\to\infty$, but  $f\not \to 0$.

\noindent Here (3.4) is clear; in $M_h f \to 0$ one can assume $f=
\gamma_{\om} g$ or $g_{\om}$, then $M_h f \to 0$ follows with
integration by parts; $f=\phi g \in C_0$ implies $\phi \to 0$,
then the Fourier coefficients $c_{\om} (\phi)=0$ for $\om \in \r$,
then $\phi \equiv 0$.

All this extends to $AP(\jj,X) g$, with arbitrary $\jj$, $X$.

 Similarly one gets $(\A_e *\h (\r,\cc) )\cap \A_e =\{0\}$,
so $\A_e *\h (\r,\cc) \not \st \A_e$.

\noindent A simpler (non-complete) example would be

$\A= C_c((0, \infty), X)$:
 $\qquad \m\A =\{ f=0 $ a.e.  on $ (0,\infty)\} $,

 $C_0 ((0,\infty),X)$, $\qquad O(e^{t^2})$ (see before Proposition  4.10),

 $X+\A_e$, $X+ P\A_e$ (see after Corollary 4.22).

\smallskip

\smallskip

Sufficient conditions for $\A\st \m\A$  are given by

\proclaim{Proposition  3.2}  Assume  $\A$  positive-invariant convex  $\st L^1_{loc}(\jj,X)$,   let   $L$
be the linear hull of  $\A$ and
$C (\jj,X) $ in $L^1_{loc} (\jj,X) $  and $||\cdot||:L\to [0,\infty]$ satisfy ($0\cdot\infty =0$) $||f+g||\le$

$ ||f||+||g||$, $||af||= |a|\cdot ||f||$, $||f_a||\le ||f||$  for $f, g\in L$, $0\le a\in \r$,
and  $||M_h \phi -\phi||\to  0$  as $0< h  \to 0 $, $\phi \in \A $;
 assume furthermore $\A$ closed in $L$  with respect  to $||\cdot||$.

\noindent Then $\A\st \m\A$.
\endproclaim

\demo {Proof} ($||\cdot||$ need not  be finite on  $\A$ or $L$): For $f\in\A$  and $h >0$ one has
$M_h f \in C(\jj,X) \st L\st L^1_{loc} (\jj,X) $. With $m\in\N$, $\delta:= h/m$,
$s_j := j\delta$, $0\le j\le m$, one gets

$||M_h f- \sum _{j=0}^{m-1} f_{s_j} \cdot (1/m)|| = ||1/h \sum _{j=0}^{m-1} (\int_{s_j}^{s_{j+1}}f_s \, ds - f_{s_j} \cdot \delta )|| =  $

$ || \sum _{j=0}^{m-1}\frac {1}{\delta m} (\int_{0}^{\delta}f_t \, dt - \delta f) _{s_j} || \le    (1/m ) \sum _{j=0}^{m-1}  ||M_{\delta} f- f||=  ||M_{\delta} f- f||  $.

\noindent Since  $ ||M_{\delta} f- f|| \to 0 $ as $m\to \infty $ and $\sum_{j=0}^{m-1} f_{s_j}/m \in \A$, one gets $M_h f\in \A$.
        \P

\enddemo
\proclaim{Corollary  3.3}  If  $\A \st C_{u}(\jj,X)$ is convex,
   positive-invariant and uniformly closed
 then $\A\st \m\A$.
\endproclaim

\proclaim{Examples  3.4}  Let $\, \A = X, \,$ $\,  P_{\tau},\,$
$\,  AP, \,$ $AAP, \,$ $EAP_{rc}$, $EAP, \,$ $UAA, \,$ $LAP_{ub},
\,$ $\T\E_{ub}, \,$ $ \E_{ub}, \,$ $C_{urc}, \,$ $C_{uwrc}, \,$
$C_{ub}, \,$ $C_{u}. \,$ Then for any $\jj$, $X$
 one has

 $\A\st \m\A \st \m^2\A \st \cdots $;

 \noindent for closed $\jj$ this holds also for $C_0$, $C_c$.

\noindent For these $\A$ also $\A* \h (\r ,\k)\st \A$ if $\jj=\r$ by (2.21). If  $\A=\A (\jj,\k) \st \m\A$,
 also $ W\A \st \m W\A \st \cdots  $.
\endproclaim

\proclaim{Example  3.5}  Proposition 3.2 can be applied to $\A = S^p AP $, $W^p AP $, $B^p AP $,
$1 \le p  <\infty$, any $ \jj$ and $X$, with  $ ||\cdot||_{S^p}$ respectively $ ||\cdot||_{W^p}$, respectively $ ||\cdot||_{B^p}$.
  For these also $\A* \h (\r ,\k)\st \A$ if $\jj=\r$.
\endproclaim
\noindent  Indeed, these $ ||\cdot||$ satisfy the conditions in Proposition 3.2,
 $ ||M_{h} \phi- \phi|| \to 0 $ as $h\to 0$ follows from the  approximation theorem for the $\A$'s using  $ ||M_{h} f||\le   c_p\cdot ||f||$
(see [50, p.  251, (42)]; see  also [21, pp. 71-78], [27,  pp.
34-36], [74, p. 14, Lemma 8]).

\noindent Furthermore,

(3.6) $\qquad \,\, AP\st \, S^p AP\st \, S^1 AP \st\, \m AP$
strictly if $ 1 < p < \infty$.

\noindent The three ``$\st$'' and the ''strictly'' for the first
two follows from the definitions.

An $f\in \m AP(\r,\r) \cap C^{\infty}(\r,\r)$, $-1 \le f \le 1$,
but $f \not \in S^1 AP (\r,\r)$, is given by

$f:= \sum_{n=2}^{\infty} \, h_n$, with $h_n \in
C^{\infty}(\r,[-1,1])$, $h_n \equiv 0$ on $[- 2 ^{n-1}, 2 ^{n-1}
-1 ]$,

 $h_n$ with $2 ^{n-2}$ oscillations on $[ 2 ^{n-1}-1, 2
^{n-1}]$ as $\sin (2 ^{n-1} \pi (t-2 ^{n-1}))$,

$h_n$ period $2 ^{n}$.

\noindent Since even $ f_1 -f \not \in S^1 AP(\r,\r)$, $g := f'
\in \m^2 AP (\r,\r)\cap C^{\infty}(\r,\r) \cap \h'_{pp}(\r,\r)$,
$g\not\in \, \m S^1 AP(\r,\r)$. See also [13, (3.8)].

\proclaim{Examples  3.6} Direct   calculations
   give    $\,\,\,
\A \st\m\A $ also for

 $\,\,\, \A = C_b (\jj,X)$,
 $\, L^{\infty} (\jj,X)$, $\,L^p (\jj,X)$, $\,S^p_{b} (\jj,X)$,
$\, \E (\jj,X)\,$, $\,\E_{0} (\jj,X)$, $\,\,\,$  any $\jj$ and
$\,X$, $\,\, 1\le p < \infty$, and $O(w)$ of Proposition 4.10 if
$w \in  \m O(w)$, e.g.  $w = t^r$ or $e^{rt}$ with real $r$.

\endproclaim
$L^p (\jj,X) \st \m L^p (\jj,X)$  follows with the continuous
Minkowski inequality [50, p. 251, Aufgabe 92 before (42)]. For
$\T\E_0$, $\T\E$ see Proposition 7.1.

\noindent The chains (3.3) can be finite:

\proclaim{Example  3.7}     For any  $n\in \N_0$  one has
      $ \A   \st  \m \A \cdots  \st \m^{n} \A =\m^{n+1 }\A= \cdots  $ if
 $\A= W^{1,n}_{loc} (\jj, X) $, and all inclusions here are strict (see also Corollary 4.6).

\endproclaim

\proclaim{Proposition  3.8}  For the real-linear positive-invariant  $\A \st C_u (\jj,X)$ assume
   $AP(\jj,X)\st \A   \st  \m \A$. Then all $``\st"$ in (3.3) are strict.
\endproclaim

\demo{Proof}  Since   $\widetilde {M}_h | L^1_{loc}(\jj,X)= M_h$
of (1.7), (1.4), one has
 for any   $ \A   \st   L^1_{loc}(\jj,X)$    or $\h'(\jj,X)$

\smallskip

 (3.7) $\,\,\,\qquad  \m^n (\A_{ Loc})=(\widetilde {\m}^n \A)_{ Loc}=(\widetilde {\m}^n (\A_{ Loc}))_{ Loc}$  for any $n\in\N$,

  $ \qquad\qquad\qquad\qquad\qquad\qquad U_{Loc}:= U\cap L^1_{loc}(\jj,X)$.
\smallskip

 \noindent Since   $  \m^{n-1} \A \st \m^n \A$, $\widetilde {\m}^{n-1} \A \st \widetilde{\m}^n \A$ by (3.3), we have to prove only
 $  \m^{n-1} \A \not= \m^n \A$, $n\in \N$. For this, choose pairwise  disjoint open intervals $I_n\st [0,1]$ of length
$|I_n| \le 2^{-n}$, then $g_n\in C^{\infty}(\r,\r) $ with period $2n+1 $  and $g_n\not=0$ in $[-n, n+1]$  only in $n+I_n$, with
$||g_n||_{\infty}=1$  and $\int_{n+I_n} g_n (t)\, dt=0$.
With $f_n :=Pg_n$  of period $2n+1$  and $f:= \sum_{n=1}^{\infty} f_n$ one has
$f\in C^{\infty}(\r,\r) $, $||f_n||_{\infty}\le |I_n| ||g_n||_{\infty} \le 2^{-n}$,
so $f\in AP(\r,\r) $, with $f'= \sum_{n=1}^{\infty} g_n$, $||f'||_{\infty}=1$.
$f'$ is not uniformly continuous on $[0,\infty)$ since $f'=g_n$ on $n+I_n$,
$|I_n|\to 0$, but sup$_{n+I_n} |g_n| =1$.

\noindent With $0\not = a\in X$,  $F := af |\, {\jj}\in
AP(\jj,X)\st \A$; so $F^{(n)}=af^{(n)}|\, {\jj}\in \m^n\A$ by
Lemma 2.2. $F^{(n)}\not\in \m^{n-1}\A$: Else with (1.11) one has
$h_1\cdots h_{n-1}af ^{(n)}*s_{h_{n-1}}*\cdots *s_{h_1}
=af'*(\Delta_{h_{n-1}}\delta )*\cdots *(\Delta_{h_{1}}\delta),$
$=G*\Delta_{h_{1}}\delta= \Delta_{h_{1}} G$ on $\jj$, with
 $ G=af'*(\Delta_{h_{n-1}}\delta*\cdots *\Delta_{h_{2}}\delta)|\jj$, and
$\Delta_{h_{1}}G\in\A \st C_u(\jj,X)$.

\noindent Since $f'$ is bounded on $\r$, $G$ is bounded on $\jj $
and $ \Delta_{h_{1}} G \in C_{ub}(\jj,X)$ for all $h_1 > 0$. By
Proposition 1.7 (i) $G$ is uniformly continuous on $\jj$.
Inductively one gets : $af'|\jj$ is uniformly continuous on $\jj$,
a contradiction.
        \P
\enddemo

\proclaim{Examples  3.9}  The conclusion of Proposition 3.8  holds
for $  \A = \,  AP $, $ \, AAP $,  $ \, EAP_{rc}$,  $ \,  EAP $, $
\, UAA $, $  \, LAP_{ub} $, $ \,  \T\E_{ub} $, $ \, \E_{ub} $,
$C_{urc} \,$, $C_{uwrc} \,$, $ \, C_{ub}$,   $ \, C_{u} $, any
$\jj$, $X$.

\noindent This  conclusion holds also for
 $\A = S^p AP (\jj,X)$, $1 \le p  < \infty$, since
 \endproclaim

 $  \m^{n-1} AP (\jj,X)\st \m^{n-1} S^p AP (\jj,X) \st \m^n AP (\jj,X)$, by [13, Corollary 3.6 respectively (3.8) hold also for $\jj\not =\r$],
$   \m^{n-1} S^p AP (\jj,X) = \m^n S^p AP (\jj,X)$ would imply
$  \m^{n}S^p AP (\jj,X)= \m^{n+1} S^p AP (\jj,X)$, and thus  $ \m^n AP (\jj,X)  =  \m^{n+1} AP (\jj,X)$.

\proclaim{Corollary  3.10} The conclusion of   Proposition 3.8
holds also  for $\W AP (\jj,X)$.
\endproclaim

\demo {Proof} If  $  f \in \W AP (\jj,X)$, $f$  is weakly
continuous on $\jj$; thus $ f\in L^1_{loc} (\jj,X)$ and $y (M_h
f)= M_h(y\circ f) \in  AP (\jj,\k)$  for all  $y\in  X^*$ (Example
3.4, [50, p. 52, Satz 3]); so $\W AP(\jj,X)\st\m \W AP(\jj,X)$,
(3.3) holds. If $ n  \in \N$, choose  $\phi\in \m^n  AP (\jj,\k)$
with  $\phi \not\in \m^{n-1}  AP (\jj,\k)$  by Proposition 3.7.
 Then if $0\not = a\in X$ and $f:= a\phi$, one has $f\in\m^{n} \W AP(\jj,X)$ but  $f\not  \in \m^{n-1} \W AP(\jj,X)$
because for any
$\phi :\jj\to\k$,  $0\not = a\in X$, $n\in \N_0$,  $\A(\jj,X)\st \jj^X $,  $\A(\jj,\k)\st \jj^{\k}$  with
  $y(\A(\jj,X))\st   \A(\jj,\k)$  and   $ X\cdot\A(\jj,\k)\st \A(\jj,X)$ for $y\in X^* $  one has

\smallskip

 (3.8) $\,\,\,\qquad a\phi\in  \m^n \A (\jj,X) \Leftrightarrow a\phi\in  \m^n \W\A (\jj,X) \Leftrightarrow \phi\in\m^n \A(\jj,\k)$.
        \P
\enddemo

\proclaim{Remark 3.11}  One can show that the codimensions of

$\m^{n-1}  \A$ in $  \m^n \A\qquad $
 and $\qquad\widetilde{\m}^{n-1}  \A$ in $ \widetilde {\m}^n \A $

\noindent  are infinite for all  $n\in \N$, $\A$ as in Examples
3.9 (except $S^p AP (\jj,X)$) or Corollary 3.10.
\endproclaim

\proclaim{Lemma 3.12}  If $U, V \st L^1_{loc} (\jj,X)$ satisfy
$\A\st \m\A$ the same holds for $U+V$, $U \cap V$.
\endproclaim

This follows immediately  from the definitions, for any $\jj$,
$X$.

\proclaim{Example 3.13}  $\A = VAP = V+AP$ satisfies $\A \st \m
\A$ if $V\st \m V$.
\endproclaim

\noindent So Zhang's $PAP $ [77-78] satisfies $ PAP \st\m PAP$,
for any $\jj$, $X$.

Here $PAP = V+AP$ with $V = C_b (\jj,X)\cap Av_n  (\jj,X)$;

\noindent with Lemma 3.12 and Example 3.6 one has to show $Av_n
\st \m Av_n$, where

\smallskip

$Av_n = \{ f \in L^1_{loc} (\jj,X): \frac{1}{T}
\int_{\al_0}^{\al_0 +T} ||f(s)||\, ds \to 0$ as $T\to \infty\}$;

\smallskip

\noindent but this is obvious with Fubini's  theorem.

Contrary to Example 3.1 with $e^{it^2}$ replaced by $\gamma_{\om}
(t)= e^{i\om \, t}$, one has

\proclaim{Remark 3.14}  If $\A \st L^1_{loc} (\jj,X)$, $\A$ has
$(\Delta)$, is  $\st   \m C_b $,  and  linear positive-invariant
uniformly closed, then $\gamma_{\omega} \A \st \m (\gamma_{\omega}
\A)$ for each $\omega\in \r$.
\endproclaim

\noindent We omit the proof, this is not used in this paper.

 \head{\bf \S 4   The $(\Delta)$ Condition}\endhead

The $(\Delta)$ condition introduced in Definition 1.4 was already
essential for the results of section 2. It is also important for a
characterization of the $\m^n\A$ $ (\S 5)$ and for Fourier
analysis $(\S 8)$, also for the "$\Delta$-spaces" which will be
introduced in a later paper. Here we show that all the $\A$'s
important in applications satisfy $(\Delta)$, usually this is
non-trivial.

A first application of $(\Delta)$ is

 \proclaim {Theorem 4.1} Assume  $\A$
linear, positive invariant,
 $\st$ $L^1_{loc} (\jj,X)$ with $(\Delta)$. If
 $ y, y^{m)}\in \A$, then $y^{(j)}\in \A$, $0< j< m$.
\endproclaim

This is an extension of the classical Esclangon-Landau result,
which was needed in the study of ap solutions of linear
differential equations (see Bohr-Neugebauer [29]) and which reads:
If $y \in C^n(\jj,\cc)$ is together with $Ly : = y^{(n)}  + \sum _
0 ^{ n-1} a_j y^{(j)}$ bounded on $\jj$ (with bounded $a_j$), then
the $y^{(j)}$ are bounded too for $0<j<n$ ([58, p. 177, Satz 1],
see also [51] and the references there). Theorem 4.1 extends this
from $C_b$ or $L^{\infty}$ (at least for $Ly = y^{(n)})$ to quite
general $\A$ by the results of this section and \S 7. Extensions
to solutions of linear differential-difference systems are
possible, this will be treated in \S 10.

\demo {Proof}  $m=1$: trivial.

$m\Rightarrow m+1$: $y^{(m+1)}\in A$, positive invariance and
linearity give $y_h^{(m+1)}- y^{(m+1)}\in A$ for all $h>0$; then
$(\Delta)$ gives $y^{(m+1)}- M_h y^{(m+1)} \in \A$ for all $h>0$.
This implies $M_h y^{(m+1)}= y_h^{(m)}-y^{(m)} \in \A$ for all
$h>0$.
  Set $z=z(h) = y_h-y $, $h> 0$. One
has $z, z^{(m)}\in \A$. The  assumptions of the induction imply
$z^{(j)}= (y_h-y)^{(j)}\in \A$, $0< j< m$, $h> 0$. By $(\Delta)$,
$y'-M_h y'= [y'- (y_h-y)]\in \A $, also if $m=1$. This implies
$y'\in \A$. One has $y', (y')^{(m)}\in \A$ and the induction
hypothesis implies $ y''\cdots y^{(m)}\in A$. \P
\enddemo

\proclaim {Proposition 4.2}  If  $\A \st C_{u}(\jj,X)$, condition
$(\Delta)$ holds for $\A$ provided one of the following conditions
is  satisfied:

 \noindent (i)   $\A \st C_{ub}(\jj,X)$, $\A$  $\Q$- convex and  uniformly closed.

\noindent (ii)   $\A \st C_{ub}(\jj,X)$, $\A$  is  a group under addition,  $A\st \m\A$ and $(L_{ub})$ holds for $\A$.

\noindent (iii)    $\A$  is  a group under addition,  $A\st \m\A$ and $(L_{u})$ holds for $\A$.

\endproclaim

\demo {Proof}
 (i) Let $\phi \in L^1_{loc} (\jj,X)$ and $\Delta_h\phi \in \A$,   for all $0< h\in \r$. Then
 $\phi \in C_u(\jj,X)$, by Proposition 1.7 (i).
It follows  $\phi-M_h\phi = \lim_{n\to \infty} (1/n)\sum_{k=1}^{n}( \phi -\phi_{s_k}) $ uniformly on $\jj$,
 where $s_k= (hk)/n$, $n\in \N$.
 By the assumptions on $\A$, $ (1/n)\sum_{k=1}^{n}( \phi -\phi_{s_k})\in \A $. Since  $\A$  is uniformly closed,
                 $\phi-M_h\phi \in\A$.

 (ii) Let $\phi \in L^1_{loc} (\jj,X)$ and $\Delta_h\phi \in \A$   for all $0< h\in \r$. Then
 $ \Delta_h (M_s\phi)= M_s(\Delta_h\phi)   \in \A$,  $\,0< h,\, s\in \r$ by  $A\st \m\A$.
 Since $\A$  is a group,$ \Delta_h (\phi- M_s\phi)=  \Delta_h \phi -  M_s(\Delta_h\phi)   \in \A$,  $\,0< h,\, s\in \r$.
 This implies $\phi- M_s\phi   \in \A$ by $(L_{ub})$ for all   $\,0<\, s\in \r$,  since $\phi- M_s\phi   \in C_{ub}(\jj,X)$ by (i) for $C_{ub}(\jj,X)$.

 (iii) Let $\phi \in L^1_{loc} (\jj,X)$ and $\Delta_s\phi \in \A \st C_{u}(\jj,X)$,   for all $0< s\in \r$. Then
 $ \Delta_h (\Delta_s\phi)   \in C_{ub}(\jj,X)$,  $\,0< h,\, s\in \r$, since   $\psi\in C_{u}(\jj,X)$ implies  $\Delta _h\psi\in C_{ub}(\jj,X)$ for $h> 0$.
So, by (i) for $C_{ub}(\jj,X)$ and  $\Delta _s\phi$ one gets
 $ \Delta_s (\phi- M_h\phi)=  \Delta_s \phi -  M_h (\Delta_s\phi)   \in C_{ub}(\jj,X)$,  $\,0< h,\, s\in \r$.
Proposition 1.7 (i) gives
  $ \psi := \phi- M_h\phi   \in C_{u}(\jj,X)$; by the assumption on $\A$ also   $\Delta_s\psi \in \A$, $s > 0$. The condition  $(L_u)$ gives $\psi \in \A$ for $h > 0$.
        \P
\enddemo

\proclaim {Proposition 4.3}  If  $\A \st \h'(\jj,X)$ satisfies
$(\Delta_1)$ [respectively  ($\Delta$)],  then  $\A$ satisfies
$(\Delta'_1)$ [respectively  $(\Delta')$]  provided  one of the
following conditions is satisfied:

 \noindent (i)   $\A \st C(\jj,X)$,

\noindent (ii)   $\A \st L^1_{loc}(\jj,X)$, $ f\in \A$ and $g=f$  a.e on $\jj$ implies   $g\in \A$,

$\qquad\qquad$ and $ X$ has the Radon-Nikodym  Property
  [35, pp. 181-182].
\endproclaim

For the proof we need an extension of Proposition 1.5:

\proclaim {Lemma 4.4}   If for   $T \in \h'(\jj,X)$,  $\e_0 > 0$,
$k\in \N_0 $ one has $  \Delta_h  T\in  C^{k} (\jj, X)$
  for  all $h\in (0, \e_0]$, then  $T \in C^k (\jj, X)$. If only  $  \Delta_h  T\in  W^{1,k}_{loc} (\jj, X)$
  for  all $h\in (0, \e_0]$ and $X$ has the Radon-Nikodym property, then  $T \in W^{1,k}_{loc} (\jj, X)$.
\endproclaim

\noindent Special case:  $  \phi\in  L^{1} _{loc} (\jj, X)$, to each
   $h\in (0, \e_0]$ exists $g\in  C^{k} (\jj, X)$ with $  \Delta_h  \phi = g$ a. e. on $\jj$ implies the existence of
an   $f \in C^k (\jj, X)$ with  $\phi =f$ a.e on $\jj$.

\demo {Proof}    The $C^k$-case follows with Proposition 1.5 as in
the proof of Lemma 14 [49, p. 226], case $\A= C^k$ since also for
general $X$ to  $T \in \h'(\jj,X)$ and $K$-compact $\st \jj$ there
exists   $g \in C (\r ,X)$  and
   $n \in  \N_0$ with   $T = g^{(n)}$ in an open neighborhood of $K$  (Theorem 2.11, $a\Rightarrow b$ for $\A = C_{ub}$).

Case     $ W^{1,k}_{loc} (\jj, X)$ with arbitrary $\jj$:  Assume
first $\phi \in C (\jj, X)$ with $ \Delta_h  \phi$ absolutely
continuous in $\jj$
  for  all $h\in (0, \e_0]$.  For $K= [a,b]$   compact $\st\jj$ there is
  $c \in (0, \e_0]$ with $K_1:= [a,b+c] \st \jj$.
For fixed $\e >0 $, Lemma 4.12 below can be applied to

 $I= K_1$,
$f=\phi |\, I$, $\gamma =c$,

      $\A=  AC_{\e} (K, X): = \{g\in  C (K, X): $ to $g$ exists $ \eta > 0$  with

      $ \sum_{i=1}^{m} ||g(t_i)-g(s_i) ||\le \e \, \,\qquad$ if

 $ \,\,\,\qquad\qquad \sum_{i=1}^{m} |t_i-s_i |\le \eta, a\le s_1 < t_1\le s_2< t_2 \cdots <t_m \le b, k\in\N \}$,

\noindent  with $\A_n :=\A$  with $\eta $ replaced by the fixed $1/n$ and  `` to $g$ exists $ \eta > 0$ with" deleted.
So there are $v, \delta$ with  $0\le v < v+\delta \le c$ and $(\phi- M_c (\phi)_v)| K \in \A$.
Since $\phi$ is continuous on $\jj$, $ M_{\delta} (\phi)_v  $ is absolutely continuous on $K$, and so $   \phi |\, K \in
 AC_{\e} (K, X)$. This holds for any $\e >0$, so  $\phi$  is absolutely continuous on $K$. $K$ being arbitrary, $\phi$
 is locally absolutely continuous  in $\jj$ for any $\jj$, $X$.  So, one has  with Proposition 1.5

\smallskip

 (4.1) $\qquad  \phi \in L^1 _{loc} (\jj, X)$,  $\Delta_h \phi\in AC_{loc} (\jj, X)$ for   $ 0 < h <\e_0 $

    $\qquad \qquad\qquad$ implies  $\phi \in  AC_{loc} (K, X)$.
\smallskip

\noindent If additionally $X$ has the Radon-Nikodym property,
$\phi'$  exists a.e in $\jj$ by [35, p. 107, Theorem 2, p. 138,
Corollary 8], implying  $\phi\in   W^{1, 1}_{loc} (\jj, X)$.

\noindent A simple induction and (2.5) give      $\phi\in   W^{1, k}_{loc} (\jj, X)$ if  $\phi \in  C ( \jj, X) $
with  $\Delta_h\phi\in   W^{1, k}_{loc} (\jj, X)$  for all  $ 0 < h <\e_0 $,
      $ k\in  \N_0 $.

\noindent  If now    $T \in \h'(\jj,X)$ with   $ \Delta_h  T\in
W^{1,k}_{loc} (\jj, X)$   for all  $ 0 < h <\e_0 $ and fixed $k
>0$,  $  \Delta_h T \in C^{k-1}(\jj,X)$, so the first part of
Lemma 4.4 gives  $   \phi \in C^{k-1}(\jj,X)$ with $T= \phi$ on
$\jj$. One has
  $  \Delta_h \phi \in W^{1,k}_{loc} (\jj, X)$   for all  $ 0 < h <\e_0 $, the above gives
 $ T= \phi \in W^{1,k}_{loc} (\jj, X)$.

Case $k=0$,  $ W^{1,0}_{loc} (\jj, X) :=L^1_{loc } (\jj, X)$:  To  $T \in \h'(\jj,X)$    exists  $S \in \h'(\jj,X)$
 with $S' = T$  on $\jj$ (see [71,  p. 51, The\'{o}re\`{m}e I for $X=\cc$]).
Then  $  \Delta_h S\in   W^{1,1}_{loc} (\jj, X)$ by (2.5),  $ S\in W^{1,k}_{loc} (\jj, X)$ by the above,
 $ T= S'\in  L^{1}_{loc} (\jj, X)$   by (2.5).
        \P
\enddemo
\demo {Proof of Proposition 4.3} (i) : If   $T \in \h'(\jj,X)$
with $  \Delta_h T\in \A$ for $h >0$, $ T=f \in  C (\jj, X)$ by
Lemma 4.4, to $h > 0$ exists  $g= g_{(h)}\in \A$ with  $ \Delta_h
f =g $ a. e. ; since  $f, g$ are continuous,  $  \Delta_h f\in
\A$. Then, by $( \Delta)$,  $  f-M_hf \in  \A$, $T-\widetilde
{M}_h T = f-M_h f\in \A$.

(ii)  $T \in \h'(\jj,X)$,  $  \Delta_h T\in \A\st   L^{1}_{loc}
(\jj, X)$  imply $T=f  \in   L^{1}_{loc} (\jj, X)$ by Lemma 4.4;
furthermore  $  \Delta_h f\in \A$ a.e., so $  \Delta_h f\in \A$ by
assumption.  $  (\Delta)$  gives  then  $  (\Delta')$. Similarly
for  $  (\Delta'_1)$.
        \P
\enddemo

\proclaim {Corollary 4.5}

$\m C^k(\jj,X) = C^{k-1}(\jj,X) + \{ f \in X^{\jj} : f=0 \,$ a.e.
$\}$
                if $\,\, k \in   \N$,

\smallskip

                  $\m C(\jj,X) = L^1_{loc}(\jj,X)$, $\qquad\,\,\,$ any $\,\, \jj, X$.
\endproclaim

\noindent       This follows from the "special case" after Lemma
4.4.

\proclaim    {Corollary 4.6}

 $\m W^{1,k}_{loc}(\jj,X) =
W^{1,k-1}_{loc}(\jj,X) $ if $k \in   \N$,

 $\m L^1_{loc}(\jj,X) = L^1_{loc}(\jj,X)$.
\endproclaim

\noindent See also Example 3.7.

\proclaim{Examples  4.7} By Propositions 4.2, 4.3, 1.5, the
following $\A$  all satisfy   $  (\Delta)$  and  $  (\Delta')$,
for arbitrary $\jj$ and $X$:

Constant functions $X$, Continuous periodic with period $\tau$
functions $P_{\tau}$, $ AP \,$, $AAP \,$, $ EAP_{rc} \, $, $ EAP
\, $, $UAA \,$, $LAP_{ub} \,$, $\T\E_{ub} \,$,  $ \E_{ub} \,$,
$C_{urc} \,$, $C_{uwrc} \,$, $C_{ub} \,$, $C_{u}, \, C_{0} $,
$C^k$ for $k\in\N_0$, $\, C^{\infty}$.

$S^p AP (\jj,X)$ has  $  (\Delta)$  for  $1\le p <\infty$,  and
$(\Delta')$ for $X$ with the Radon-Nikodym property (use the
Bochner transform of [13, p. 132, (3.1),(3.2)]).

\noindent Furthermore $\W \A (\jj,X)$ has   $ (\Delta)$  and
$(\Delta')$ for $X$ if $ \A=\A (\jj,\k)$ has  $  (\Delta)$.

\endproclaim

\proclaim       {Example 4.8} $ \A= C_c(\jj,X)$, $\h (\jj,X)$ have
$(\Delta)$, any $\jj, X$.
\endproclaim
\demo{Proof} We prove the case $\jj \not =\r$, $\jj \, =\r$
follows similarly. Let $\phi \in L^1_{loc}(\jj,X)$ and $\Delta_h
\phi \in C_c(\jj,X)$ for all $h
>0$. By Proposition 1.7 (i), $\phi \in C_u (\jj,X)$. By the
assumptions  supp $\Delta_1 \phi \st [a,b]$ for some $a, b \in
\jj$. It follows $\phi (t+1) = \phi (t)$ for all $t\ge b$. The
uniform continuity of $\phi$ and $\Delta_h \phi = 0$ far out for
       small $h$  implies then $\phi(t) = \phi(b)$ for  $t>b$.
        It follows $\phi - M_h \phi \in C_c$ for closed $\jj$.

     If  $\jj = (\al ,\infty)$ one can even show  $\phi =$ constant on $\jj$,
$(\Delta)$ follows.

 $(\Delta)$ for $\h(\jj,X)$ follows from  $(\Delta)$ for $C_c$ and
$C^{\infty}$ of Example 4.7. \P
\enddemo

\proclaim{Proposition 4.9} $L^{p}_w (\jj,X) $ satisfies $(\Delta)$
for $1\le p \le \infty$, arbitrary   $ \jj, X$, and measurable
$w:\jj \to [0,\infty)$ for which there is a real c with

\smallskip

 (4.2) $\qquad w(t+h) \le c \, w(t)$ and  $\qquad w(t) \le
c \, w(t+h)$ a.e. in $t\in \jj$, $0 < h \le 1$.
 \endproclaim
\smallskip

\noindent Here $\qquad L^{p}_w (\jj,X) := \{ f$ measurable  $\in
X^{\jj}:\qquad w |f|^p \in L^{1} (\jj,\r)\}$

\noindent with

\smallskip

$||f||_{p,w} = ||\om\,f||_p = [\int _{\jj} w(t) ||f(t)||^p\,
dt]^{1/p}$ if $1\le p < \infty$,

\smallskip

$L^{\infty}_w (\jj,X) := \{ f $ measurable $\in X^{\jj}: \qquad w
|f|$ measurable bounded a.e.  $\}$

\noindent with

\smallskip

$||f||_{\infty ,w} = ||\om \,f||_{\infty} = $  ess\, sup $_{t\in
\jj} w(t)\, ||f(t)||$,

\smallskip

\noindent  with $\qquad \om := w^{1/p}\qquad $ if $p <\infty$,
else $w $.

\smallskip

\noindent $Examples$ of $w$ with (4.2) are $1+|t|^r$ or $e^{r\,
t}$ for any real $r$ and $\jj$ or $t^r$ if $0\not \in $ closure of
$\jj$.

\demo{Proof}
 Let $\,\phi\,\in\, L^{1}_{loc}(\jj,X)$ and $\,\Delta_h\phi\, \in
\, L^{p}_w(\jj,X)$, $1\le p $ fixed $ <\infty$. The measurability
of $\, \phi\,$ gives the $(t,h)$- measurability of $\,f:= \om
\Delta_h\phi$ on $\,\, \jj\times \r^+$. By assumption  $g(h) :=
\,\,||\om \Delta_h\phi||_p < \infty$ for $h\in \r ^+$; $g:  \r^+
\to [0,\infty)$ is measurable by a suitable version of the
continuous Minkovski inequality  [50, Aufgabe 92, \S7, p. 251].
   So there is $n\in\N$ such that $K_n := \{ h\in (0,1]: g(h) \le n\}$ has positive measure.
 Then the difference set
$K_n-K_n$  contains $(-\e,\e)$ with $0< \e \le 1$ (see [50, p.
189, Aufgabe 116, \S6]): To $0 < h < \e$ there are $u$, $v \in
K_n$ with $h= v-u$ and $g(u) \le n $, $g(v)\le n$, or $||\om
(\phi_u- \phi_v)||_p\le 2n$. With (4.2) and $\rho := c^{1/p}$ one
gets

$||\om (\phi_ {h+1}- \phi_1)||_p\le  \rho ||(\om (\phi_ {h}-
\phi))_1||_p \le$

$ \rho^2 ||(\om (\phi_ {h}- \phi))_u||_p \le  \rho^2 ||\om (\phi_
{v}- \phi_u)||_p \le 2\rho^2\, n$.

\noindent Furthermore $g(1) = ||\om (\phi_ {1}- \phi)||_p =: q <
\infty$, this and (4.2)implies similarly

$||\om (\phi_ {h +1}- \phi_h)||_p \le \rho \, q $.

\noindent Assuming $c\ge 1$, these three inequalities yield

$g(h) = ||\om \Delta_h \phi||_p \le 2 \rho^2  (n+q) =: \sigma$ for
$0< h \le \e$.

\noindent Since $||\om (\phi_ {(n+1)h }- \phi)||_p \le ||(\om
(\phi_ {h }- \phi))_{nh}||_p + ||\om (\phi_ {nh }- \phi)||_p \le
\rho^n \sigma + ||\om (\phi_ {nh }- \phi)||_p$, induction gives
the boundedness of $g $ on each $[0,T]$, $T > 0$.

The continuous Minkovski inequality gives  with  $ \psi (t): =
\int _0 ^T  \om (t )\Delta_s \phi (t) \,ds$ that $ \psi \in L^p
(\jj,X)$, or (with the measurability of $\phi$) $M_T \phi - \phi
\, \in L^p_w (\jj,X)$ for $T >0$, i.e. $(\Delta)$ holds for $
L^p_w $. All this works also for $p=\infty$. \P

\enddemo

\noindent For any $w: \jj \to [0,\infty)$ define $\qquad
O(w)(\jj,X):=\,\, \{ f\in L^1_{loc} (\jj,X):$

$\qquad$ to $f $ exists $\,\, m\in \N$ with $||f(t)|| \le m w(t)$
a.e. $\,\,$ for $|t| \ge m, \, t\in \jj\}$.

\proclaim{Proposition 4.10} $O(w) (\jj,X) $ satisfies $(\Delta)$,
 for any   $ \jj, X$ and measurable $w:\jj
\to [0,\infty)$.

   \endproclaim
\demo{Proof}
 Let $\,\phi\,\in\, L^{1}_{loc}(\jj,X)$ and $\,\Delta_h\phi\, \in
\, O(w)(\jj,X)$ for all $h> 0$. We treat the case $\jj \not =\r$,
$\jj \, =\r$ follows in the same way.

For $m\ge m_0 \in \N \cap \jj$ define  $\psi (t,h):= m \,w(t)-
||\Delta_h \phi (t)||$ and $q_m (h):= \mu _L -$ inf $\{ \psi
(t,h): \qquad t\ge m\} $, where $ \mu _L =$ Lebesgue measure. Then
$q_m : [0,\infty) \to [\infty,-\infty) $ is well defined, it is
measurable since $||\Delta_h \phi||$ is $(t,h)$ measurable (see
[50, p. 140, Aufgabe 92]).

So $\Omega_m := \{ h\in (0,1]: q_m (h) \ge 0 \}$ is measurable for
$m\ge m_0$. Since $[0,1]= \cup _{n=m_0}^{\infty}\Omega_n$, there
is $m_1$  with $\mu_L(\Omega_{m_1}) > 0$, then as in the proof of
Proposition 4.9 there is $\e \in (0,1]$ with $[0,\e] \st
\Omega_{m_1}-\Omega_{m_1}$.

To $h \in (0,\e]$ there are thus $u, v\in \Omega_{m_1}$ with
$h=v-u$, implying $F(t,h):= 2 m_1 \,w(t)- ||\Delta_h \phi (t)||\ge
0 $ if $t \ge m_1 +1$, $t\not \in$ nullset (depending on $h$).

Now $M_T := \{ F(t,h)< 0: m_1 +1 \le t \le T, 0\le h\le \e  \}$ is
measurable, all $M_T|_h$ are null sets by the above, $0\le h\le \e
$;  with Fubini's theorem, $M_T$ is a null set, so there is a null
set $P_1 \st [m_1 +1, \infty)$ such that  for $t\not \in P_1$,
$t\ge m_1+1$ one has $ F(t,h) \ge 0$ for almost all $h\in [0,\e]$.
As in the proof of Proposition 4.9 one gets from that a null set
$P \st [m_1 +1, \infty)$ such that to $a >0$ there is $m\in \N $
with $||\Delta_h \phi (t)|| \le m w(t)$ for almost all $h \in
[0,a]$, $t\ge m$, $t \not \in P$. Integrating with respect to $h$
one gets

$||\phi - M_a \phi|| \le m w(t)$ for almost all $t\ge m$, that is
$\phi - M_a \phi \in O(w)(\jj,X)$ for all $a >0$. \P
\enddemo

\proclaim{Remark  4.11}  For any $\jj$, $X$,
 $C_{rc}(\jj,X)$ has $(\Delta)$.
\endproclaim

\noindent We omit the proof, this is not used here. See section
10, question 5 for $C_{wrc}$ and $L^{\infty}_{wrc}$.

To get  $  (\Delta)$ for additional spaces, and already for Lemma
4.4 and Proposition 4.3,   we need

\proclaim {Lemma 4.12} Assume $I\st \r$ arbitrary interval, $0
<\gamma < |I|$,  $f \in  L^1_{loc}(I,X)$. Assume further  $\A=
\cup_{n=1}^{\infty} \A_n \st  L^1_{loc}(I^{-\gamma},X) $  with the
following properties:

 \noindent (i)   $ (\Delta_s f)| I^{-\gamma} \in \A$ for all $0 \le s \le \gamma $,

\noindent (ii) If   $ (\Delta_{s_m} f)| I^{-\gamma} \in$ same
$\A_n$ for  $s_m\in [0, \gamma]$  and $s_m\to r$, then  $
(\Delta_{r} f)| I^{-\gamma} \in \A_n$.

\noindent (iii) If   $ (\Delta_{s} f)| I^{-\gamma} \in \A_n$ for
$ 0\le u \le s \le u+\rho \le \gamma$  with  $\rho >0 $, then

 $ -( (1/\rho) \int_{u}^{u+\rho}\Delta_{s} f \, ds)| I^{-\gamma} \in \A$.

\noindent Then there exist $v, \delta$ with  $ 0\le v < v+\delta \le \gamma$ such that   $ (f-(M_h f)_v) |  I^{-\gamma} \in \A$  for $0 < h \le\delta$.

\endproclaim

\noindent Here $ I^{-\gamma} = (\alpha, \beta - \gamma)$ [respectively  $(\alpha, \beta - \gamma]$]
 if $I =(\alpha, \beta )$ [respectively   $(\alpha, \beta ]$],
 $-\infty \le  \alpha < \beta   < \infty$, similarly for   $I =[\alpha, \beta )$, $[\alpha, \beta ]$ and
  $  I^{-\gamma}= I$  if $\beta = \infty$.

\demo {Proof} Define $M_n: = \{ s\in [0,\gamma ] :   (\Delta_{s}
f)| I^{-\gamma} \in  \A_n \}  $. Then by (i), $[0,\gamma] =
\cup_{n=1}^{\infty} M_n$. By (ii), the $M_n$ are closed in
$[0,\gamma]$. So by the  Baire category theorem [75, p. 12] there
exist
 $ n_0, v, \delta $  with   $ 0\le v < v+\delta \le \gamma $ and $ [v , v+\delta ] \st  M_{n_0}$.
By (iii), the integral expression there is in $\A$ for $u=v$ and
 $ 0 < \rho \le \delta $.  Then formula

\smallskip

(4.3) $ \qquad\qquad\,   f(t) = (M_{\delta} f)(t+v)  -( (1/\delta)
\int_{v}^{v+\delta}\Delta_{s} f (t) \, ds)$ for $ t  \in
I^{-(v+\delta)}$

\smallskip

 \noindent gives $ (f-(M_h f)_v) |  I^{-\gamma} \in \A$  for $0 < h \le \delta$.
        \P
\enddemo

\proclaim{Remark 4.13} If in (iii) the $\A$ can be replaced by the
$\A_n$ of the assumption of (iii), then also the $\A$ in the
conclusion of Lemma 4.12 can be replaced by
                         $ \A_{m_0}$  with $m_o$ independent of $h\le \delta$.
\endproclaim

\proclaim{Corollary  4.14}  If  $\A$,  $\A_n$ are as  in Lemma 4.4
but with $I=\jj $ (so  $I^{-\gamma}=\jj)  $, and if additionally
$\A$ is  real-linear and $\A \st \m\A$,
 then  $\A$  satisfies $  (\Delta)$.

\endproclaim

\demo {Proof} Since  $\jj = (\alpha,\infty ) $,  $[\alpha,\infty )
$  or $\r$,  $  I^{-\gamma}=\jj  $. By Lemma 4.4, if $f\in
L^1_{loc} (\jj,X)$ with $\Delta_s f\in \A$ for $s >0$, one has $
(f-(M_h f)_v )\in \A  $          for $0 < h \le\delta$ with
suitable  $ \delta, v$. Now for $0 < h \le \delta$

 $ (f-M_h f) = f +M_h (\Delta_v f) -M_h (f_v) =  f  -(M_h f)_v +M_h (\Delta_v f), \in \A  $

\noindent  by the  assumptions on $\A$. Furthermore

 $ 2 (f-M_{2h} f) = 2f - M_h f -M_h (f_h) =  f  -M_h f + f-M_h f- M_h (\Delta_h f),  \in \A  $.

\noindent This gives $ (f-M_ {2^n h} f)   \in \A$  for $0 < h \le \delta$, $n\in \N$, and thus $(\Delta)$ for $\A$.
        \P
\enddemo

\noindent Special case: If $\A \st C(\jj,X)$, (ii) of  Lemma 4.12
holds  if  the $\A_n$ are  closed with respect to pointwise
convergence on $\jj$ (Proposition 1.5).

\smallskip

                In the following

   $Lip (\jj,X):=\{f \in  C (\jj,X): $ to $  f$ exists $L $ with $ ||f_h-f||_{\infty} \le L |h|$ for all $ h > 0 \} $.

\smallskip

    $A C (\jj,X)$ is the space of absolutely continuous functions.



\smallskip

\proclaim{Examples  4.15}
  $\A = C_b (\jj,X)$, $A C (\jj,X)$
   and $Lip (\jj,X)$     satisfy   $  (\Delta)$  and  $  (\Delta')$.
\endproclaim

\noindent Here  $Lip (\jj,X) =\cup _{n=1}^{\infty} \A_n$  with
$\A_n := \{ f \in C (\jj, X) : |\Delta _h f| \le nh $
 on $\jj$ for all $ h>0  \}$, similarly  for $C_b (\jj,X)$.

\noindent      (ii) and (iii) are always fulfilled (Proposition
1.5), with the integral    in  (iii) of Lemma 4.12 even in $A_n$.

\noindent $AC$ can be treated with $A_n^{\e}$ as in the proof of
Proposition 7.3.

\smallskip

\proclaim {Example 4.16}
  $S^p_b  (\jj,X)$ has  $  (\Delta)$  for  $1\le p <\infty$.
\endproclaim

\noindent  This can be reduced   to $  (\Delta)$  for $C_b(\jj,X)$
of Example 4.15 with the Bochner-transform of [13, p. 132, (3.3)
and Lemma 3.2]. A direct proof follows from Lemma 4.11, Corollary
4.13 with

$A_n := \{ g\in S^p_b : ||g||_{S^p} \le n \}$; by approximating
$f$ on $[t, t+2+r]$ in the $L^p$-norm with  continuous functions,
$t$ fixed, one gets (ii), (iii) follows with the  continuous
Minkovski inequality as in the proof of Proposition 4.9.

\proclaim{Proposition 4.17} (a) For any $\jj$, $X$, $n\in\N_0$, if
$U\,$ and $V\st X^{\jj}$  have $(\Delta)$, then also $U \cap V$,
$\m^n U$.

(b) If $U$ is real linear, positive invariant,  $\st
L^1_{loc}(\jj,X)$ with $(\Delta)$, then $U\st \m U$.
\endproclaim

\demo{Proof} (a) Follows immediately from the definitions  and
$M_h M_k \phi =M_k M_h \phi$ of (2.6).

(b) If $f\in U$, then  $\Delta_h f\in U$, thus $f-M_h f\in U$,
$f\in \m U$.
\enddemo

\noindent Examples: $ L^p \cap C_0$, $AP \cap C^k$; but
$O(e^{t^2})\not \st \m(O(e^{t^2}))$, though it has $(\Delta)$ by
Proposition 4.10

\noindent The case $U+V$ is not so simple:

 \proclaim{Proposition 4.18} Assume $U\,$ and $V$ are
additive groups, positive-invariant, with  $(\Delta)$, $\st
L^1_{loc}(\jj,X)$; assume further $U \cap V=\{0\}$, $\jj =[\al,
\infty)$ or $\r$  and for $\Phi \in L^1_{loc}(\jj,X)$ with
$\Delta_h \Phi \in U+V$ for $h>0$ in the decomposition $\Delta_h
\Phi (t)= u(t,h) +v(t,h) $ the $v(\al,h) $ in $h\in (0,\infty)$ is
measurable; finally, if $\jj \, =\r$, $\,\al =0$ and additionally
$U$, $V$ are  invariant. Then  $ U+V$ satisfies $(\Delta)$.
\endproclaim

\demo{Proof} If $w=u+v$ with $u\in U,\,$ and $v\in V$, the $u,\,
v$ are unique, so for $\Phi$ as in the assumption the $u(t,h),
v(t,h)$ are well defined for $t\in \jj$, $h >0$. Now
$\Delta_{h+k}\Phi =(\Delta_{h}\Phi )_k +\Delta_{k}\Phi $ give with
the $U$, $V$ positive-invariant for $t\in \jj$, $h$, $k\, \ge 0$

\smallskip

(4.4)$\qquad u(t,h+k)= u(t+k,h)+ u(t,k), \,\, v(t,h+k)= v(t+k,h)+
v(t,k)$.

\smallskip

\noindent Case $\jj= [\al,\infty)$: $t=\al$ and $\phi (s):= v(\al,
s-\al)$, $s \ge \al$ gives for $k=s-\al$

$\Delta_h \phi (s) = v(s, h)$, $s\in\jj$, $h \ge 0$;

\noindent since $ v(\cdot ,\, h)\in V$ for $h \ge 0$ and $V$ is
positive-invariant, this gives $\Delta_h \phi\in V$ for $h > 0$.

By assumption $v(\al ,\,\cdot) $  is measurable on $\r^+$, then
$\phi$ on $\jj$, with $\Delta_h  \phi \in V \st L^1_{loc}(\jj,X)$.
 Lemma 4.19 below gives $\phi \in  L^1_{loc}(\jj,X)$.

Since $V $ satisfies $(\Delta)$ one has $\phi- M_h \phi \in V$ for
$h >0$. With $\psi :=\Phi -\phi$, this implies $ \Delta_h\psi =
\Delta_h \Phi - \Delta_h\phi =u(t,h)+ v(t,h) - (v(\al, t+h-\al)-
v(\al, t-\al)), = u(t,h)$ by (4.4), for $t\in \jj$, $h> 0$.

\newpage

So $\Delta_h \psi\in U$ for $h > 0$ with  $\psi \in
L^1_{loc}(\jj,X)$. $(\Delta)$ for $U$ gives  $\psi- M_h \psi \in
U$, with the above one gets $\Phi- M_h \Phi \in U+V$.

Case $\jj=\r$, $\al=0$ (for example): $\Delta_{-h} \Phi=-
(\Delta_h \Phi)_{-h}$ and the full invariance of $U$, $V$ yield
 $ \Delta_{h}\Phi \in U+V$ for all real $h$. Then the above $\phi$
 is defined on $\r$ with $\Delta_h \phi\in V$ for $h > 0$;
  this gives, with the assumption  $v(0,h)= \phi(h)$ measurable on $ \{0\}\times \r^+$ and
  $V\st L^1_{loc} (\r,X)$, measurability of  $\phi$ on $\r$. One can
 proceed as in the case $\jj= [\al,\infty)$, now with $\jj=\r$,
 getting $(\Delta)$ for $U+V$. \P
\enddemo

 \proclaim{Lemma 4.19}  Assume $\phi :\jj \to X$ is
measurable,   $\jj$, $X$ arbitrary and $\Delta_h \phi\in
L^1_{loc}(\jj,X)$ for all $0 <h \le \e_0$. Then $\phi\in
L^1_{loc}(\jj,X)$.
 \endproclaim

\demo{Proof}(See Remark 1.6) With $[a,b]\st \jj$, $\psi(h): =\int
_a^b || \Delta_h \phi (t)||\, dt$ is measurable in $h\in [0,\e]\st
[0,\e_0]$ by the Fubini-Tonelli theorem. As in the proof of
Proposition 4.9 there exists $n_0 \in \N$ and $\delta >0$ such
that to $0\le h \le \delta$ there are $u, v \in [0,\e]$ with
$h=v-u$ and $ \int _a^b ||  \phi_v (t)- \phi_u (t)||\, dt \le 2
n_0$, or $ \int _{a+\e}^b ||  \phi_h (t)- \phi (t)||\, dt \le 2
n_0$ for $0\le h\le \delta$. The Fubini-Tonelli theorem gives then
the existence of $ \int _{0}^{\delta}   \phi (t+h)\, dh$ for
almost all $t \in [a+ \e, b]$; this implies $\phi\in
L^1_{loc}([a+2\e, b],X)$. $b\to \infty$, $\e \to 0$, $a \to$ end
point of $\jj$ respectively $\al \to -\infty$ gives $\phi\in
L^1_{loc}(\jj,X)$, if $\al \in \jj$, with $\Delta_{\e_0}\phi \in
L^1_{loc}(\jj,X)$. \P
\enddemo

 \proclaim{Corollary 4.20} If $\A\st
L^1_{loc}(\jj,X)$ has  $(\Delta)$, and if $\phi :\jj \to X$ is
only measurable with $\Delta_h \phi\in \A$ for all $h
>0$, then $\phi -M_h \phi\in \A$ for $h >0$.
 \endproclaim

\noindent For  applications of Proposition 4.18 see section 7.

Next we show that $(\Delta)$  does not always hold.

 \proclaim{Example 4.21}
 $ \A_e =AP(\jj,X). e^{it^2}$  of Example 3.1   does not satisfy $(\Delta)$ of
Definition 1.4, though it is  a linear uniformly  closed invariant
subspace  of $C_b(\jj,\cc)$.
    \endproclaim
\demo{Proof} Assuming the contrary, since $ \A_e $ is   linear and
invariant, $\A_e \st \m \A_e$ by Proposition 4. 17 (b). This
contradicts Example 3.1. \P
\enddemo
\noindent Example 4.21 implies, (even for arbitrary $\jj$, $X$),
with

$PU= \{ Pf:\,\, f\in U\}\,\,$
 indefinite integrals.

 \proclaim{Corollary 4.22} $\A _1$, $\A _2$ do  also not satisfy
 $(\Delta)$ and $\A \st\m \A$, where

(4.5) $\qquad\, \A _1 = X+ AP\cdot e^{it^2}$,$\,\,\qquad\, \A _2 =
X+ P(AP\cdot e^{it^2})$.
 \endproclaim
\demo{Proof} Since both $\A_1$, $\A_2$ are linear and positive
invariant, $(\Delta)$  would imply $ \A_i\st \m\A_i$ by
Proposition 4.17 (b), $i=1, 2$.

$ \A_2\st \m\A_2$ and  $M_h P\psi- PM_h \psi=  M_h P\psi (\al_0)$,
$\psi  \in L^1_{loc} (\jj,X)$, $h>0$ imply,

\newpage

for $f\in AP\cdot e^{it^2}$ and $h>0$, $PM_h f = a +P\phi$ with
$a\in X $ and $\phi \in AP\cdot e^{it^2}$; differentiation yields
$M_hf= \phi $,  or $f \in \m (AP\cdot e^{it^2})$, contradicting
Example 3.1 if $f\not \equiv 0$.

$ \A_1\st \m\A_1$ would imply $\phi := M_h  e^{it^2}= a+\psi\cdot
e^{it^2}$ with $a\in X $ and $\psi \in AP$, and $\phi \in C_0$ by
Example 3.1, then $0= m (\phi )= a+ m(\psi\cdot e^{it^2})=a$,
$\psi\cdot e^{it^2} \in C_0$, $\psi=0$ by Example 3.1; so $M_h
e^{it^2}=0$ for $h>0$, a contradiction. \P
\enddemo

 Conversely one has, using Proposition 1.5 and  $M_h(\phi') = (M_h(\phi))'$

 (4.7)   $\qquad\,  X + P\A\,\,$ has $\,\,(\Delta)$,

\noindent if $\A$ has $(\Delta)$ and $\A  \st   C(\jj,X)$.

\noindent  Contrary to Example 4.21 and  Proposition 7.7 one has

 \proclaim{Remark 4.23}  For any $\jj$,
$X$,
 $\A = AAP\cdot e^{it^2} = C_0+AP\cdot e^{it^2}$ satisfies $(\Delta)$ and $\A\st \m \A$.
\endproclaim

\noindent Since our proof is too long, we omit it.

\noindent With Remark 4.23 one gets
 \proclaim{Remark 4.24} If $U$ is an additive group with
 $AP\cdot e^{it^2} \st U \st L^1_{loc} (\jj,X)$ and $U\cap C_0
 (\jj,X)=\{0\}$, then $U$ does not satisfy $(\Delta)$.
 \endproclaim

\demo{Proof} Assuming $U$  has $(\Delta)$ and  $\phi \in L^1_{loc}
(\jj,X)$ with $\Delta_h \phi \in AP\cdot e^{it^2} \st U, \st
AAP\cdot e^{it^2}$. With Remark 4.23 one gets

$\phi- M_h \phi = u \in U$, $u=v+w$ with $w\in AP\cdot e^{it^2}$
and $v\in C_0$, implying $v= u-w \in U$, then $v=0$, $u =w \in
AP\cdot e^{it^2}$ contradicting Example 4.21. \P

\enddemo

\proclaim {Examples 4.25} $U= UAA\cdot e^{it^2}$ or $LAP_{ub}\cdot
e^{it^2}$ or $ r(\phi,\r, X)\, |\,\jj\cdot e^{it^2}$ do not
satisfy $(\Delta)$, with $UAA$ respectively $LAP_{ub}=$ almost
automorphic respectively Levitan almost periodic functions, $
r(\phi,\r, X)=$ recurrent functions (see [9, p. 10]).
\endproclaim

 \noindent    Further $(\Delta)$ results can be found in \S7, Propoposition 8.9 and
         Remark 8.10.

\smallskip

\head{\bf \S 5  Mean classes, derivatives and uniform continuity
}\endhead

\noindent In this section we complement (2.17), (2.19), Theorem
2.11(b), and also the discussion of strict inclusions.
 Furthermore a far-reaching
generalization of results of Bochner and Upton will be given in
Proposition 5.6 see [22, p. 442], [74, p. 15, Theorem 11].

\proclaim {Proposition  5.1}

\noindent (i) If  $\A \st L^1_{loc} (\jj,X)$ or $ \st \h' (\jj,X)$
satisfies $(\Delta_1)$ and is invariant with respect to
multiplication by   positive numbers, then $\m\A \st\A_{Loc}
+\A'_{Loc}$;

  if such an $\A$ satisfies $(\Delta'_1)$ instead of
 $(\Delta_1)$, $\widetilde {\m}\A \st\A +\A'$.

\smallskip

\noindent (ii)  If  $\A  \st L^1_{loc} (\jj,X)$ or $ \st \h'
(\jj,X)$ is real-linear,
 positive-invariant and satisfies $(\Delta)$ respectively $(\Delta')$, then for $n\in \N$

\smallskip

(5.1 ) $ \,\,\,\qquad \m^n\A = (\m^{n-1}\A)_{Loc} +
(\m^{n-1}\A)'_{Loc}\,\,\,\,$ respectively

  $\,\,\,\qquad\,\,\,\qquad\widetilde {\m^n}\A =   \widetilde {\m}^{n-1} \A+ (\widetilde {\m}^{n-1} \A)'$.

\smallskip

\endproclaim
\noindent Here

 $\A_{Loc} : = \A\cap L^1_{loc}(\jj,X)$,
$\qquad \A^{(n)} : =\{T^{(n)}: T\in \A  \}$,

 $\A'= \A^{(1)} $,       $\qquad \A'_{Loc} :=(\A')_{Loc}=\A^{(1)}_{Loc}$,

\smallskip

(5.2 ) $ \,\,\, \A^{(n)}_{Loc} = \{\phi\in L^1_{loc}(\jj,X):$

$\qquad\qquad\qquad$ to $  \phi$ exists $ \psi\in \A \cap
W^{1,n}_{loc} (\jj,X)$ with $\phi= \psi^ {(n)} a.e.  \}$,
 $n\in \N_0$.
\smallskip

\noindent For relations between  $(\Delta_1)$ and $(\Delta'_1)$
see after Definition 1.4;
 $(\A')_{Loc}$ is strictly contained in $ (\A_{Loc})'$ in general.

\demo {Proof} (i)
 Case $\m\A$ : Let $\phi\in \m\A, \st L^1_{loc}$. Then $\phi=(\phi -M_1\phi) +M_1\phi$
 with $M_1\phi \in  \A$.    By the  assumption,
$(P\phi)_h - P\phi = hM_h \phi \in\A$ for all       $0< h\in \r$.
Hence by $(\Delta_1)$,
 $(P\phi-M_1 P\phi)\in\A$. Now for $ \psi \in  L^1_{loc}(\jj,X)$,  $h>0$,

\smallskip

(5.3 ) $ \qquad\qquad\,     M_h P\psi -  P M_h \psi = c_h = (M_h
P\psi )(\al_0) \qquad$ on $ \jj$

\smallskip

\noindent (proof by differentiation), so $ \psi := P(\phi - M_1
\phi) - C_1 \in\A \cap  W^{1,1}_{loc}(\jj,X)$ by [53, Theorem
3.8.6, p. 88], with $\psi' = \phi -M_1  \phi$ a.e.
 This means   $\psi' = \phi - M_1 \phi \in \A'_{Loc}$.

  Case $ \widetilde {\m}\A$ :
 Let $T\in \widetilde{ \m}\A$. Then $T=(T -  \widetilde{ M}_1)T +  \widetilde{M}_1 T$.
By [71, p. 51, Th\'{e}or\`{e}me I] there exists     $S\in \h'(\jj,
X)$ such that $S'=T$. One has $S_h-S= h  \widetilde{ M}_h T \in \A
$ for all $h > 0$. Hence by $(\Delta_1)$, $(S-  \widetilde{M}_1
S)\in\A$. Since $(S -  \widetilde{M}_1 S)' =  (T-  \widetilde{M}_1
T)$ and   $  \widetilde{M}_1 T \in  \A$,
 one gets  $   \widetilde{ \m}\A \st\A +\A'$.

(ii) for $n=1$ follows by Lemma 2.2, $\m \A$, $ \widetilde {\m}\A$
are real-linear, with

 \smallskip

 (5.4)  $ \,\,\,\qquad \A$ as in Proposition 5.1 (ii) implies

$\,\,\,\qquad\,\,\,\qquad\A_{Loc} \st \m \A$
 respectively $\A  \st \widetilde {\m}\A $.

 \smallskip

 \noindent (5.1) for general $n$ follows from this for   $\m^{n-1} \A$ respectively  $ \widetilde {\m} ^{n-1}\A$ with
Lemma 2.3.
        \P
\enddemo

\proclaim {Corollary  5.2}  If   $\,\,\,\A \st L^1_{loc}(\jj,X)$
or $ \,\,\,\st \h' (\jj,X)$ is  real-linear, positive-
 invariant,
  satisfies $(\Delta)$ and   $n\in \N$, then

(5.5) $ \qquad\qquad\,   \m^n\A =  \sum_{j=0}^{n} \A^{(j)}_{Loc}$.
\endproclaim
\demo {Proof} This follows from (5.1) inductively   with
$(\A^{(j)}_{Loc})'_{Loc}=\A^{(j+1)}_{Loc} $, Lemma 2.2 and (1.11),
(5.4).\P
\enddemo
\noindent {\bf{Remark}}. (5.5) is false for $n=1$, $\A=\A_{e}$ of
Example 3.1.

For $\jj = \r$  one can show more:

\proclaim {Corollary  5.3}  Assume   $\A  \st \h' (\r,X)$, real
linear, positive-invariant and   $\A * L_c^{\infty} (\r,\r) \st
\A$, $n\in\N$; then

\noindent (i) $\,\,$  If  $\A$  satisfies $(\Delta')$,  $\,\,$
$\widetilde {\m^n}\A =  \A + \A^{(n)}$,

\noindent (ii)  $\,\,$ If  $ \A\st L^1_{loc}(\r,X)$ satisfies
$(\Delta)$,  $\,\,$  $ \m^n\A =  \A + \A^{(n)}_{Loc}$.
\endproclaim

\noindent Here $L_c^{\infty} (\r,\r):= \{ f\in L^{\infty}
(\r,\r):\, $ supp $f$ is compact $\}$.

 \demo {Proof} (i) By Proposition 5.1, Lemma 2.2 and
induction one gets  $\,\,\, \widetilde {\m^n}\A =$

\noindent  $ \sum _{j=0}^{n} \A^{(j)}$. For
 $T \in  \A$ and  $ 1\le j \le n$ Lemma 2.8 gives  $ T= U_j^{(n-j)} +\zeta _{n-j} *T$, with $U_j \in \A$ since all
$\gamma F_{n-j }\in  L_c^{\infty} (\r,\r) $;  so      $ T^{(j)} =
U_j^{(n)} +\zeta _{n-j}^{(j)} *T \in  \A^{(n)} + \A  $ for $1  \le
j < n  $.

(ii)  As in (i) one gets
 $ \m^n\A \st (\sum_{j=0 }^n  \A^{(j)} )_{Loc}  $ and if $f\in \A$, $f^{(j)} = g_j^{(n)}+ h_j $
with  $g_j, h_j  \in \A$. So, for            $\phi \in\m^n \A$
one gets
 $\phi  = F+ G^{(n)} $ with  $ F, G \in \A$. Since here  $ \phi, F\in  L ^1_{loc} (\r,X)$,
  also the distributional derivative
 $ G^{(j)}\in  L ^1_{loc} (\r,X)$  for $1  \le j < n  $ and thus   $\in \A^{(n)}_{loc}$ by (2.5).
        \P
\enddemo

 In a similar way, one can show the following extension  to   $\jj\not = \r$:

\noindent If   $\A \st L^1_{loc}(\jj,\r)$, is  real-linear,
positive-invariant with    $(\Delta)$, and if there exists  $
\widetilde {\A} \st  L ^1_{loc} (\r,X)$ with
  $\A = \widetilde { \A}|\,\jj$, and         $ \widetilde {\A} * L_c^{\infty} (\r,\r) \st  \widetilde {\A}$, then
$ \m^n\A = \A + \A^{(n)}_{Loc} $  for $ n\in \N$.

\smallskip

\noindent Examples are $\A =AP(\jj,X), AAP(\jj,X), S^p AP(\jj,X),
\E_{ub} (\jj,X)$.

\proclaim {Corollary 5.4} Any Stepanov $ S^p$-almost periodic
function $f$ can be written in the form

(5.6) $\qquad f= u+v'$ a.e., with $\,\, u,v \in AP(\jj,X)$ and
$\,\,v\in W^{1,1}_{loc} (\jj,X). $
\endproclaim
\noindent This follows from (5.5) and (3.6).

 \proclaim {Example 5.5}  In general   $ \m \A (\r, X)|\,
\jj$ is strictly contained in  $   \m(\A (\r,X)|\,\jj)$.
\endproclaim

We indicate  this for $AP(\r,\r)$, $\jj=[0,\infty )$ : Construct
first  a $g  =\sum_{n=1}^{\infty} g_n \,\in AP(\r,\r)$ with $g
|\,\jj\in C^{\infty}(\jj,\r)$, but for no $\e >0 $ is $g$ of
bounded variation in $[-\e, 0]$. This one  gets with $g_n \in
C^{\infty}(\r,\r)$, periodic with period $2^{n}+1$, $\equiv  0$ on
$[0, 2^n]$, $||g_n||_{\infty} \le 2^{-n}$,  and $g_n (t)= 2^{-n}
\sin (\omega_n t)$ on $[-4^{-n} +\e_ n  ,-4^{-n-1} - \e_n ]$, supp
$g_n |\,[-1,0] \st I_n := [-4^{-n}  ,-4^{-n-1}  ] $ with $\e_n =
|I_n |/(2n)$, $\omega_n = \pi \cdot 4^{2n+2}$. Then also $g(0)=0$
and $g'(0 -) $ does not exist.

\noindent If $f:= (g|\, \jj)'$, $f\in \m AP(\jj,\r)$ by Lemma 2.2.

\noindent  $f\not \in \m AP(\r,\r) |\, \jj$ : Assuming the
contrary, by Proposition 5.1  and Example 4.7  there exists $u,
v\, \in  AP(\r,\r)$ with  $v(0)=0$, $v\in W^{1,1} (\r,\r)$ and $f
=u+ v'$ a.e. on $\jj$. Integration  gives $g=w+v$ on $\jj$ with
$w(t)= \int_{0}^{t} u(s)\, ds$. This means that $|w|\le c_1 <
\infty $ on $\jj$. For $t > 0$ and a suitable $\tau \in \T (u,
\frac {1}{t})$ one gets
 $|w (-t)|\le 2c_1 +1  $, $w$ is  bounded on $\r$. Then
 $w \in  AP(\r,\r)$ by  the Bohl-Bohr theorem (see the introduction and before Proposition 6.2),
 implying $g=w+v$ on $\r$. But then $g$ would be absolutely continuous on $[-1, 0]$, a
contradiction.

\smallskip

\smallskip

   A result of Bochner [22, p. 442], [2,
   p. 6, VI], [76, p. 24] states that if $f  \in  AP(\r,\cc)$ has
   uniformly continuous derivative $f'$, then $f'$ is also ap; an analogue for
   $f ''$  can be found in [33, p. 525, problem 12], [76, p. 25],
  analogues for asymptotic ap functions respectively almost automorphic functions  respectively $\lambda$-classes in
  [76, p. 38, Theorem 2] respectively [46, p. 26] respectively [9, Corollary
  1.4.3].
 Another result of Stepanoff [73],  [21, p. 81], [2, p. 78, VII] states that if
  $f$ is Stepanoff-ap  in $S^p AP(\r,X)$ and uniformly continuous, then
  $f$ is already Bohr-ap.
 All these results are subsumed by
Proposition 5.6:

\proclaim {Proposition 5.6} If  $\A \st L^1_{loc}(\jj,X)$ or $ \st
\h' (\r,X)$ is  uniformly closed   and  $n, k\in\N_0$, then \ $
(\m^{n+k}  \A) \cap \m^{k} C_u (\jj,X) =(\widetilde { \m}^{n+k}
\A)\cap \m^{k} C_u(\jj,X) \st \m^{k} \A\, $ and
 $\, \h' _{ \A}  \cap \widetilde  {\m}^{k}C_u (\r,X) $
 $\st \widetilde { \m^{k}}\A$  for $\jj=\r$.
\endproclaim

\demo {Proof (see [13, p. 134 Proposition 3.5])} If $\phi \in \m
\A \cap C_u$, then $M_h \phi \in    \A  $ and  $ M_h\phi\to \phi$
as $h\to 0$
 uniformly, $\phi \in     \A $. Inductively, this gives   $ (\m^n  \A) \cap C_u (\jj,X) \st  \A$ and with $\A$
  also  $ \m^n  \A$ is  uniformly closed;
 $ (\widetilde { \m}^n  \A) \cap L^1_{loc} (\jj,X) =  \m^n  (\A \cap L^1_{loc}(\jj,X))$ (see  (2.18)). To
 $ \phi   \in \h' (\jj,X)    \cap C_u (\jj,X) $ there exists $(\varphi_n)\st  \h   (\r,\k )$ with
 $\phi * \varphi_n   \to   \phi$ uniformly.

\noindent For general $k$ one gets with this
 $ (\m^{n+k}  \A) \cap \m^{k}C_u = \m^k (\m^{n} \A \cap  C_u (\jj,X))  = \m^{k} \A $, similarly for
$\h'_{\A}$.
        \P
\enddemo

    With Proposition 5.6 an extension of the Bohl-Bohr-Kadets  result on
integration of ap functions can be obtained (see \S 6 after
Corollary 6.6, [54], [8]); also with Proposition 5.6 a  result of
Upton [74, Theorem 11, p. 15] can be extended.

\head{\bf \S 6   Almost periodic distributions
 and indefinite
integrals}\endhead

By the above, especially examples 3.4, 3.5, 4.7, all the results
of section 2 apply to these spaces. For example, since
 $ AP (\jj,X) \st S^p  AP (\jj,X) \st  \m AP (\jj,X)$, each inclusion being strict
 by (3.6), with (2.17), (2.19) and Examples 4.7  almost periodic distributions are scaled by
(6.1), (6.2):
\smallskip

\smallskip

(6.1) $\, \qquad\h'_{AP} (\r, X) \,\,=    \cup _{n=0}^{\infty}
\widetilde  {\m}^{n}  AP (\r,X)$

$  \qquad \qquad\qquad\qquad\qquad = \cup_{n=0}^{\infty}
\widetilde  {\m}^{n} S^p AP (\r,X) = \h'_{ S^p AP } (\r,X)$,

\smallskip

(6.2) $ \qquad\h'_{AP} (\r, X) \cap L^1_{loc} (\r,X) =    \cup
_{n=0}^{\infty}  {\m}^{n}  AP (\r,X)$

  $ \qquad \qquad\qquad\qquad\qquad = \cup_{n=0}^{\infty}   {\m}^{n}S^p  AP (\r,X)= \h'_{ S^p AP } (\r,X)\cap L^1_{loc} (\r,X) $,

\smallskip

\noindent for $1\le p <\infty$, with $ {\m}^{n}  AP (\r,X)\st
{\m}^{n}  S^p AP (\r,X) \st  {\m}^{n+1}  AP (\r,X)$,
 similarly for $\widetilde {\m}$, any $X$.

\smallskip

Almost periodic distributions can  also be characterized by
translation or compactness properties, that is Bohr's, Bochner's
and von Neumann's definition all give $\h'_{AP}$:

\proclaim{Theorem 6.1} For  $T \in \h'_{ L^{ \infty}}(\r,X)$  and
$ \Phi_T(l) :=T_l $ for $ l \in \r$, the following statements are
equivalent:

\noindent (a)  $T \in  \h'_{ AP} (\r,X)$;

\noindent (b)  $T \in  \h'_{S^p AP} (\r,X)$ for some  $p\in [1,\infty)$ [or equivalently for all $p\in [1,\infty)$ ];

\noindent (c)  For any $V =$ neighborhood of $0$  in   $ \h'_{ L^{
\infty}}(\r,X)$, the set   $\Cal {T }(T,V)$  of $V$-periods

 of $T$ is relatively dense in $ \r$;

\noindent  (d)  To each sequence $(a_m)_{m\in \N} \st \r$ there exists  a subsequence $(a_{m_n} )_{n\in \N}$ and

    $S \in \h'_{ L^{ \infty}}(\r,X)$   with $ T_{a_{m_n}} \to S$ in  $ \h'_{ L^{ \infty}}(\r,X)$  [or  in  $(\h_{ L^ 1})'(\r,X)$ ];

\noindent  (e) $\, =$  (d), with "subsequence $(a_{m_n}))_{n\in \N}$" replaced by  "subnet $(a_{n(i)})_{i\in I}$";

\noindent  (f) $\, =$  (e), with "sequence" replaced by  "net";

\noindent  (g)  $\{ T_l: l\in \r\}$ is totally bounded [= relatively compact] in   $ \h'_{ L^{ \infty}}(\r,X)$;

\noindent  (h )   $\Phi _T \in AP(\r, \h'_{ L^{ \infty}}(\r,X))$  [or  $\in  AP(\r, (\h_{ L^1})'(\r,X))$];

\noindent  (i) there exists  $f, g \in  AP(\r,X)$ and $m\in \N_0$  such that  $T = f+g ^{(m)} $  on $\h (\r,\k)$;

\noindent  (j) there exists   a sequence $(f _n)\st   AP(\r,X)$
[or equivalently a net $(T_i)_{i\in  I}\st$

$ \h'_{ AP}(\r,X)$ ]
  with $f_n$ [respectively $ \, T_i$] $ \to T$ in  $ \h'_{ AP}(\r,X)$.
\endproclaim

\noindent Here   $\Cal {T }(T,V):= \{ \tau \in \r:  T_{\tau}- T \in  V\}$;  the topology of  $ \h'_{ AP}(\r,X)$ and $\h'_{ L^{ \infty}}(\r,X)$ is given by
the seminorms  $|| T|| _{U} := $ sup $\{ || T(\va)|| : \va \in U\}$, with  $ U\st \h(\r,\k)$,
  $U$ bounded in  $ \h_{ L^1}(\r,\k)$ (see \S 1).

\demo {Proof}  $ (a) \Leftrightarrow (b)$  follows by (6.1).

$(a) \Leftrightarrow (i) \Leftrightarrow (j)$: Theorem 2.11,
Proposition 2.9.

$T_i \to S$ in $(\h_{L^1})'(\r,X) $ is equivalent with $\Phi_{T_i}
\to \Phi_{S} $ in $(\h_{L^1})'(\r,X) $ uniformly on $\r$ and  also
 equivalent with $T_i \to S$ in $\h'_{L^{\infty}}(\r,X) $, that is
uniformly on $\h_{L^1}$-bounded $U \st \h (\r,\k)$; see $(i)
\Rightarrow (h)$.

$(i) \Rightarrow (h)$: If $g \in AP(\r,X)$, by definition there is
a sequence of trigonometric polynomials $g_n \in \Pi (\r,X) $ with
$g_n\to g$ uniformly on $\r$. If now $U$ is bounded $  \st Y:
=\h'_{L^{\infty}}(\r,X) $, for $ \va \in U$ and $m\in\N_0$ and
$s\in \r$ one has

$||\Phi_{g_n^{(m)}}(s)(\va) -\Phi_{g^{(m)}}(s)(\va)||= ||\int
_{\r} (g_n (t)-g (t)) \va_{-s}^{(m)} (t)\, dt ||\le $

 $\qquad\qquad\qquad\qquad\qquad\qquad ||g_n-g||_{\infty} \cdot ||\va^ {(m)}||_{L^1} \le \e$

\noindent for $n \ge $ suitable $n_{\e}$, since sup $\{ ||\va
^{(m)}||_{L^1}: \va \in U\} < \infty$ for fixed $m$.

\noindent Since elementary calculations show $\Phi_{g_n^{(m)}}\in
\Pi (\r,Y)$, the above means $\Phi_{g^{(m)}}\in AP (\r,Y)$ by our
definition in \S 1, especially  $\Phi _f \in AP (\r,Y)$ if $f \in
AP (\r,X)$. So (h) holds if $T$ is as in (i).

$(h) \Rightarrow (g)$: With  the periodicity and uniform
continuity of $\gamma_{\om}$, the $a\cdot \gamma_{\om}:\r\to Y$
has totally  bounded range in $Y$ if $ a\in  Y$, $\om \in \r$.
Since furthermore with $f(\r)$ and $g(\r)$ also $(f+g)(\r)$ is
totally bounded, all trigonometric polynomials $\in \Pi (\r,Y)$
have totally bounded range. $\Phi _T \in  AP (\r,Y)$  can be
uniformly approximated by such, implying $\Phi _T (\r)= \{ T_s:
s\in\r\}$ is totally bounded in $Y$. For ``relative compact'' see
Remark 6.3 and (6.3) below.

$(g) \Leftrightarrow (\tilde {g}), =$ set of translates $\{
(\Phi_T)_t: t\in\r\}$ is totally bounded in $ V(\r,Y)$, where

$ V(\r,Y): = Y^{\r}$ equipped with the locally convex topology of
uniform convergence on $\r$:

\noindent Given $U \st \h (\r,\k)$, $U$ bounded in $\h _{L^1}
(\r,\k)$, the $U':= \{\va_s : \va\in U, s \in\r\}$ is also
bounded. So by $(g)$ the set $\{T_t: t\in \r\}$ can be covered by
$W_1, \cdots, W_n$ with $W_j =\{T\in Y: ||T(\va)-T_{t_j}(\va) ||
\le \e \,\,$ for $\, \va \in U' \}$. This implies that
$\{(\Phi_T)_t: t\in \r\}$ is covered by the

$\{ \Phi \in Y^{\r}: ||(\Phi (s)- (\Phi_T)_{t_j} (s)) (\va)||\le
\e $ for all $\va \in U, s\in\r \}$,

\noindent  giving $(g) \Rightarrow (\tilde {g})$. $({g})
\Rightarrow (\tilde {g})$ follows with $s=0$.

$(g) \Leftrightarrow (c) \Leftrightarrow (e) \Leftrightarrow (f)$.
This follows from the following abstract characterizations of
almost periodicity (for the definitions of nets and subnets see
[56, p. 28]):
 \proclaim{Proposition 6.2} If  $ G$ is a locally compact abelian group, and
$Y$ a topological abelian group, $V:= C(G,Y)$ is equipped with the
topology ${\frak {T}}_u$ of uniform convergence on $ G$, then for
$\Phi \in V$ the following statements are equivalent

(i) The set of translates $\{ \Phi_t: t\in G\}$ is totally bounded
in $(V,{\frak {T}}_u)$ (von Neumann ap [25, p. 22]).

(ii) To each net $(t_{\al})_{\al \in A}$ from $G$ there exists a
subnet $(t_{\al (\beta}))_{\beta \in B}$ with
$(\Phi)_{t_{\al(\beta)}}$ Cauchy in $(V,{\frak {T}}_u)$.

(iii) To each sequence  $(t_{m})_{m \in \N}$ from $G$ there exists
a subnet $(t_{m(\beta}))_{\beta \in B}$ with
$(\Phi)_{t_{m(\beta)}}$ Cauchy in $(V,{\frak {T}}_u)$.

(iv) For each neighborhood $U$ of $0$ in $Y $ the set ${\Cal
{T}}(\Phi,U)$ of $U$-periods is relatively dense in $G$ (Bohr ap).

(v) $\Phi$ is Maak ap (see [61, p. 26]): $\phi: G\to Y$ is called
Maak almost periodic  if for each neighborhood $U(0)$ of $0$ of
$Y$ there exists a partition $P_1,\cdots, P_n$ of $G$ such that

 $G=\cup_{k=1}^n P_k$, $\,\,\phi (a +s)- \phi (b+s)\in U$ for $ a$ and $b
\in$ same $P_k$, $s\in G$.
\endproclaim

\noindent

\noindent $(V,{\frak {T}}_u)$ stands here for $V$
 equipped with the
topology ${\frak {T}}_u$ of uniform convergence on $ G$.

\noindent In (iv), $M$-relatively dense in G means there exists a
compact $K\st G$ with

$\qquad M\cap (x+K)\not = \emptyset$ for each $x\in G $.

\proclaim{Remark 6.3} ($\al$) $G$ can be any abstract  abelian
group in the discrete topology, $K$-compact means then $K$ finite
in (iv).

($\beta$) If in Proposition 6.2 the $Y$ is topologically complete,
e.g. complete, then (i) is equivalent with

$(i') \, \,\{ \Phi_t: t\in G\}$ is relatively compact in
$(V,{\frak {T}}_u)$;

 (ii) is equivalent with

$(ii')$ To each net $(t_{\al})_{\al \in A}$ from $G$ there exists
a subnet $(t_{\al (\beta}))_{\beta \in B}$ with
$(\Phi)_{t_{\al(\beta)}}$ convergent  to some $\Psi \in V$
uniformly on  $ G$;

 $(iii)$ is equivalent with

$(iii')$ To each sequence  $(t_{m})_{m \in \N}$ from $G$ there
exists a subnet $(t_{m(\beta}))_{\beta \in B}$ with
$(\Phi)_{t_{m(\beta)}}$ convergent  to some $\Psi \in V$ uniformly
on  $ G$.

($\gamma$) All the above holds for $Y=$ uniform space.

($\delta$) All the above  can be extended to non-abelian $G$.

\endproclaim

We indicate a proof of Proposition 6.2  ("folklore" ) only for the
special case ($\beta$) of Remark 6.3, only this is used here.
First one has

\smallskip

\smallskip

(6.3) $\qquad$ For any $X$, $Y =\h'_{L^{\infty}}(\r,X)$ is a
complete locally convex space.

\demo{Proof} (See [71, p. 201]): Given a Cauchy-net $(T_i)_{i\in
I}$ from $Y$, then for fixed $\va \in \h (\r,\k)$,  $(T_i
(\va))_{i\in I }$ is a Cauchy-net in $X$; $X$ being complete,
there exists $T: \h (\r,\k)\to X$ with $T_i (\va) \to T(\va)$ for
$\va \in \h (\r,\k)$; $T$ is obviously linear. By definition of
$Y$, $T_i (\va) \to T(\va)$ uniformly on $U \st  \h (\r,\k)$ which
are bounded in $ \h_{L^1} (\r,\k)$, especially on $\{ \va_m : m\in
\N\}$ if  $\va_m \to 0$ in $\h_{L^1} (\r,\k)$, $\va_m \in \h
(\r,\k)$. $\h_{L^1} (\r,\k)$ being metrizable, this gives the
continuity of $T$, i.e. $T\in \h'_{L^{\infty}}(\r,X)$. \P
\enddemo

\noindent (6.3) means that for Proposition 6.2 we can assume $Y$
complete, thus also $C(G,Y)$ is complete.

``Totally bounded'' in $(g)$ or $(\tilde {g})$ is then equivalent
with relatively compact ([57, p. 36, (2) a)]). Furthermore, the
equivalence of $(i),\, (ii),\, (iii)$ is then via $(i'),\, (ii'),
\, (iii')$ a general topological result ([57, p. 36, (2)b), (3)]).

\noindent The equivalence of  $(i)$ and $(iv)$ follows directly
from the definitions.

\noindent $(iv) \Rightarrow (i)$ follows as in [2, p. 5, VI],
where the case $G=\r$, $Y=$ Banach space $X$ is treated; similarly
$(i) \Rightarrow (iv)$ follows from [2, p. 8, ($\e$) or p. 9,
VIII]. As in $(\tilde {g}) \Leftrightarrow (g)$ above the
equivalence of $(g)$, $(c)$, $(e)$, $(f)$ follows then from
Proposition 6.2 for $\Phi_T$.

$(g) \Rightarrow (d)$.

\noindent This also follows as above from a general topological
results:

\proclaim{Lemma 6.4} If $G$, $Y$, $(V,{\frak
 {T}}_u)$ are as in
Proposition 6.2, and if additionally $Y$ satisfies von Neumann's
countability axiom $(A_0)$, then for $\Phi \in V= C(G,Y)$,

$(i') \{ \Phi_t: t\in G\}$ is relatively compact in  $(V,{\frak
{T}}_u)\qquad$ is equivalent with
\smallskip

$(vi')$ To each sequence  $(t_{m})_{m \in \N}$ from $G$ there
exists a subsequence $(t_{m_n})_{n \in \N}$  with
$(\Phi)_{t_{m_n}}$ convergent to some $\Psi \in V$ uniformly on
$G$ (Bochner's criterion [23, p. 154]).
\endproclaim

\noindent Here $(A_0)$ for $Y$ means  there exist countably many
open  $U_n$ with $ 0\in U_n \st Y$ and  $\cap _{n=1}^{\infty} U_n
=\{0\}$ (see [63, p. 4, Definition 2b (2)]).

If  $Y$ has $(A_0)$  and is additionally topologically complete,
then by the approximation theorem in [25, p. 37, Theorem 27],
$(i)-(vi)$ $are\,\, also$ $equivalent\,\, with$

 $(vii) \qquad \Phi \in AP(G,Y)$

\noindent  with $AP(G,Y)$ defined as in $\S 1$ (using continuous
bounded characters $\gamma: G\to \cc$ instead of $\gamma_{\om}$
and   replacing "sequence  $(\pi_n)  \st \Pi(\r,X$) such that
$||\pi_n - \phi|| \to  0 $ as $ n \to \infty$"
 by  "net $(\pi_i)  \st  \Pi(\r,Y)$  with  $\pi_i \to \phi$  uniformly on $G$").

Lemma 6.4 can be applied in our case  to $\Phi_{T}$: $\{
(\Phi_{T})_t: t\in\r\}$ is relatively compact by $(g) \Rightarrow
(\tilde {g})$ and the remarks after (6.3) above.

Though $Y:= \h'_{L^{\infty}}(\r,X)$ does not satisfy the first
axiom of countability $(A_1)$ (it is also not separable), it
satisfies $(A_0)$ for arbitrary $X$: $\h (\r,\k)$ is separable in
the induced topology of $\h_{L^1}(\r,\k)$

\noindent ( ${\Cal {P}}_r := \{\rho_n\cdot p: n\in \N,\,\, p$
polynomial with rational coefficients $\in \k \}$ is countable,
with fixed $\rho_n \in \h (\r,\k)$, $0\le \rho_n \le 1$, $\rho_n
=1$ on $[-n,n]$, supp $\rho_n \st [-n-1,n+1]$; approximating
$\va^{(m)}$ by rational polynomials, $m\in \N_0$ and $\va \in \h
(\r,\k)$ fixed, it can be shown that ${\Cal {P}}_r$ is dense in
$\h (\r,\k)$ as stated);

\noindent Taking $U_{m,n} := \{T\in \h'_{L^{\infty}}(\r,X):
||T(\va_n)|| \le 1/m \}$ with $\{\va_n: n\in \N\} \,$
$\h_{L^1}$-dense in $\h(\r,\k)$, one gets $(A_0)$ for
$\h'_{L^{\infty}}(\r,X)$, and thus $(g) \Rightarrow ({d})$.

 \demo{Proof of Lemma 6.4}
 $(i') \Rightarrow (vi')$: Given $U'_n$ neighborhoods of $0$ in $Y$
with  $\,\cap _{n=1}^{\infty}\, U'_n $ $=\{0\}$ there are closed
neighborhoods $U_n$ of $0$ in $Y$ with $U_n+ U_n \st U'_n$  and
$U_{n+1} \st U_n$ for each $n\in \N$. If $(t'_{m})_{m \in \N}$  is
a sequence  from $G$, with $M =  \{ \Phi_t: t\in G\}$ totally
bounded ([57, p. 36, (2)a)]) recursively there exist subsequences
$(t_{n,m})_{m \in \N}$ of $(t'_{m})$ for $n\in \N$ with
$(t_{n+1,m})_{m \in \N}$ is a subsequence of $(t_{n,m})_{m \in
\N}$  and $\Phi_{t_{n,k}} (s)- \Phi_{t_{n,l}} (s ) \in U_n$ for
$s\in G$ and $k, l \in \N $, any $n\in \N$. If $(t_m): =$ diagonal
sequence of $(t_{n,m})_{m \in \N}$, one has:
  $\,\,(t_m)$ is a subsequence of $(t'_{m})$, and

(6.4) $\qquad \Phi_{t_m} (s)- \Phi_{t_n} (s ) \in U_m$ for $s\in
G$ if  $n\ge m$, $m, n\in \N$.

\noindent $(\Phi_{t_{m}})_{m\in \N}$ is Cauchy in $(V,{\frak
{T}}_u)$:

\noindent If not there exist a neighborhood $U_0$ of $0$ in $Y$,
strictly monotone sequences $(m_k)$, $(n_k)$ in $\N$ and $(s_k)
\st G$ with

(6.5) $\qquad \Phi_{t_{m_k}} (s_k)- \Phi_{t_{n_k}} (s_k )\not \in
U_0$ for $k\in \N$.

\noindent By $(i')$, there exists a subnet $(k(i))_i$  of $(k)$
 and $f\in V $ with

$\Phi_{t_{m_{k(i)}}} \to f$ uniformly on $G$.

\noindent Similarly there are a subnet $(l(j))_j$  of  $(k(i))_i$
 and $g\in V $ with

$\Phi_{t_{n_{l(j)}}} \to g$ uniformly on $G$.

\noindent Choosing a neighborhood  $W_0$ of $0$ in $Y$ with
$W_0=-W_0$, $W_0+W_0 \st U_0$ one gets $k_0 \in \N$ with

$\,  \Phi_{t_{m_{k_0}}} (s)- f(s)\in W_0$ and  $\Phi_{t_{n_{k_0}}}
(s)- g(s) \in W_0$ for $s\in G$.

\noindent (6.5) gives $f\not = g$.

\noindent  Now by (6.4) one has, for fixed $q\in \N$ and $m\ge q$,

 $ \Phi_{t_m} (s)- \Phi_{t_{n_{l(j)}}} (s )
\in U_q$ for $s\in G$ and $j \ge j(m)$;

\noindent  since $ U_q$ is closed in $Y$, this implies

$ \Phi_{t_{{m}}} (s)- g (s ) \in U_q$ for $s\in G$ and $m \ge q$.

\noindent Similarly one gets $ f(s)- g (s ) \in U_q$ for $s\in G$.
Since  $ q$ is arbitrary, $(A_0)$ implies $f=g$, a contradiction.

\noindent $\{ \Phi_t: t\in G\}$  relatively compact gives then a
$\Psi \in V$ with $ \Phi_{t_{m}}\to \Psi$ uniformly on $G$, i.e.
$(vi')$.

\noindent We omit the proof of $(vi')\Rightarrow (i')$, since in
our case it follows
         from (6.3) and Remark 6.3.
\P
\enddemo

(In fact, the uniformity on $\{ \Phi_t: t\in \r\}$  induced by
$(V,{\frak{T}}_u)$ satisfies the first axiom of countability
$(A_1)$).

$(d) \Rightarrow (a)$: As in the proof of $(g) \Rightarrow (\tilde
{g})$, $T_{a_{m_n}} \to S$ in $\h'_{L^{\infty}} (\r,X)$ implies
$\Phi_{T_{a_{m_n}}} \to \Phi_S$ uniformly on $\r$. With definition
 $T*\va (t)= T(\check {\va}_t)$ [75, p. 156, (2)] this gives for fixed $\va \in\h (\r,\k)$, with
$\Psi := T*\va \in C_b (\r,X)$: $\Psi (s+a_{m_n}) \to S*\va (s)$
as $n\to \infty$, uniformly in $s\in \r$. Since this holds for any
sequence $(a_m)$ from $\r$, Bochner's criterion for $X$-valued ap
functions on $\r$ [23, p. 154], [2, p. 9, VIII, p. 15, I] gives
$\Psi \in AP(\r,X)$ as defined in $\S 1$. $\va$ being arbitrary,
(a) follows.

\noindent The proof of Theorem 6.1 is complete. \P
\enddemo

 \noindent By (6.1), (6.2) the Stepanoff $S^p AP (\r,X)$
functions and distributions are subsumed   by almost periodic
distributions $\h'_{AP} (\r,X)$.

\noindent For the    functions $AAP (\r,X)$  and
  ${ EAP} (\r,X)$ this is not the case. One can  even show  $C_c  (\r,X)\cap  \h'_{AP} (\r,X)= \{0\}$,
so   $\h'_{AP} (\r,X) \st \h'_{AAP} (\r,X) \st \h'_{EAP} (\r,X)$, each inclusion being strict; however,
already $\m AP \not \st EAP $.

\noindent  Similarly, there exists $ g\in \m AP (\r,\r)$  with  $g \not \in B^p AP (\r,\r)$. Conversely, not even
 $W^p AP (\r,\r) \st\h'_{AP} (\r,\r)$  in the following strong sense:
There exists Weyl ap function  $ f\in  C^{\infty} (\r,\r)$ with $
f\in W^p AP (\r,\r)$ for all $p\in [1,\infty)$ such that  $
[f]_{B^1} \cap \h'_{AP} (\r,\k)= \emptyset$, where $[f]_{B^1} :=$
Besicovitch equivalence class $\{ F\in L^1_{loc } (\r,\k):
||F-f||_{B^1} =0 \}$. (For $f$ one can use a slight modification
of the $F$ of "main example IV" of [27, Appendix, pp. 131-133],
suitably refining the arguments and lemmas there; there is no
space for details).
\smallskip

\smallskip

The following is a generalization of the Bohl-Bohr-Amerio-Kadets
theorems [21, p. 7], [37], [54], [8 and references therein] (see
the introduction, and after Corollary 6.6).

 \proclaim
{Proposition 6.5} If  $ U, \A \st L^1_{loc}(\jj,X)$ or $ \st \h'
(\jj,X)$,  $\A$ satisfies $(L_U)$ then  $  \m  \A$ respectively
$\widetilde { \m}   \A$ respectively   $ \h'_{\A}(\r,X)$ for
$\jj=\r$ satisfy  $(L_{\m U})$ respectively $(L'_{\widetilde { \m}
U})$ respectively $(L'_{\h'_U})$ (for $(L_U)$ etc. see \S 1).
\endproclaim

\demo {Proof} The definitions and  $\Delta_h M_s = M_s \Delta_h$,
$\Delta_h  \widetilde { M_s} = \widetilde { M_s} \Delta_h$. \P
\enddemo

To apply this to indefinite integrals, we need an extension of $P$
defined
 for functions in section 1 to   vector-valued distributions defined on any open interval $I$:

\smallskip

(6.6) $\qquad$ There exists  $ \widetilde {P}:   \h' (I,X) \to \h'
(I,X) $  with $  (\widetilde {P} T)' =T$, $T\in  \h' (I,X)$,

  $ \qquad\qquad\qquad\qquad\qquad\qquad \widetilde {P}$ linear, continuous.
\smallskip

\noindent This follows as in  [71, p. 52]; for $f \in
L^1_{loc}(\jj,X)$ one has however only $ \widetilde {P} f =Pf +c
$, with the constant $c$ depending on $f$.

\proclaim {Corollary 6.6}  $ U, \A   \st \h' (\r,X)$, $\A$
satisfies $(L_U)$, is linear and contains all the constants, then
 $ T\in   \h'_{\A} (\r,X)$, $\, \widetilde { P} T  \in   \h'_U (\r,X)$ implies
$   \widetilde { P} T  \in   \h'_{\A} (\r,X) $.
\endproclaim

\demo {Proof} Differentiation and (2.8) give $\Delta_h \widetilde
{P} T = h   \widetilde {M}_h T   + $ constant,
 $\in  \h'_{\A} (\r,X)$
by Proposition 2.4 (a).
        \P
\enddemo

\noindent Special case :
  $\A= \ AP (\r,X)$ with $c_0 \not \st X$, $U= L^{\infty}( \r,X)$, then  $\A$  satisfies $(L_{ub})$  by
   [13, p. 120] and then $(L_b)$ with Proposition 1.7 (i).

\smallskip

\noindent $ \A = UAA(\r,X)$ and $LAP_{ub}(\r,X)$ are also possible
(see [13, p. 120]).

\smallskip

\noindent Similar results hold for  $  {\m}^n\A$  and $ \widetilde
{\m}^n \A$ instead of   $\h'_{\A} (\r,X)$, this will be treated
somewhere else.

\noindent Corollary 6.6 contains the classical Bohl-Bohr-Kadets
theorem:

$If$  $\phi \in AP (\r,X)$ with $c_0 \not \st X$ $and\,\, only $ $
P\phi\in  \h'_{L^{\infty}} (\r,X)$  ($for\,\, example$  $ P\phi
\in S_{b}^{1}(\r ,X)$), $then$
 $ P\phi \in AP(\r,X)$.

Indeed,
  $ P\phi \in\h' _{AP} (\r,X)$  by the special case above and the remark  after
   (6.6). But then  $ P\phi \in AP(\r,X)$ by Proposition 5.6.

\noindent  Corollary 6.6 can be extended to the half line, this
will be shown in a future note.

\head {\bf \S 7  Ergodic classes}\endhead

 This section  is devoted to the study of ergodic classes.  We give inclusion relations between the
various $\m^n\A$, $\A = \E$, $\E_{ub}$, $\T\E$, $\T\E_0$, and
apply results of section 4 to show $(\Delta)$ for $\E(\jj,X)$,
$\E_0(\jj,X)$, $\E_n(\jj,X)$, $Av_0(\jj,X)$,  $Av_n(\jj,X)$.

\proclaim {Proposition 7.1} For any $\jj$, $X$, the following
inclusions hold and are strict

\smallskip

 (7.1) $\qquad\qquad
\T \E_{ub} \st\E_{ub} \st \m\E_{ub}  \st \E   \st \m^2   \E_{ub}
\st  $

 $ \qquad\qquad\qquad\,\, \st\m  \E  \st \m^3  \E_{ub}  \cdots  \st\m^n  \E \st\m^{n+2}  \E_{ub} \st \cdots    $,

 $ \qquad\qquad\qquad\,\, \st \h' _{ \E_{ub}} (\r,X) = \h'_{  \E} (\r,X)$                    for $\jj=\r$.

\smallskip

\noindent (7.1) holds also, if there everywhere    $ \E$ is
replaced respectively  by  $ \E_0$, $\T\E_0$ or $\T \E$.

\noindent  Furthermore, $   \, C_{ub} (\r,X) \cap \h'_{ \T \E}
(\r,X) =  \T\E_{ub} (\r,X)$ and  $  C_{ub} (\r,X) \cap \h'_{ \T
\E_0} (\r,X) = \T \E_0 (\r,X)\cap  C_{ub} (\r,X)$; but $\E_{ub}
\not\st \cup_{0}^{\infty} \m^n\T\E = \cup_{0}^{\infty}
\m^n\T\E_{ub}$.
\endproclaim
\demo {Proof} Case $\E$:   If $\phi \in\m \E_{ub} (\jj,X) $,
$\phi=\psi+\xi$ with $\psi \in  \E _{ ub}(\jj,X) $ and  $\xi \in
(\E _{ ub})'_{Loc}(\jj,X) $ by Proposition 5.1 (i) and Examples
4.7.
 This means that $\xi \in  \E(\jj,X) $ and proves  $\phi \in  \E (\jj,X) $.
  If $\phi \in  \E(\jj,X) $,  $M_{h}\phi \in  \E(\jj,X)$ for each $h> 0$ by (2.6).  Now

      (7.2) $\qquad\qquad
\E (\jj,X) \st \m C_{b} (\jj,X)$,

 \noindent since $M_T \phi$ and  $M_{T+h}\phi$ are bounded for suitable  $T $, also for open $\jj\not =\r$. This
 implies that  $ M_{\tau}M_{h}\phi \in  \E_{ub}(\jj,X)$ for all $\tau,  \, h > 0 $.

If $f \in C_{ub} (\jj, \r)$ is    defined by $f = 1$ on  $ I_{2n}$
and $f = 0$  on $I_{2n+1}$, where  $I_n =[10^n  +1, 10 ^{n+1} -1]$
then  $\gamma_{\om} f\in \E_{ub} (\jj, \r)$ for all $\om \not = 0
$, but $f\not \in\E$. This gives
 $\T\E_{ub}  (\jj,X)  \not = \E_{ub}  (\jj,X)$.

If $f(t)=  \sin t^2 $,    $ f'\in \E(\jj,\r)$; Lemma 2.2 gives   $
f^{(n)}\in \m^{n-1 }\E(\jj,\r)$. If   $f^{(n)} \in \m^n\E_{ub}
(\jj,\r)$,   (2.8) gives $\Delta_{h_n }\cdots \Delta_{h_1 } f \in
\E_{ub} (\jj,\r) \st C_{ub} (\jj,\r)$, then Proposition 1.7 (i)
inductively
 $f\in   C_{ub} (\jj,\r) \st C_{u} (\jj,\r)$, a contradiction.
This means  that the inclusions   $ \m^n\E_{ub}(\jj,X)
\st\m^{n-1} \E (\jj,X)$  are strict for all $n\in \N$.

\noindent We omit the examples for $ \m^{n-1}\E (\jj,X)  \not =
\m^{n+1} \E_{ub} (\jj,X)$.

\noindent  $\m^n\E_{ub}\st \h'_{\E_{ub}}$ by Corollary 2.5 and
Example 4.7,
 $ \h'_{\E_{ub}} = \h'_{\E}$  then by Corollary 2.14, Example 3.4 and
 $\E_{ub}\st\E\st\h'_{\E_{ub}}$ of (7.1).

\noindent Case $\E_0$ can be proved similarly.

 \noindent  For  the $\T\E,\, \T\E_0 $-cases   we need first $ \T\E\st
\m\T\E$, $ \T\E_0\st \m\T\E_0$:

  If $\phi \in\T  \E(\jj,X) $ (respectively  $\T  \E_0(\jj,X)$),
 $M_h \phi \in C_b (\jj,X)$  by  (7.2). This implies $\gamma_{\om} M_h \phi \in C_b (\jj,X)$ for all  $\gamma_{\om}$.
 Therefore $(\gamma_{\om}M_h\phi)'\in \E_0(\jj,X)$.
Since  $\gamma'_{\om} \,M_h\phi=  (\gamma_{\om}M_h\phi)'-
\gamma_{\om}(\Delta_h\phi)/h$, one gets
 $\gamma'_{\om}M_h\phi \in \E(\jj,X)$ (respectively  $ \E_0(\jj,X)$). This gives  $(\gamma_{\om} M_h\phi)\in \E(\jj,X)$(respectively  $ \E_0(\jj,X)$)
  for all $ {\om}\not = 0$ and hence    $M_h \phi \in\T \E (\jj,X) $  (respectively  $\T  \E_0(\jj,X)$) with   $ \E_{ub}\st \E \st \h'_{\E_{ub}}$
   ($ \E_{0,ub}\st \E_0 \st \h'_{\E_{0,ub}}$) of (7.1).

\noindent (7.1) for  $\T \E$  (respectively  $\T \E_0(\jj,X)$)
follows then as for $ \E$, especially $\h'_{\T \E_{ub}} =\h'_{\T
\E} $ ( respectively $\h'_{\T \E_{0,ub}} =\h'_{\T \E_0} $).

\noindent $ C_{u} \cap \h'_{\T \E} =\T \E_{ub} $ follows with
Proposition 5.6 and the above.
        \P
\enddemo

\proclaim {Remark 7.2}

(i) There exists even $f\in C_{ub}(\jj,X) \cap  C^{\infty}(\jj,X)
$ with $\gamma_{\om} f \in
 \E_{ub}(\jj,X)$ for $\om \not = 0$ but $f\not \in \T \E_{ub}(\jj,X)$.

 (ii) For any $\jj$, $X$,
 $k\in \N_0$ one has
 $\m^{k} L^{\infty}\cap  \m^{k} \E =\m^{k} L^{\infty}\cap  \m^{k+1} \E_{ub} $.
\endproclaim

 \proclaim{Proposition  7.3}
  $\A = \E_0$, $\E_n$, $Av_0$, $Av_n$  have     $  (\Delta)$    for arbitrary   $  \jj, X$.
\endproclaim

\demo {Proof} Case  $\A = \E_0$:

 Given $\phi \in L^1_{loc} (\jj,X)$,  with  $  \Delta_h \phi\in \E_0(\jj,X)$ for all $h > 0$. For fixed $\e > 0$
define

$ A^{\e}_{m}:= \{ f \in L^1_{loc} (\jj,X):   ||M_T f (t) || \le
\e, $ for all $t\in \jj $, $T\ge m \}$, $E_{\e} := \cup _{m\in\N}
A^{\e}_{m}$.

\noindent One can apply Lemma 4.12 ($M_T\Delta _s \phi = \Delta_s
(M_T \phi)$ is continuous in $s$), so $ \phi -(M_h \phi) _v \in
E_{\e}$ for $0 < h \le \delta $. $\Delta_v  \phi \in \E$  implies
$\Delta_v  \phi \in E_{\e}$ for all $\e$, then $ M_h\Delta_v  \phi
\in E_{\e}$,  yielding $\phi- M_h  \phi \in E_{2\e}$ for all $ 0
<h \le \delta$.  With

 $  ( \phi-M_{2^{n+1} h}  \phi) =  ( \phi-M_{2^{n} h}  \phi)- \frac{1}{2}  M_{2^{n} h}  \Delta_{2^ n h} \phi$

 and $\e_n >0$ with  $\sum _{n=1}^{\infty} \e_n < \e  $

 \noindent one gets analogously  $ \phi -M_h \phi \in E_{3\e}$  for all $h > 0$. Thus
 $ \phi -M_h \phi \in \cap _{\e >0} E_{\e}\st \E_0$.

The proofs of the cases $\E_n$, $Av_0$, $Av_n$ are similar,
replacing the above $A^{\e}_m$  respectively by

$  \{ f \in L^1_{loc} (\jj,X):   M_T (|f) |(t) \le \e, $ for all
$t\in \jj $, $T\ge m \}$,

$  \{ f \in L^1_{loc} (\jj,X):   || M_T (f) (\al_0)|| \le \e, $
for
 $T\ge m \}$,

$  \{ f \in L^1_{loc} (\jj,X):   M_T (|f|) (\al_0) \le \e, $  for
$T\ge m \}$. \P
 \enddemo

\smallskip

\proclaim{Proposition  7.4} If $U\,$ and $V\,$  are additive
groups, positive invariant, with  $\,\,\,(\Delta)$,$\st
L^1_{loc}(\jj,X)$, with closed $\jj$ and any $X$, if furthermore
$U \cap V=\{0\}$, $V \st C(\jj,X) $ and to $U$ exists $m_0 \in
\N_0$ with  $U\st \m^{m_0} A v_0 (\jj,X) $ and finally, $U,$ $V$
are  invariant if $\jj=\r$, then $ U+V$ satisfies  also
$(\Delta)$.
\endproclaim
\demo{Proof} With Proposition 4.18 we have to show that in
$\Delta_h \Phi (t)= u(t,h) +v(t,h) $ the $v(\al,h) $ in $h\in
(0,\infty)$ is $h$-measurable on $ (0,\infty)$, $\Phi$ as there.
Now, by  Proposition 4.18 one has, also if $\jj=\r$,

(7.3)$\qquad \Delta_h \phi (t)= v(t, h)$, $t\in\jj$, $h \ge 0$
respectively $h\in\r$  if $\jj=\r$,

\noindent with $\phi (t):= v(\al, t-\al)$, $t\in\jj$ and
$v(\cdot\, , h ) \in V \st C(\jj,X)$ for $h\in \r^+ $ respectively
$\r$.

We now first assume $X=\k$, $=\r$ or $\cc$. By Corollary 7.5
below, a refinement of a result of de Bruijn, there exists a $
g\in C(\jj,\k)$ and an additive $H:\r\to \k$ such that

(7.4)$\qquad \phi =g +H$ on $\jj$.

\noindent with $\Delta_h \Phi= u+ v$ and (7.3) one gets, with $F:=
\Phi -g $, $\in L^1_{loc}(\jj,\k)$

(7.5)$\qquad\Delta_ h F (t) =u(t,h) +H (h)$, $\qquad t\in \jj$.

\noindent Applying the means $M_n$, with respect to the variable
$t$, to both sides of (7.5), one gets

$  G_n (t,h):= M_n (\Delta_ h F) (t) = M_n (u) +H (h)$, $h\in
\r_+$,
 where $G_n \in C(\jj\times \r_+,\k)$.

\noindent Choosing $a\in (\r_+)^{m_0}$ and applying $M_a $ to $G_n
$, one gets

$ M_a G_n (\cdot\,,h)= M_a  M_n (u) +H (h)= M_n M_a u +H (h)$,
$h\in \r_+$.

Since by assumptions $ M_a u\in A v_0 (\jj,X) $, $M_n  M_a u (\al)
\to 0 $ as $n\to \infty$ by the definition of $A v_0 $.

This gives  $q_n (h):= M_a G_n (\cdot\,,h) (\al) \to H(h)$ for
$h\in \r_+$ as $n\to \infty$. Since the $G_n$  are
$(t,h)$-continuous on $\jj\times \r_+$, the same holds for $M_a
G_n$, so the $q_n$ are continuous on $ \r_+$, with $q_n (h) \to
H(h)$ pontwise on $ \r_+$.

\noindent This implies the measurability of $H$ on $\r_+$; the
additivity of $H$ on $\r$ gives the  measurability of $H$ on $\r$.
$H(s+h)-H(s)=H(h)$ for $s,h\in \r$ and Lemma 4.19 imply $H \in
L^1_{loc}(\r,\k) $; integrating with respect to $h$ over $[0,1]$
gives the continuity of $H$ on $\r$ (even H(t)=ct). (7.4) shows
that $\phi$ is continuous on $\jj$.

 For general $X$, choose $y\in$ dual  $X^*$ and consider $\psi :=
y\circ \phi$ with $\phi$ of (7.3). One has

$\Delta_h \psi \in C(\jj,\k)$ for $h\in \r_+$, $\psi =g+H$ by the
above, $\Delta_h (y\circ \Phi-g)= y\circ u +H (h)$, with $g$, $H$
depending on $y$.  Since $u\in \m^{m_0} Av_0 (\jj,X)$ implies
$y\circ \Phi \in  L^1_{loc}(\jj,\k) $ ($y (M_h f)= M_h (y\circ f)$
for $ f\in  L^1_{loc}(\jj,X) $ by Hille's theorem  [75, p. 134,
Corollary 2], by the above one gets the continuity of $\psi=y\circ
\phi$ on $\jj$.

\noindent $y$ being arbitrary $\in X^*$, this means that $\phi$ is
weakly continuous on $\jj$; but then $\phi (\jj)$ is weakly
separable and  then (norm)-separable. Since $\phi$ is weakly
measurably, so by Pettis' theorem ([75, p. 131] or [50, p. 158,
Satz 5]) $\phi$ is (Bochner-Lebesgue) measurable on $\jj$. (7.3)
gives the desired measurability of $v(\al, \cdot )$ on $\r_+$.
\enddemo

To complete the proof of Proposition 7.4, we use the following
result of de Bruijn [30, p. 197, Theorem 1.3]

\proclaim{Theorem }  If $I\st \r$ is an arbitrary interval, $f:
I\to \r$ is an arbitrary function such that for each $h\in \r$ the
difference

$(\Delta_h f) \,|\, (I\cap (I-h))$ is continuous  on $I\cap
(I-h)$, then there exists $g$, $H$ with

\smallskip

(7.6)$\qquad \, f= g+H$ on $I$,

\smallskip

\noindent  $g\in C(I,\r)$ and $H:\r\to \r$ additive, i.e.
$H(t+s)=H(t)+H(s)$ for all $s,\, t\in \r$.
 \endproclaim

We need here the following reformulation, $\k=\r$ or $\cc$.

\proclaim{Corollary 7.5}  If $I\st \r$ is an arbitrary interval,
$f: I\to \k$ is arbitrary  such that for each $h >0$ the
difference $(\Delta_h f) \,|\, (I\cap (I-h))$ is continuous  on
$I\cap (I-h)$, then there exists $g$, $H$  with $g\in C(I,\k)$,
and $H:\r\to \k$ additive with (7.6).
\endproclaim
For $I= [0,\infty)$ or $\r$ the assumptions mean: $\Delta_h f$ is
continuous on $\jj$ for all $h >0$.

 \demo{Proof} The extension of (7.6) to complex valued $f$,
 $g$, $H$ is obvious. By assumptions $(\Delta_h f) $ is defined and continuous  on
$I\cap (I-h)$ if $h >0$: The substitution $t=s-h$ gives
$-(\Delta_{-h} f) (s) $ is defined and continuous  on $I\cap
(I-h)+h=  (I+h)\cap I = I\cap (I - (-h))$, so $(\Delta_h f) $ is
continuous on $I\cap (I-h)$ also for negative $h$ as needed for
(7.6). \P
\enddemo
\proclaim{Remark 7.6} Under additional assumptions extensions of
Proposition 4.18 and Proposition 7.4 to open $\jj =(\al,\infty)$
are possible: See  Corollary 7.8 and Remark 7.12 below.
\endproclaim
\proclaim{Corollary 7.7} Let $U\,$  be an additive group, positive
invariant,  $\st \m^{m_0} Av_0 (\jj,X)$
 with $m_0\in \N_0$, $\jj$ closed, and
$U$  invariant if $\jj=\r$. Then  $ U+X$ has $(\Delta)$ if and
only if $U$ has $(\Delta)$.
\endproclaim

\demo{Proof} ``if'': One has $U \cap X =\{0\}$, since $a\in X$,
$a\in U \st  \m^{m_0} Av_0 (\jj,X)$ implies, for $h_0 \in
(\r_+)^{m_0}$,

$a= M_{h_0} a \in Av_0 (\jj,X)$, so $a =  \lim _{T\to \infty}
\frac{1}{T} \int_{\al}^{\al+T} a\, dt = 0$, we can apply
Proposition 7.4.

\noindent We omit the proof of ``only if''. \P
\enddemo
\proclaim{Corollary 7.8} $\E$ and $Av:= Av_0 +X$ satisfy
$(\Delta)$ for any $\jj$, $X$, so also $\E \cap L^{\infty}$, $\E
\cap C_b$, $\E_u :=\E \cap C_u$. For $\T\E$ and $\T\E_0$ see
Proposition 8.9.

\endproclaim

\demo{Proof} $\E= \E_0 +X$ and  Corollary  7.7 give $(\Delta)$ for
$\E$ for closed $\jj$ ; similarly for $Av$.

If $\jj= (\al,\infty)$, $\phi \in L^1_{loc} (\jj,X)$ with all
$\Delta _s \phi \in \E (\jj,X) $, $\psi := \phi-M_h \phi $ with
fixed $h>0$, $\jj_1 := [\al+1,\infty)$, then $\psi|\, \jj_1 \in \E
(\jj_1,X)$. Since $ \E (\jj,X) \st \m C_b (\jj,X)$ by (7.2),
Example 4.15 and Proposition 4.17 (a) yield $\psi \in \m C_b
(\jj,X)$ for any $\jj$, so $M_1 \psi$ is bounded on $\jj$. This
and $\psi|\, \jj_1 \in \E (\jj_1,X)$ give $\psi \in \E (\jj,X)$

 $(\Delta)$ for $Av$
and $\jj= (\al,\infty)$ follows with Remark 7.9.

 $\E \cap L^{\infty}$, $\E
\cap C_b$, $\E \cap C_u$ satisfy $(\Delta)$ with Propositions
4.17, 4.9,  Examples 4.15 and 4.7. \P
\enddemo

\proclaim{Remark 7.9}  $Av(\jj,X)= \{f\in  L^1_{loc}(\jj,X): \,\,
\lim_{T \to \infty}\frac {1}{T} \int_{r}^{r+T}\,f(t)\, dt$ exists
$= a \in X$ for $r=\al_0\}$; if the limit exists for one $r \in
\jj\not = \r$, it exists for each $r\in \jj$, and is  independent
of $r$; $m(f):=a$ extends the $m|\, \E(\jj,X)$.

\noindent Similarly for $\jj=\r$, with $\lim_{T \to \infty}\frac
{1}{2T} \int_{r-T}^{r+T}\,f(t)\, dt$ exists for one $r$.
\endproclaim

\noindent Other examples  for Corollary 7.7 are :

\noindent $O(w)+X$ satisfies $(\Delta)$ with  $O(w)$ of
Proposition 4.10, if $w \in$ some $\m^{m_0} Av_n $, e.g $w(t)\to
0$ as $|t| \to \infty$, $t\in \jj$.

\noindent $L^p+X$ satisfies $(\Delta)$ if $1\le p \le \infty$:
$p=\infty $ is Proposition 4.9 ; else $L^p  \st  \m C_0$ for
closed $\jj$, and $C_0 (\jj,X) \st Av_0 (\jj,X) $.

\noindent A simple application of Proposition 7.4 would be
$AAP=C_0 +AP$, but this is covered already  by Proposition 4.2
(i), for arbitrary $\jj$.

The main applications of Proposition 7.4  are:

\proclaim{Example  7.10} $PAP (\jj,X)$ satisfies $(\Delta)$ for
closed $\jj$, any $X$.
\endproclaim
\demo {Proof} In Proposition 7.4 take $V= AP(\jj,X)$ and $U= C_b
(\jj,X)\cap Av_n (\jj,X)$. Then $V$ satisfies $(\Delta)$ by
Examples 4.7 and $U$ satisfies $(\Delta)$ by Proposition 4.17,
Examples 4.15, Proposition 7.3.
\enddemo

\proclaim{Example  7.11} Similarly $ U+ AP(\jj,X)$ satisfies
$(\Delta)$  if $U= Av_n $, $\T Av_0$, $Av_0 \cap C_b$, $Av_0$,
$C_0$, $EAP_0$, $\E_{n,ub}:= \{ f\in \E_{ub}(\jj,X): |f|\in
\E_{0}(\jj,\r)\}$, $\E_n$, with $\jj$ closed,  any $X$.
\endproclaim

\noindent One has

$AAP= C_0 +AP \st EAP_0 +AP= EAP \st \T\E_0 +AP\st \T (Av_0) +AP$,

 $AAP= C_0 +AP \st  \E_{n,ub} +AP\st \E_n \cap C_b +AP$

$ \qquad \qquad \qquad \qquad\st PAP \st Av_n +AP \st \T
(Av_0)+AP$.

 \proclaim{Remark 7.12} Examples 7.10 and 7.11 hold also for $\jj= (\al,\infty)$, any $X$.
\endproclaim

\demo{Proof}          With $\jj_r : = [\al +r,\infty)$ one has
$(\phi - M_h \phi)|\, \jj_r  =
              u_r + f_r$,  $f_l = f|\,\jj_l$ with unique $f  \in  AP(\r,X)$, then
              $f_r = f|\,\jj_r$;  $u : = \phi - M_h \phi  - f|\,(\al ,\infty)$  gives
$u|\,\jj_r =
              u_r  \in  U|\,\jj_r$  for $r>0$, then $u  \in U $ with Remark 7.9,
Proposition           4.17 and Examples 4.7, 4.15. \P
\enddemo

\proclaim{Proposition 7.13} If $U$ and $W$ are uniformly closed
additive groups
   $\st  X^{\jj}$  with  $U  \st  Av_n(\jj,X)$ and $W  \st  r_+(\jj,X) \cap C(\jj,X)$, then
   $U + W$  is uniformly closed and the sum is direct.
\endproclaim

\noindent Here $Av_n$ is defined in section 1, recurrent functions

$ r_+(\jj,X) : = \{\phi  \in  X^{\jj} : T_+(\phi,1/n,n)$
relatively dense in $\r_+$ for each $n \in  \N \}$, and

$T_+(\phi,\e,n) : = \{ \tau \ge 0 : ||\phi(t+\tau)-\phi(t)||\le
\e\,\,$ for all  $\,\,|t|\le n, \,\,t  \in \jj \}$.

\demo{Proof}  Follows immediately from the (any $\jj$, $X$) Porada
inequality:

\noindent If  $u \in  Av_n(\jj,X)$ and $w \in  r_+(\jj,X) \cap
C(\jj,X)$,
          then

(7.7) $\qquad         ||w||_{\infty}  \le  ||u + w||_{\infty}$.

\noindent (7.7) has been introduced by Porada [64, p. 247, (i)]
for $U = EAP_0$ and $W = AP$, $\jj=\r$, $X=\cc$, then treated in
[77, Lemma 1.3] for $U = PAP_0$, $W = AP$.

\noindent (7.7) follows in turn, for $u$, $w$ as there, from

(7.8) $\qquad   w(\jj)  \st$  closure of  $(u+w)(\jj)$ in $X$.

\noindent (7.8) has essentially been shown in  [12?, Proposition
2.4, (2.2)], $w \in C_{ub}$ or $u  \in  C_b$  are not needed in
the proof. \P
\enddemo
 \proclaim{Remark 7.14} One can also show that $\E\cdot
e^{it^2}$ and $\T (\E\cdot e^{it^2})$ satisfy $(\Delta)$, see also
Remark 8.10.
\endproclaim

 \head{\bf \S 8 Fourier analysis }\endhead
In this section $X$ will be a complex Banach space, $\k = \cc$ (if
$\k = \r$, everything works with $\sin \om t$, $\cos \om t$
instead of $e^{i \om t} $).

To get Fourier coefficients and a formal Fourier series for
elements of a class $\A$, two properties are sufficient: $\A$ is
closed with respect to multiplication by characters, and there is
a (hopefully invariant and continuous) mean on $\A$:

\proclaim {Proposition   8.1}    (i)  If  $\A \st \h'(\r,X)$
satisfies $(\Gamma)$, so does  $\h'_{\A}(\r,X)$.

(ii) If   $\A \st  L^1_{loc}(\jj,X)$ (respectively $ \h'(\jj,X)$) is linear,  positive-invariant,
 satisfies $(\Gamma)$ and $(\Delta_1)$ (respectively    $(\Delta'_1)$),
then $\m^n\A$ (respectively $\widetilde {\m}^n\A$ ) satisfies $(\Gamma)$, $n\in \N$.
\endproclaim

\demo {Proof} (i) If  $T \in \h'_{\A}(\r,X)$ and $\gamma_{\om}(t)=
e^{i\om t}$,  $\gamma_{\om} T \in  \h'(\r,X)$ and  $(\gamma_{\om}
T)*\va (x)=  (\gamma_{\om} T) (\check {\va}_{-x} ) = T (
\gamma_{\om}  \check {\va}_{-x})=T ( e^{i\om x}  \check
{\psi}_{-x})$, where $\psi :=  \gamma_{-\om} \va \in \h (\r,\cc)$,
so $( \gamma_{\om} T) *\va = \gamma_{\om} \cdot (T*\psi)$, $\in
\A$  for $\va\in \h(\r,\cc)$.

(ii) follows by induction : With Proposition 5.1 (i) one has
$\m^n\A \st U+U'_{Loc}$, $U :=\m^{n-1}\A$; if $f= u+v'$ with $u\in
U$, $v\in U\cap W^{1,1}_{loc}$, then $\gamma_{\om} f =\gamma_{\om}
u - \gamma'_{\om} v +(\gamma_{\om} v)' \in  U+U'_{Loc}$ if $U$
satisfies $(\Gamma)$, since also $\gamma_{\om}v \in
W^{1,1}_{loc}$. $W^{1,1}_{loc}$ and $U$ satisfies $(\Delta)$ by
Lemma (2.3). Similarly for $\widetilde {\m}^n\A$.\P

\enddemo

\noindent Proposition 3.8 of [13] is the  special case $n=1$,
$\A\st C_{ub } (\r,X)$.

\proclaim {Example  8.2}      $\A =AP$, $S^pAP$,  $W^pAP$,
$B^pAP$, $AAP$, $UAA$,  $LAP_{ub}$,  $EAP_{rc}$, $EAP$,
$\T\E_{ub}$  and  $\T\E_b$  all satisfy even  $ AP(\jj,\cc)\cdot\A
\st \A$.
\endproclaim
\noindent This is well known respectively follows from the definitions (also for $EAP$).

\proclaim {Lemma  8.3}
 $\h'_{\T\E_{ub}}= \h'_{\T\E} \st \h'_{\E_{ub}}\, = \h'_{\E} \st \h'_{C_{ub}} = \h'_{L^{\infty}}\st \f'|\h(\r,\k)  $ for $\jj=\r$, any $X$.
\endproclaim

\demo {Proof}   Propositions 2.9, 7.1 and Theorem 2.16. \P

\enddemo
\proclaim {Proposition  8.4}  The ergodic mean
 $ m:\E (\r,X)\to X$ can be extended uniquely to an $\tilde {m}: \h'_{\E} (\r,X)\, \to X$ which is
linear and $(\h_{L^1})'$-continuous. This
$\tilde {m}$ is  translation-invariant and satisfies

\smallskip

 (8.1 ) $ \qquad  \qquad  \,   m(T*\va)   = \tilde {m} (T ) \cdot \int \va (x) \, dx$,  $ \qquad T\in \h'_{\E}(\r,X)$, $\va\in \h(\r,\k)$.

\smallskip

\noindent Furthermore $ \,\,\widetilde {M}_T (S)\to \tilde {m}
(S)$ as  $T\to \infty$ in the  $(\h_{L^1})'$-Topology $|\,
\h(\r,\k)$.
\endproclaim

\noindent $Proof$ (See Schwartz [71, p. 207]): The definition of
$\E (\r,X)$ and $m$ gives, with (1.11) and (2.3), for $f\in
\E(\r,X)$  and  $\va\in \h(\r,\k)$,

 $    | M_T(f*\va)- m(f)\cdot  \int \va (x) \, dx | \le \,  \e \cdot ||\va||_{L^1}$ on $ \r$  for $T\ge T_{\e}$.

\noindent $T\to \infty$ yields then for $f\in \E(\r,X)$,  $\va\in \h(\r,\k)$

     \smallskip

 (8.2 )$ \qquad  \qquad  \,\va*f\in \E(\r,X)$, $ \qquad m(\va*f)= m(f)      \cdot  \int \va(x) \, dx$.

\smallskip

The linear map $ {m}: \E (\r,X)\to   X$  is
$(\h_{L^1})'$-continuous, i.e.  if $(f_n) \st \E(\r,X)$  with $\int f_n (x)\cdot\va (x)\,dx\to 0$ uniformly in $\va\in U$
for each $U \st \h (\r,\k)$ which is bounded in $\h_{L^1}$, then $m(f_n)\to 0$:

For this,  choose $\va\in  \h(\r,\k)$  with   $\int \va (x) \, dx=1$. Then $U=  \{ ((s_T * \va  \check {)\,} )_ {-t} : t \in  \r , T>0 \} \st \h(\r,\k)$
is bounded in $ \h _{L^1}(\r,\k)$, since
$||( ((s_T * \va  \check {)\,} )_ {-t} )^{(j)}||_{L^1}  =|| s_T * (\va )^{(j)}||_{L^1} \le ||s_T|| _{L^1}\cdot ||\va^{(j)} ||_{\infty} = ||\va^{(j)} ||_{\infty}$
for all $ t\in \r$, $T> 0$. So

$  M_T(f_n*\va)(t)= f_n* (s_T*\va)(t)= \int f_n  (s) ((s_T * \va
\check {)\,} )_ {-t} (s)\, ds \to 0$ as $n\to \infty$, uniformly
in $t \in \r$, $T>0$.  $T\to \infty$  gives   $m(f_n*\va) \to 0$,
with (8.2) one gets $m(f_n) \to 0$.

With (8.2), Theorem 2.11 (e) can be applied to $\A=\E (\r,X)$, so
$\E (\r,X)$ is sequentially $(\h_{L^1})'$-dense in
$\h'_{\E}(\r,X)$. The standard extension process gives a linear
$\tilde {m}: \h'_{\E}(\r,X) \to X$.

For this  $\tilde {m}$ one has (8.1): If $(g_n) \st \E(\r,X)$, $
T\in \h'_{\E}(\r,X)$, $g_n\to T$ in $(\h_{L^1})'$, $\va\in
\h(\r,\k)$, as above $g_n *\va \to T*\va$ uniformly on $\r$,
$T*\va\in \E (\r,X)$. Thus
 $m(g_n *\va) \to m( T*\va)$. (8.2) for $(g_n)$ gives (8.1).

This $\tilde {m}$ is   $(\h_{L^1})'$-continuous in $0$: If $T_n
\to 0$ in $(\h_{L^1})'$, $T_n\in \h'_{\E} (\r,X)$, then $T_n*\va
\to 0$ uniformly on $\r$, so    $m(T_n*\va) \to 0$; (8.1) gives
$\tilde {m} (T_n) \to 0$ if $\int \va (x) \, dx \not = 0$.

$\tilde {m}$ is invariant  since $m |\,\E (\r,X)$ is invariant by definition.

 $\tilde {M}_T (S)\to \tilde {m} (S) =$  constant function follows using
$g_n\in  \E (\r,X)$  with  $g_n\to S$, since with $U$   also

 $\{ M_{-T} \va : \va\in U, T\ge 1 \}$ is $\h_{L^1}$-bounded.\P

\smallskip

\smallskip

\proclaim {Remark} (8.1)  holds also for $\va \in \h_{L^1}$
($\h'_{\E} \st (\h_{L^1})'| \h (\r,\k)$ by Proposition 2.9).
\endproclaim

\proclaim {Corollary  8.5}
 If $ T\in  \h'_{\E} (\r,X)$ and $ n\in  \N $,
 then $\tilde {m} (T^{(n)}) =0$.
\endproclaim

\demo {Proof} (8.1) for  $T^{(n)}, \in \h'_{\E} (\r,X)$ as   $T$,
with $T^{(n)} *\va =T*(\va^{(n)})$.\P
\enddemo

\noindent With the mean $ \tilde {m}$   of Proposition 8.4 one can
now define Fourier coefficients etc.:

\proclaim {Definition  8.6} For  $ S\in \h'_{\T\E}(\r,X)$ we
define  ($\gamma_{\om} (t)= e^{i\om t}$)

     \smallskip

 (8.3)$ \qquad  \qquad  \, c_{\om} (S) := \tilde {m } (\gamma_{- \om } S) $, $ \,\,\, \om \in\r$      $\qquad$ (Fourier coefficients of $S$).

     \smallskip

 (8.4)$ \qquad  \qquad  \,\sigma_{B} (S):= \{ \om\in \r :  c_{\om} (S) \not = 0 \}$        $\qquad$ (Bohr-spectrum of $S$).

     \smallskip

\noindent and the formal Fourier series

     \smallskip

   $ \qquad  \qquad   \qquad  \,\sum _{\om\in \sigma_B (S)} \, c_{\om} (S)\, \gamma_{\om}$.
\endproclaim

\smallskip

\noindent By Proposition 8.1,
 $\gamma_{-\om} S \in  \h'_{\E}(\r,X)$, so everything is well defined.

\smallskip

\noindent For $f\in W^p AP (\r,X) \st \E (\r,X)$, one gets the
usual Fourier coefficients, series and spectrum, with $AP\st S^p
AP \st \ W^p AP$, $1 \le p <\infty$.

\noindent $\T\E_0$ contains the $f\in \T\E $ with Fourier series
$0$.

\smallskip

\noindent We say that $U\st  \h'_{\T\E}(\r,X) $ has $countable\, spectra$, if for each $S\in U$
the Bohr spectrum $\sigma_{B} (S)$ is at most countable.

\proclaim {Proposition  8.7} If   $ S\in \h'_{\T\E}(\r,X)$,
$\va\in \h _{L^1}(\r,\cc)$ and $n\in\N$, then $S*\va \in \T\E
(\r,X)$,
 $ S^{(n)}\in \h'_{\T\E}(\r,X)$ and

     \smallskip

(8.5) $\qquad\qquad\, c_{\om} (S*\va) =  c_{\om} (S)\cdot \hat
\va(\om), \qquad  c_{\om} (S^{(n)}) =(i\om)^n  c_{\om} (S)$,
$\qquad \om \in\r$.

     \smallskip

\noindent Therefore   $ \h'_{\A}(\r,X)$ has countable spectra,  if $\A$ has,  $\A\st \h'_{\T\E} (\r,X)$.

\endproclaim
\demo {Proof}
               $\,\,\, \h'_{\T\E}(\r,X) \st  \h'_{L^{\infty}}(\r,X)=  (\h_{L^1})'(\r,X)|\,   \h(\r,\k)\,\,\,$ by
 Lemma 8.3 and

 \noindent Proposition 2.9, so $(S*\va)(x):=$ extension $\widetilde {S}(\check {\va}_{-x})\in C^{\infty}(\r,X)$.
Since $\h(\r,\k)$  is dense in $\h_{L^1}(\r,\k)$ there are
$(\va_m) \st \h(\r,\k)$  with  $\va_m \to \va\in\h_{L^1}(\r,\k)$
in $ (\h_{L^1})'(\r,\k)$; by definition  $S*\va_m \in \T\E$.
Continuity of $\widetilde {S}$ gives $ S*\va_m \to  {S} *\va$
uniformly on $\r$, so $S*\va\in\T\E$ since $\E (\r,X)$ is
uniformly closed. This also gives (8.1) for $\va\in\h_{L^1}$, then
(8.5) with $T=\gamma_{-\om} S$, Proposition 8.1 and
 $\gamma_{-\om}\cdot (S*\va)= (\gamma_{-\om}\cdot S)*(\gamma_{-\om} \va)  $, $\gamma_{-\om}\va\in \h_{L^1} $.

\noindent  Since $\hat {\va}$ is entire if $\va\in\h (\r,\k)$,
(8.5)  gives the countability of $\sigma_{B}(S)$ for $S\in
\h'_{\A}$.\P
\enddemo

\proclaim {Examples  8.8} All of the above can be applied to $\A
=AP$,   $AAP$,   $S^pAP$,  $W^pAP$,  $EAP_{rc}$, $EAP$ and the
corresponding $\h'_{\A}$, with $\h'_{AP} =
\h'_{S^pAP}\st\h'_{W^pAP} \st\h'_{\T\E} $ by [(8.3)], since all
these are  $ \st \T\E$  and have countable spectra. (See (1.2f),
[69,
 Theorem 2.4] for $EAP$).
\endproclaim

\noindent Definition 8.6 is thus  also meaningful for all mean
classes $ \m^n\A$, $\widetilde { \m}^n\A$, $ \,\A$ as in Examples
8.8, since $ \,\A \st  \T\E $ implies $ \m^n\A \st \widetilde {
\m}^n\A \st \widetilde {\m}^n\T\E \st\widetilde { \m}^n\h'_{\T\E}
=\h'_{\T\E}$. The last equality follows from Propositions 2.4 (b)
and 8.9; one even has  (2.17), (2.19) for $\T\E (\jj,X)$, since
$\T\E * \h \st \T\E$ by (8.2) and $\gamma_{\om} (f*\va) =
(\gamma_{\om} f)* \gamma_{\om}\va $: Proposition 8.9 is needed for
(2.17), (2.19) for $\T\E$.

\proclaim {Proposition 8.9} For any $\jj$, $X$, $  \T \E (\jj,X)$
satisfies $(L_{\E}) $ and $(\Delta) $,
 $  \T \E_0 (\jj,X) $ satisfies  $(L_{\E_0}) $ and
$(\Delta) $.
\endproclaim
\demo {Proof}  $(L_{\E})$ (respectively $(L_{\E_0}) $):
   Given $\phi \in \E (\jj,X)$  with  $  \Delta_h \phi\in\T \E(\jj,X)$ (respectively $  \T \E_0 (\jj,X) $) for all $h > 0$,
it follows  $M_h \phi \in C_b (\jj,X) $ for all $h > 0$ by (7.2).
Examples 4.15 and Proposition 8.1 (ii) for $C_b$ give
 $M_h \gamma_{\om}\phi \in C_b (\jj,X) $ for all  $\gamma_{\om}$. This implies
 $(M_h (\gamma_{\om}\phi))' \in \E_0(\jj,X)$ for all $ \gamma_{\om}$, and then
   $\,\Delta_h (\gamma_{\om}\phi) \in \E_0(\jj,X)$.
The identity      $(\gamma_{\om} (h) -1)\gamma_{\om} \phi=
(\gamma_{\om} \phi)_{h}- (\gamma_{\om}\phi) - \gamma_{\om}
(h)\gamma_{\om} \Delta_h\phi \,$ gives $\,\,\gamma_{\om}\phi \in
\E(\jj,X)$ (respectively $ \E_0 (\jj,X) $)  for all ${\om}\not
=0$. This and the assumption proves $(L_{\E})$ (respectively
$(L_{\E_0})$).

  $(\Delta) $:    Given $\phi \in L^1_{loc} (\jj,X)$  with  $  \Delta_h \phi\in\T \E(\jj,X)$ (respectively
 $\T \E_0 (\jj,X) $) for all $h > 0$,
    $\phi-M_h\phi\in  \E(\jj,X)$ ( respectively  $\E_0(\jj,X)$) for all $h > 0$ since $ \E(\jj,X)$
( respectively  $\E_0(\jj,X)$) satisfies $( \Delta)$ by Corollary
7.8 (respectively  Proposition 7.3). The identity
$(\phi-M_h\phi)_s - (\phi-M_h\phi)= \Delta_s \phi  - \Delta_s
(M_h\phi)=  \Delta_s \phi-M_h (\Delta_s\phi)  $ and $\T\E \st
\m\T\E $  ( respectively $\T\E_0 \st \m\T\E_0 $ ) of (7.1) for
$\T\E$ (respectively $\T\E_0$ ) show that $(\phi-M_h\phi)_s -
(\phi-M_h\phi) \in \T \E(\jj,X) $ ( respectively  $\T\E_0(\jj,X)$)
for all $s>0$. Therefore $(\phi-M_h\phi) \in \T \E(\jj,X) $
(respectively $  \T \E_0 (\jj,X) $)   by $(L_{\E}) $ (respectively
$(L_{\E_0})$). This proves $(\Delta) $ for $\T\E$ (respectively $
\T \E_0 $) . \P
\enddemo

\proclaim {Remark 8.10}  $\gamma_{\om }\E$, $ \gamma_{\om } \E_0$,
$ \gamma_{\om }Av$, $ \gamma_{\om }Av_0 $ all have $(\Delta)$ for
all $\om \in \r$.
\endproclaim

\noindent Since for example $\T\E= \cap _{\om\in\r}
\gamma_{\om}\,\E$, this generalizes Proposition 8.9.

 \proclaim {Proposition 8.11} For $S\in \h'_{\T\E}
(\r,X)$ (see Example  8.8) one has

 $\qquad \sigma_{B} (S)\,\,  \st \, $  supp $\,\, \widehat {S}$.

\endproclaim

\demo    {Proof}   $ \h'_{\T\E} (\r,X) \st \f' (\r,X)$ by Lemma
8.3, so $\widehat {S}$ is defined. If $\om \not \in $ supp
$\widehat {S}$, there is  $\va\in \h (\r,\cc)$ with $\va (\om)
\not = 0$ and $\va\cdot \widehat {S} =0$. To $\va$ exists $\psi
\in \f (\r,\k) \st \h_{L^1}$ with $\va= \hat {\psi}$. Since
$\widehat {S*\psi} = \hat {\psi} \cdot\widehat {S}$ also for $S\in
\f' (\r,X)$,  $ \psi \in \f (\r,\cc)$ one gets $S*\psi =0$. (8.5)
gives $c_{\om} (S)=0$.   \P
\enddemo

\noindent For the $\A$     as in Examples 8.8, only in  the case
$\A =AP$ or   $\A =S^p AP$ is $S$  uniquely defined by its Fourier
series: If  $S_1, S_2\in \h'_{AP} (\r,X)$ with $c_{\om}(S_1)=
c_{\om}(S_2)$ for $\om\in \r$, with $S= S_1-S_2$ and Proposition
8.7 one gets $c_{\om}(S*\va )=0$ for all $\om$ and $\va\in
\h(\r,\k)$; since $S*\va \in AP(\r,X)$, this implies $S*\va =0$ by
[2, p. 25, VI]. Since this holds for all   $\va\in \h(\r,\k)$, one
gets $S_1=S_2$.

\noindent This is nothing new however: By Theorem 6.1,
$(a)\Leftrightarrow (h)$, every $S \in \h'_{AP}$ can be considered
as an almost periodic function $\Phi_{S} \in AP(\r,Y)$ with $Y=$
locally convex topological complete  vector space $(\h_{L^1})'$,
so by  the general Bochner-von Neumann theory [25], there exist
even summation methods for the Fourier series of $\Phi_{S}$,
converging uniformly on $\r$ to $\Phi_{S}$. Using the uniqueness
of
 a linear, invariant, normalized and continuous (with respect to uniform convergence)
map:  $AP(\r,X)\to X$, one can show

     \smallskip

 (8.6)$ \qquad  \,  \, m_{N} (\Phi_{S}) = \tilde {m }  (S) $,$ \,\,$ (constant distribution)       $\,\,$  for  $S \in \h'_{AP} (\r,X)$.

     \smallskip

\noindent where $m_N$ denotes the Bochner-von Neumann mean [25,
pp. 28-29] on ${AP} (\r,Y)$. With this one gets

     \smallskip

  $\,\,\,\qquad \, \, c_{\om} (\Phi_{S}) = m_N (\gamma_{-\om} \Phi_{S}) =  \gamma_{\om} m_N ( \Phi_{ \gamma_{-\om} S}) = \gamma_{\om} c_{\om} (S) $.

     \smallskip

\noindent so if    $\sum _{\om} s_{n,\om} \, c_{\om} (\Phi_{S})
\,\gamma_{\om} \to\Phi_S$  uniformly  on $\r$ for some summation
method $(s_{n,\om})$ (see [13, Theorem 3.10]), $t=0$ gives

\smallskip

 (8.7)$ \qquad  \,  \, \sum _{\om} s_{n,\om}  \, \, c_{\om} (S)\,\, \gamma_{\om} \to S \,\,\,$  as $n\to \infty$ in  $ \,\,(\h_{L^1})'$,      $  \,\,\, S \in \h'_{AP} (\r,X)$,

\smallskip

\noindent where $(s_{n,\om})$ depends only on $\sigma_{B}(S)$.

\noindent This holds especially for  $f \in \m^n {AP}$ or  $S\in
\widetilde {\m}^n AP$, but here a stronger convergence holds with
Corollary 5.2:

\noindent For  $f \in \m^n {AP}(\jj,X)$, any $\jj$,

           \smallskip

 $\,\,\,\qquad \qquad\sum_\om \,\, s_{n,\om}\,\, c_{\om} (f)\,\, \gamma_{\om} \to f$

\smallskip

 \noindent in the locally convex topology defined on $\m^n L^{\infty} (\jj,X)$
by the seminorms $\,\, ||g||_h :=$

\noindent  $ ||M_{h_n} \cdots M_{h_1} g ||_{\infty}$,
$h=(h_1,...,h_n)$, $h_j > 0 $. Here also the ergodic mean on
$\E(\jj,X)$ can be extended uniquely and continuously to $\cup
_{n=0} ^{\infty} \m^n \E(\jj,X) $. This generalizes theorem 3.10
of [13].

\smallskip

\noindent For   $S\in \widetilde {\m}^n AP(\jj,X)$ similar results hold, with $\widetilde {\m}_{h_j}  $   in the definition of  $||g||_h $.

\smallskip

\noindent With  Proposition 8.11, (8.7) and $\hat {\gamma}_{\om} =
2\pi \delta_{-\om}$  one gets

\smallskip

 (8.8)$ \qquad  \qquad  $ If     $   S \in \h'_{AP} (\r,X)$, then   closure  $ \, {\sigma_{B} (S)} =  $supp $ \widehat {S}$.

\smallskip
\smallskip

\noindent Furthermore with Theorem 2.10 one can show

\smallskip

 (8.9)$ \qquad  \qquad  $      $    \h'_{\A} = \h'_{AP} \oplus \h'_{\A_0}$

\smallskip

\noindent for  $\A= AAP$, $EAP_{rc} $ and $EAP$,  where $AAP_0=
C_0$, $EAP_0 =$ nullfunctions in $EAP$ (see [69,, $W_0 (\r_+, X)$
p. 424]), similarly for $(EAP_{rc})_0$. Therefore to           $ S
\in \h'_{\A}$  there is exactly one
 $   U \in \h'_{AP}$    with $c_{\om} (S) =c_{\om} (U)$, namely the $U$  in
$S= U+V$ of (8.9).

\noindent  Let us finally remark that the ergodic mean $m$ can
also be extended to

$\cup_0^{\infty}  \m^n \E(\jj,X)\,\,$  for $\,\,\jj \not = \r$,

\noindent with corresponding consequences; this will be done in a
forthcoming note.

\head {\bf\S 9 Applications to the study of asymptotic behavior of
solutions of  differential equations}\endhead

 In many cases  the   study of the
  asymptotic behavior of  (a uniformly continuous)  solution $u$ of a  differential equation or system,
that is demonstrating  that $u$  belongs  to  some  given class $\A$
 of functions (for example asymptotic almost periodic  ones), reduces to show the following

\text{\bf {A}}. ``Spectrum''     $sp_{\A} u \st M$, where $M$ is
the $spectrum $ of the equation (see [38, Lemma 1], [65, p. 289],
[10, Theorem 3.3], [13, Lemma 2.4, Theorem 2.5], [32, (3.2)]).

\text {\bf {B}}.   Special cases respectively analogues of the
following  Proposition 9.5 respectively Theorem 9.7 (see [38],
[59, p. 92, Theorem 4, p. 94, Theorem 5], [9, p. 22, Theorem
4.2.6],
 [10, Theorem 2.6, p. 62], [67, Theorem 3.9], [5, Theorem 2.3], [6,
Corollary 3.3], [13, Theorem 4.1], [32, Theorem 2.5]).

\smallskip

 The notion of spectrum $sp_{\A}\phi $ needed here has been introduced in [9, p. 20] in the case $\A\st C_{ub}(\jj,X)$, where $\jj= [0,\infty)$ or $\jj=\r$.
For $\A= AP(\r,\cc)$, $\phi \in C_{ub}(\r,\cc)$,   $sp_{\A}
(\phi)$ was implicitly    used by Loomis [60].  For spectrum of
elements of  $L^1 $-Banach modules see [15].  For $\A= AP(\r,X)$,
 $\phi \in C_{ub}(\r,X)$, $sp_{\A} \phi$  is used in  [59, p. 91 (set of non-almost periodicity)], also
[72, Theorems 2.3, 3.2]. For the cases  $ \A= AP(\r,X)$,
$AAP(\r^{+},X)$,  $EAP (\r^{+},X) $ and $\phi \in C_{ub}(\r,X)$,
$sp_{\A}\phi $ is also defined in [67].

\smallskip

 In Proposition 9.3 below we compare $sp_{\A} \phi$ with   the singular  set $ \sigma_{+} (\phi)$  in $ i\,\r$ of the Laplace transform
 $\tilde  {\phi} \, (\la):=\int_{0}^{\infty} \,e^{-\la t} \phi (t) \,dt $, $\la \in \cc^{+} :=\{\la \in \cc:  \text {Re} \,\la >0\} $,    where
$ \sigma_{+} (\phi):=\{u \in \r: \tilde {\phi}$ is not analytically extendable to some ball $B_{\delta} (iu)=\{ \la \in \cc: |\la-iu| < \delta \}\}$
 (= $Sp_{\r_+}\phi$ in [5, p. 293]). It is obvious that $\sigma_+ (\phi)\st \sigma (\phi)$ the singular set
of the Fourier-Carleman transform of $\phi$ (see [65, (0.46), p.
 19]).

\smallskip

 For $\A=\{0\}$, $\jj=\r$ and $\phi\in L^ {\infty} (\r,X)$, $ sp_{\A} \phi $  $= \sigma (\phi)=$  Beurling spectrum of $\phi\,=\,$
 supp $\hat {\phi}$ (see [55, p. 147, Definition 4.3, p. 171, Theorem], [9, p. 20], [65, p. 19, (0.46),  p. 22, Proposition 0.5]  and (9.4) below).

\smallskip

  We need Spectrum  $sp_{\A}$ also for distributions:

\smallskip

 The assumptions on  $\A$ are,  with   $\jj$, $X$ as in section 1 but always $\k = \cc$:

\smallskip

(9.1)  $\,\,\,\qquad \A$ linear, $\st L^1_{loc}(\jj,X)$,  $\,\,\A$
uniformly closed, $\,\, \A$ $C_{ub}$-invariant.

\proclaim {Definition 9.1} For  $ \A\st  L^1_{loc}(\jj,X)$,
usually with (9.1),   $V \st L^1(\r,\cc)$ and   $S\in  \f' (\r,X)$
we define

        \smallskip

(9.2)   $ \,\,\,\qquad sp_{\A,V} S:=\{\om \in\r: \va\in V,
(S*\va)|\,\jj \in \A$ implies $\hat {\va}(\om)=0 \}$,

        \smallskip

\noindent provided    the convolution  $S*\va$ and $(S*\va)|\,\jj$ are     defined for all $\va\in V$.
\endproclaim

 If $S=\phi  \in  L^{\infty}(\r,X)$ and $V = L^1(\r,\cc)$, we write $sp_{\A}\phi$, the definition of [9].

\smallskip

  The other cases needed here are $ S \in (\h_{ L^1})'(\r,X)$, $V =\h_{ L^1}(\r,\cc)$, and
 $ S \in \f'(\r,X)$, $V =\f(\r,\cc)$; in both cases
 $S*\va$ is well    defined and  a $C^{\infty}$-function and belongs to $ \h_ {L^{\infty}} (\r,X)$ respectively
    $\h_{P}(\r,X) $ ($``P``$ means all derivatives have polynomial growth).

\noindent $sp_{\A,V}$ is always closed (the complement is open), also $sp_{\A} \st sp_{\A,\h_{L^1}}\st sp_{\A,\f}$ if defined.

By (9.4), if $0\in \A$, one has   $  sp_{\A} \phi = sp_{\A, \f}
\phi\st \sigma (\phi) =$ supp $\hat {\phi} $  for  $\phi \in
L^{\infty}(\r,X)$.

 For any $ \A\st L^1_{loc}(\jj,X)$, $\phi \in L^{\infty}(\r,X)$,  $ S \in (\h_{ L^1})'(\r,X)$ one has with $(\h_{L^1})'* \h_{L^1}\st \h_{L^{\infty}}$,

\smallskip

(9.3)  $\,\,\,\qquad     sp_{\A} \phi = sp_{\A_{ub}}\phi,\,\,$
$sp_{\A,\h_{L^1}} S= sp_{\A_{ub},\h_{L^1}} S,\,$

$\,\,\,\qquad\,\,\,\qquad$ with $\A_{ub}=\A\cap  C_{ub}(\jj,X)$.

\proclaim {Proposition 9.2} For   $ \A$ with (9.1),
 $\phi\in L^{\infty} (\r,X) \, $ and $ S \in (\h_{ L^1})' (\r,X)$,  one has

\smallskip

 (9.4)  $\,\,\,\qquad  sp_{\A} \phi = sp_{\A, \h_{L^1}} \phi= sp_{\A,\f}\phi, \qquad $ $sp_{\A, \h_{L^1}} S= sp_{\A,\f} S$.
\endproclaim
\demo {Proof} With Wiener's local inversion theorem
  [18, p. 22, Proposition 1.1.5 (b)]
to   $ f\in L^1 (\r,\cc)$ with $\hat {f} (\om) \not =0 $ there
exist $g \in L^1(\r,\cc)$, $\va \in \f (\r,\cc)$ with $\hat
{f}\cdot\hat {g} =1$ and $\hat {f}\cdot\hat {\psi} =1$ on a
neighborhood of $\om$,
  $ \psi := g*\va \in L^1*\h_{L^1}\st \h_{L^1}$, then $\chi \in \f(\r,\cc)\st \h_{L^1} (\r,\cc)$
  with $ f*\psi *\chi=\chi $, $\hat {\chi} (\om)\not =0$.

  For the $S$-case
one has  to use $(\h_{L^1})' *L^1 *\h_{L^1}$ is well defined  and associative, we omit the details.       \P
\enddemo
\proclaim {Remark}
 $sp_{\A} \phi = sp_{\A, \h_{L^1}} \phi$  holds even for $\A$ when $C_{ub}$-invariant is
replaced by positive-invariant.
\endproclaim
The following result is not used here and may be omitted. But it
relates $sp_{C_{0}} \phi$ used in [9], [10], [13] with the set of
singularities of Laplace transform $\sigma_{+} \phi$ of $\phi \in
C_{ub} (\r,X)$ used in [5]-[7], [17].

 \proclaim {Proposition 9.3}
 Let  $\phi \in C_{ub} (\r,X)$, $\A =C_{0} (\r_+,X)$.
If  $sp_{C_{0}} \phi$ is countable, then  $sp_{C_{0}} \phi
\st\sigma_{+} (\phi)$.

\noindent The inclusion is in general strict: There is $\psi \in
C_{ub} (\r_+)$  with $sp_{C_0} \psi =\emptyset$ but $\sigma_{+}
\psi\not = \emptyset$.
\endproclaim

\demo {Proof} Case $\sigma_+ (\phi)=\emptyset $. We show $sp_{C_0}
\phi= \emptyset$: By Theorem 4.4 [3, p. 847] or [7, Theorem 4.4.1,
p. 275], $P\phi |\,\r_+ \in X \oplus C_0 (\r_+,X)$ giving $\phi
|\,\r_+ \in C_0(\r_+,X) $, by [9, Theorem 1.4.1].  By [9, Theorem
4.2.1],
 $sp_{C_0}\,\, \phi= \emptyset$.

Case  $sp_{C_0}\,\,\phi= \{0\}$. We show   $0 \in
 \sigma_+ (\phi)$:
Assuming $0 \not \in  \sigma_+ (\phi)$,  $P\phi |\,\r_+$ is
bounded by  a result of Ingham  [16, Theorem 1.4] and [39, Theorem
20, p. 66] (see also [4, Remark 3.2, p. 419]. Hence $h M_h \phi=
\Delta_h P \phi |\,\r_+ \in \E_{ub} (\r_+,X)$ for all $h
>0$. Since $h M_h (\phi|\,\r_+) \to \phi |\,\r_+$ as $h \to 0  $
and $\E_{ub} (\r,X) $ is closed, $\phi |\,\r_+ \in \E_{ub}
(\r_+,X)$. By [9, Theorem 4.2.3], $\phi |\,\r_+ \in  C_0
(\r_+,\r)$. This contradicts the assumption and proves
 $0 \in \sigma_+ (\phi)$.

Case $0$ is an isolated point of  $sp_{C_0} \phi  $. We show
 $0 \in \sigma_+ (\phi)$: The assumption implies  $sp_{C_0} \phi*f=\{0\}$
for any $f \in L^1(\r,\k)$ with $\hat{f}(\om)=1$ in some
neighborhood of $0$ and supp $\hat {f} \cap  sp_{C_0}\,\, \phi
=\{0\}$. It follows
 $0 \in \sigma_+ (\phi*f) $ by the above. Since
 $ \sigma_+ (\phi*f) \st  \sigma_+ (\phi*f - \phi ) \cup \sigma_+
 (\phi)\st \sigma (\phi*f - \phi ) \cup \sigma_+
 (\phi)\st$
$\, [\sigma (\phi) \cap $ supp $(\hat {f}-1)] \cup \sigma_+
(\phi)$ (see [65, p. 25, Proposition 06] and $0 \not \in $ supp
$(\hat {f}-1)$, one gets $0 \in \sigma_+ (\phi)$. Since This is
true for any isolated point of $sp_{C_0}\,\, \phi$ and both
$sp_{C_0}\,\,\phi$,
 $\sigma_+ (\phi)$ are closed the inclusion follows.

\noindent The inclusion is strict:  For  $\psi (t) :=\frac {sin
(log (1+t))}{1+t}$, $t\ge 0$, $\psi (t):= 0$, $t\le 0$, one has
$sp_{C_0} \psi =\emptyset$ and its indefinite integral $P\psi$ is
bounded. But by [68, Example 4.1]
  $P\psi |\,\r_+ \not \in X\oplus  C_0 (\r_+,\r)$.  Using Theorem 4.4 of [3],
one can show  $\sigma_{+} (\psi) \not =\emptyset$. \P
\enddemo

\proclaim {Lemma 9.4} If   $ \A$ satisfies   (9.1),
 $\phi\in L^{\infty} (\r,X)$ and $\phi|\, \jj  \in \m \A$, $f \in L^1(\r,\k)$, then $(\phi *f)|\, \jj  \in  \A$.

\endproclaim
\demo {Proof}  Because of (9.1)  and since step functions are
dense in $L^1(\r,\k)$, it is enough to show $(\phi *\chi_{I})|\,
\jj  \in  \A$, $I$ compact interval of $\r$; with $C_{ub}$-
invariance  one can assume $I= [-h,0]$, $h>0$. But then $\phi*
\chi_I= h (\phi*  s_h) =hM_h \phi $ on $\jj$, $\in \A$. \P
\enddemo

\proclaim {Proposition 9.5} If  $\A$ and $ U$  satisfy (9.1),
   $\A$ has $(L_U)$,
   $\phi \in L^{\infty}(\r,X)$,   $ \gamma_{\pm \om} \A \st \A $ and  $\gamma_{-\om} \phi|\, \jj \in \m U $ for all
   $\om
   \in sp_{\A}\phi$, then
   $sp_{\A}\phi$  is perfect (= closed without isolated points).
\endproclaim
\demo    {Proof}     If  ${\om}$ is an isolated point of $sp_{\A}
\phi$, there exists $k\in \f(\r,\cc)$ with  $\hat {k}(\om)=1$ and
supp $\hat {k} \cap sp_{\A} \phi =\{\om\}$. With $\A_{ub}= \A\cap
C_{ub}(\jj,X)$, one has $\A_{ub}\st \m\A_{ub}$ by Corollary 3.3,
so Corollary 2.3(C)  of [13] (valid also for open $\jj$)  gives $
sp_{\A} \phi*k= sp_{\A_{ub}} \phi*k= \{\om\}$. With $\gamma_{-\om}
\A \st \A $  one gets   $ sp_{\A_{ub}}\gamma_{-\om} (\phi*k)\st
\{0\}$. Corollary 2.3 (B) of [13] (also for open $\jj$) can be
applied and gives    $ \Delta_h (\gamma_{-\om} (\phi *k)|\, \jj
\in \A $, for all $h>0$, with (9.1) also if $\alpha >0$ in $\jj$.
Now $\gamma_{-\om} (\phi *k)|\, \jj = (\gamma_{-\om} \phi)
* (\gamma_{-\om} k)|\, \jj \in U $ by Lemma 9.4. $(L_U)$ gives
$\gamma_{-\om} (\phi *k)|\, \jj \in \A  $. With   $\gamma_{\om} \A
\st \A $ one gets $ \phi *k|\, \jj \in \A $.  The definition of
spectrum implies $\hat {k}(\la)=0$, a contradiction. \P
\enddemo

\proclaim {Remark} $(L_{U_{ub}})$ for $\A$ is enough in
Proposition 9.5.
\endproclaim

\proclaim {Lemma 9.6} For any  $ \A\st L^1_{loc}(\jj,X)$ with
(9.1),
 if  $S \in (\h_ { L^1})'  (\r,X)$ and $h> 0$, then     $ \widetilde {M}_h S\in  (\h_ { L^1})'  (\r,X)$ and

 (9.5)  $\,\,\,\qquad sp_{\A,\h _{L^1}} S = \cup_{h>0} sp_{\A,\h_{L^1}}(\widetilde {M}_h S)$.

\endproclaim

\demo {Proof} If  $\va \in \h_ { L^1}  (\r,\cc)$, $M_{-h}\va \in
\h_ { L^1}  (\r,\cc)$ for all  $h>0$. Since $\widetilde {M}_h
S(\va) = S({M}_{-h} \va)$ by (1.7), $M_{-h}$ is continuous on
$\h_{L^1} (\r,\cc)$ and $S \in (\h_ { L^1})'  (\r,X)$, one gets $
\widetilde {M}_h S\in  (\h_ { L^1})'  (\r,X) $. If  $\om\in
sp_{\A,\h_{L^1}} \widetilde {M}_h S$,  $\va \in \h_{L^1}(\r,\cc)$
and $S*\va |\, \jj\in\A$, then $(S*\va)*s_h |\, \jj\in\A$ if    $
\A \st\m \A $. One can prove associativity and commutativity of
the convolutions  in
 $(\h_ { L^1})'  (\r,X)*  L^1(\r,\k)*\h_ { L^1}  (\r,\k)$ (see (2.6)   of [13]), with
   $(\h_ { L^1})'  (\r,X)*  L^1(\r,\k)\st (\h_ { L^1})'  (\r,X)$, so
   $(S*s_h)*\va |\, \jj\in\A$. Since $S*s_h=\widetilde {M}_h S$, one gets
      $\hat {\va}(\om)=0$. This gives $\om\in  sp_{\A, \h_ { L^1}} S$.
For general $\A$  one can use (9.3), since $\A_{ub} \st\m\A_{ub}$,
by Corollary 3.3.

\noindent  Since we do not need $``\st`` $ of (9.5) we omit the
proof.  \P
\enddemo

\proclaim {Theorem 9.7} Assume  $\A$, $U$ satisfy (9.1), $(L_U)$
holds for  $\A$, $m\in \N_0$,  $\phi\in  \m^m L^{\infty} (\r,X)$
with  $sp_{\A,\h_{L^1}} \phi$  at most countable, $\phi |\,\jj \in
\m^m C_u (\jj,X)$ and   $\gamma_{\pm \om} \A \st \A $,
$\gamma_{-\om} \phi|\, \jj \in \m^{q(\om)} U $ with  some $q(\om
)\ge  m$ for all $\om
   \in sp_{\A,\h_{L^1}}\phi$, then  $\phi |\, \jj\in \m^{m} \A$.
\endproclaim
\demo{Proof}
 $\phi\in  \m^m L^{\infty} (\r,X)$     implies
$\phi \in L^1_{loc}(\r,X)$, so all restrictions on $|\,\jj$ are well defined.
 $ \m^m L^{\infty} (\r,X)  \st  \m^{m+1} C_{ub} (\r,X) \st (\h_ {L^1})'(\r,X)$
by Examples 4.7, (2.19) and Proposition 2.9, so $\phi \in (\h_
{L^1})'(\r,X)$, $\sigma:= sp_{\A,\h_{L^1}}\phi$ is well defined.

 $m=0$ : By the trivial  $sp_{\A}\phi\st sp_{\A,\h_{L^1}}\phi$ of (9.4)
and proposition 9.5,  $sp_{\A}\phi$ is perfect,  since $\phi |\,
\jj \in C_u (\jj,X)$  gives  $(\gamma_{-\om}\phi) |\,\jj \in
U_{ub}\st \m U_{ub} $ with Proposition 5.6 and Corollary 3.3. Then
$sp_{\A} \phi$ countable implies
 $sp_{\A}\phi =\emptyset$ (see [1, p. 221, Satz 23]). Since  $\phi |\,\jj \in  C_{ub}(\jj,X)$,
 $sp_{\A}\phi = sp_{\A_{ub}}\phi$, by (9.3). With $\A_{ub}\st\m \A_{ub}$ by Corollary 3.3.
Corollary 2.3 (A) of [13]  also for open $\jj$ can be applied to
$\phi$, $\A_{ub} $, yielding  $ \phi |\,\jj \in \A$.

 $m>0$:  For     $  h=( h_1, \cdots, h_m) \in \r^{m}$ with $h_j >0$
we consider     $\psi :=M_{h} \phi :=M_{h_{m}}\cdots M_{h_{1}} \phi$,  $\in   L^{\infty} (\r,X)$, with
 $\psi |\,\jj \in C_{ub} (\jj,X)$. Assuming
that  $ \gamma_{-\om}  \phi |\, \jj \in  \m^{q(\om)} U$ implies
  $ \gamma_{-\om}  \psi |\, \jj \in  \m U$, which we show below, one gets that
 $sp_{\A}\psi$ is perfect. Now
$ sp_{\A}\psi=  sp_{\A} M_h \phi\st  sp_{\A,\h_{L^1}} M_h\phi$ by
the trivial part of (9.4), Lemma 9.6 gives then inductively $
sp_{\A}\psi \st  sp_{\A,\h_{L^1}}\phi$. So $ sp_{\A}\psi$ is
perfect and countable, thus  is empty (see the case $m=0$). With
(9.3) one can apply Corollary 2.3 (A) of [13] also for open $\jj$
to $\psi$ and $\A_{ub}\st \m \A_{ub}$, obtaining $ \psi|\,\jj
\in\A$. Since $h>0$ is arbitrary, $ \phi |\,\jj\in \m^m \A$
follows.

Now we show   $ \gamma_{-\om}  \psi |\, \jj \in  \m U$: Since  $
L^{\infty} (\r,X)$ has $(\Delta)$ by Proposition  4.9, Proposition
8.1 gives

$F:= M_k (\gamma_{-\om}\phi)\in   L^{\infty} (\r,X)$ for
$\om\in\r$ and
 $  k=( k_1, \cdots, k_m) \in \r^{m}$ with $k_j >0$, $M_k = M_{k_{m}}\cdots M_{k_{1}}$.

\noindent Similarly with Examples 4.7, $F|\,\jj\in C_{ub}(\jj,X)$.
For fixed $\om \in sp_{\A,\h_{L^1}}\phi$, one has       $M_r
F|\,\jj\in U$ for $0< r\in \r^{q-m}$ or $F|\,\jj\in  \m^r U \cap
C_{ub}(\jj,X)\st U_{ub}:= U \cap C_{ub}(\jj,X)  $ by Proposition
5.6. Since $U$ is $C_{ub}$-invariant, $U_{ub}$ is
$C_{ub}$-invariant.
 This and Lemma 9.4 for
$U_{ub}$, $f= \gamma_{-\om}\cdot s_{h_1}$ gives $\int_{0} ^{h_1}
\gamma_{\om} (u_1)  F(u_1 +\cdot)\, d u_1 |\,\jj\in U_{ub}$;
inductively we get

$G(t):=\frac {1}{ h_1\cdots h_m } \int _0 ^{h_m}\cdots \int _{0}^{
h_1} \gamma_{\om} (u_1+\cdots + u_m)   F (u _1 +\cdots + u_m + t)
\, du_1 \cdots d u_m $

\noindent satisfies $G|\,\jj\in U_{ub}$. But Fubini's theorem
gives $G= M_k (\gamma_{-\om} M_h \phi) $, so
  $\gamma_{-\om}\psi |\,\jj= \gamma_{-\om} M_{h} \phi |\, \jj \in U_{ub } \st  \m U_{ub } \st \m U $     as claimed.      \P
\enddemo

\proclaim {Corollary 9.8} If  $\A$, $U$ satisfy (9.1), $(L_U)$
holds for  $\A$, $S\in  \h' _{L^{\infty}} (\r,X)$ with
$sp_{\A,\h_{L^1}} S$  at most countable,
   $\gamma_{\pm \om} \A \st \A $,  $(\gamma_{- \om} S)|\, \jj \in \h' _{U} $  for all
   $\om
   \in sp_{\A,\h_{L^1}} S$, then  $S |\, \jj\in \h'_ {\A}$.
\endproclaim
 Here, for distributions  $T\in \h' (\r,X)$ and any $\jj \st \r$,  $V \st L^1_{loc }(\jj,X)$,   $ T|\, \jj \in \h' _{V} $ means
$(T*\va) |\, \jj\in V$ for all  test functions $ \va  \in  \h (\r,\k)$ in accordance with Definition 1.3.

\demo {Proof} Theorem 9.7, $m=0$ can be applied to $\phi= S*\va
\in \h_{L^{\infty}}\st C_{ub}(\r,X) $, since $\gamma_{-\om}(\phi *
\va)= (\gamma_{-\om}\phi) * (\gamma_{-\om}\va) $ and $
sp_{\A,\h_{L^1}} S *\va \st sp_{\A,\h_{L^1}} S$ with Lemma 9.4 and
associativity  in
 $(\h_ { L^1})'  (\r,X)*  L^1(\r,\k)*\h_ { L^1}  (\r,\k)$ (see Corollary 2.3 (c) of [13]
for  $ S \in  L^{\infty}  (\r,X)) $. \P
 \enddemo

\proclaim {Remark 9.9} In many cases
  $\A$, $U$ satisfy  $(\Gamma )$ and $(\Delta)$ (see Section 4, 7), then
   $\gamma_{\pm \om} \A \st \A $ holds  and     $(\gamma_{- \om} \phi)|\, \jj \in \m^{q(\om)} U  $  respectively $\in \h'_U$  can be replaced by
  $ \phi|\, \jj \in \m^q U  $  for some $q\ge m $ respectively $\in \h'_U$ (Proposition 8.1).
\endproclaim

\proclaim {Remark 9.10} The assumption   $\phi |\,\jj \in \m^m C_u
(\jj,X)$ in Theorem 9.7  is essential, already  for  $m=q=0$, $
X=\cc$ and  $\jj=\r$:
\endproclaim

  There exists $\phi =f_0$ Stepanoff  $S^p$-almost periodic  function  for $1\le p < \infty$, $\in C^{\infty} (\r,\r)$
and bounded so in $\m C_{ub} (\r,\r)$, but not uniformly
continuous on $\r$. One has    $ sp_{\A,\h_{L^1}} f_0 =\emptyset$
for $\A =AP(\r,\cc)$ with Lemma 9.6 and  $S^p AP \st \m AP$, this
$\A$ satisfies $(L_U)$ for $U= C_b (\r,\cc)$ (even   for $U=
L^{\infty}$,  by Proposition 1.7 and [13, p. 120]), so   $\A$, $U$
satisfy (9.1), also  $(\Gamma )$ and $(\Delta)$    by Examples
4.7, 4.15, but $f_0\not \in \A = AP $.

\smallskip

{\bf {Special\,\, cases}} of $\A$ and $U$ in Proposition 9.5,
Theorem 9.7 and Corollary 9.8 where besides (9.1) also the main
assumption $(L_U)$ for $\A$ are fulfilled:

(I)  $\qquad \A= AP(\r,X)$,  $LAP_{ub} (\r,X)$, $UAA(\r,X)$,  $ U=C_{ub}(\r,X)$ and $c_0 \not \st X$.

\noindent For $(L_U)$ see [13, p. 120, case (L.1)].  Here Remark
9.9 applies.

\smallskip

\smallskip

 (II) $\qquad \A$ as in (I),  $U=C_{uwrc}(\r, X)$,  and now $X$
arbitrary; again Remark 9.9  applies  (for  $(L_U)$ see [13, p.
 120, case (L.2)]).

\smallskip

(III) $\qquad  \A $, $\jj$, $X$ arbitrary with (9.1), $(\Delta)$
for $\A$ and $\{$constants$\}$ $\st \A$, $U=\E (\jj,X)$.

\noindent Here  $(L_{U})$  holds practically by the Definition 1.4
of   $(\Delta)$,  only ``$\A$ uniformly closed'' of (9.1) is
needed for $(L_{\E})$.
 By Proposition 4.2, $(\Delta)$ holds for linear $\A \st C_{ub}(\jj,X)$ if $\A$ is only uniformly closed, a generalization of Theorem 3.1.2  of [9],
one gets  generalizations of Theorems 4.2.5 and 4.2.6 there.

\smallskip

(IV) $\qquad  \A $, $\jj$, $X$ arbitrary with only (9.1) and
$(\Delta)$, $U=\E_0$ (Especially  $\A=C_0 $, $(EAP_{rc})_0 $,
$(EAP)_0$).

\noindent Here  $(L_{U})$ follows    again   directly from $(\Delta)$.

\smallskip

(V)  $\qquad \A=S^p AP(\r,X)$, $1\le p<\infty$,   $c_0 \not \st X$, $U=L^{\infty}(\r,X)$ or even $\m^j L^{\infty}(\r,X)$ for some $j\in\N$.

\noindent Here $(L_U)$ can be derived from $(L_{ub})$ for $AP(\r,X)$, this will be treated somewhere else.

\smallskip

 \smallskip

 So for example Theorem 4.1 of [13]   with all its subclasses is  subsumed by the above.

\smallskip

 One also gets:

 Theorem 4.1 of [13] is also true for distribution solutions $\Omega$, if there
everywhere $\m$  is replaced by ``$\widetilde {\m}$''.

\smallskip

 Similarly, this Theorem 4.1 remains true with $\jj=\r $, $b_j \in \h'_{\A}$,
$\Omega_j \in \h'_{L^{\infty}}\cap\h'_{\T\E} $ in (i), $\Omega_j
\in \h'_{C_{uwrc}} $ in (iii), yielding $\Omega_j \in \h'_{\A}$.

\smallskip

 With  the special case (IV) one gets directly the following analogue to the Hardy-Littlewood tauberian theorem
[4, Theorem 2.7 pp. 415-416], [16,  p. 23, Theorem 1.2].
 \proclaim {Remark 9.11}
If   $\A\st L^1_{loc}(\jj,X)$ satisfies $(\Delta)$ and is uniformly closed, then $  \phi \in  \m \A (\jj,X)$   and $P\phi  \in  \E (\jj,X)$ implies $ [P\phi -m(P\phi)] \in \A $.
For $\A= C_0 (\jj,X)  $ this gives $P\phi (t)\to m(P\phi)$ as $|t|\to \infty$.
\endproclaim

\demo{Proof} $(\Delta)$ implies $(L_{\E_0})$ for $\A$ and $
[P\phi- m(P\phi)] \in \E_{0} (\jj,X)$; Proposition 4.2 (i) gives
$(\Delta)$ for $C_0$. \P
\enddemo

 \head{\bf\S 10. Differential-difference  equations }\endhead

\noindent For $m$, $n$, $q\, \in \N$ we consider the following
differential-difference  system operator

(10.1) $\qquad\, Ly := y^{(m)}+ \sum_{j=0}^{m-1}\,\sum_{k=1}^{q}
a_{j,k}\, y^{(j)}_{r_k}$;

\noindent here the $a_{j,k}$ are $n\times n$ matrices $\in
\k^{n\times n}$ with coefficients from $\k$, $= \r$ or $\cc$,
$y=(_1 y, \cdots, _n y)$ is a row vector, the $r_k$ are reals with
$r_1< \cdots < r_q$ and $0$ among them, $y_s (t):= y(t+s)$ where
defined with $_j y : \jj\to X $ with $X$ a Banach space over $\k$,
$\jj= [\al, \infty)$ with $\al \in \r$. With $r:= r_1 \le 0 $ we
define

(10.2) $\qquad\, \jj' := [\al +r, \infty)$.

\noindent $y\in \A$ etc. means $_j y \in A$ for $1\le j\le n$.

  {For} $\jj, X$ as above, $\A \st X^{\jj}$,
 real $r\le 0 $, {we say that} $\A$  {is} $r'$-$invariant$ if $f\in
C^1(\jj',X)$ with $f'|\, \jj\in \A$ and $f_r|\, \jj\in \A$ imply
$f'_r|\, \jj\in \A$, with $\jj'$ of (10.2).

\proclaim {Theorem 10.1} Assume $L$, $r$ as above, $\A$ linear,
positive invariant, $r'$-invariant, $\st$ $L^1_{loc} (\jj,X)$ with
$(\Delta_1)$ and $\A\st \m\A$. If $y \in W^{1,m}_{loc} (\jj',X)$,
$y_r \in \A$ and $L y\in \m\A$, then $y^{(j)}_r \in \A$, $0< j<
m$, $y^{(m)}| \jj \in \m\A$.
\endproclaim
\demo {Proof} $y \in W^{1,m}_{loc} (\jj',X)$ implies $y_{r_k} |
\jj\in W^{1,m}_{loc} (\jj,X)$ for $1\le k\le q$, so $ (Ly)| \jj$
is well defined, $\st L^1_{loc}(\jj,X)$.

$m=1 :$  $y_r \in\A$ and $\A$ positive invariant give $y_{r_k}
|\jj \in \A$, so $ (Ly)| \jj \in \m\A$ implies $y' | \jj \in
\m\A$.

$m\Rightarrow m+1$: One can write $Ly = (Qy)' + \sum _{k=1}^{q}
a_{0,k} y_{r_k}$ with  $Qu = u^{(m)}+
\sum_{j=0}^{m-1}\,\sum_{k=1}^{q} a_{j+1,k}\, u^{(j)}_{s_k}$; with
$z:= (Qy)|\jj$ one gets $z  \in W^{1,1}_{loc} (\jj,X)$ and $z'\in
\m\A$, since $y_r$, $y_{r_k} |\jj \in \A $,  $\st \m\A$ by
assumptions. Now if $h>0$, $\Delta_h z/h = M_h z' \in \A$; so
$\Delta_h z = \Delta_h  (Qy| \jj) =Q (\Delta_h  y) | \jj$, $\in
\A$. $w:= \Delta_h  y\in W^{1,m+1}_{loc} (\jj',X)$ satisfies $w_r=
\Delta_h (y_r)\in \A$ and  $(Q w)| \jj \in \A$, $\st \m\A$,  by
assumptions on $\A$. The induction hypothesis gives $w_r^{(j)}\in
\A $ for $0\le j < m$, thus $(\Delta_h  (y_r))^{(j)} | \jj \in
\A$;  since $Q w \in \A$ by the above, also $(\Delta_h y)^{(m)} |
\jj \in \A$. This implies $\Delta_h (y' | \jj) \in \A$ for $h
>0$; $(\Delta_1)$ gives $(y'- M_1 (y')) | \jj \in \A$; since
$M_1 (y') =  \Delta_1 y$ and $ y | \jj \in \A$, then  $ y_1 | \jj
\in \A$ ( one $r_k =0$) one gets $ y' | \jj \in \A$. $y\in
W^{1,m+1}_{loc}$ and   $m \ge 1$ give  $y\in C^1 (\jj', X)$. $y_r
\in \A$ by assumption, so the $r'$-invariance  of $\A$ gives $
y'_r | \jj \in \A$.

With $u: = y' \in W^{1,m}_{loc} (\jj',X)$ one has therefore $u_r
\in \A$ and $Qu = (Qy)'\in \m\A$. Again using the induction
hypothesis, one gets $y_r ^{(j+1)}= u_r ^{(j)}\in \A$ for $0\le j
< m$, then $y ^{(m+1)}\in \m\A$ with (10.1). \P
\enddemo

\proclaim {Corollary 10.2} If in Theorem 10.1 the $L y\in \m \A$
is replaced by  $L y\in \, \A$ a.e, then $y_r ^{(j)}\in \A$ for $0
< j < m$ and  $y ^{(m)}| \jj \in \,\A$ a.e.
\endproclaim

\noindent Here $f\in \A$ a.e. means there is $g\in \A$ with $f=g$
a.e. on $\jj$.

\smallskip

Theorem 4.1 follows with $Ly = y^{(m)}$  ($r=0$, the
$r'$-invariance holds trivially).

The following  shows that $Ly \in \m \A$ is decisive in Theorem
10.1, $y_r \in \A$ can be relaxed considerably:

\proclaim {Corollary 10.3} If $\A$ and $L$ are as  in Theorem
10.1, $s\in \N_0$, $y\in W^{1,m}_{loc} (\jj',X)$ with only $y_r\in
\m^s \A $,  $L y\in \m \A$ , then $y_r ^{(j)}\in \A$ for $0 \le j
< m$ and $y ^{(m)}| \jj \in \m\A$.
\endproclaim

\proclaim {Corollary 10.4} If $\A$ and $L$ are as  in Theorem
10.1, $k, p, s\in \N_0$ with $m\le p$, $y\in W^{1,p}_{loc}
(\jj',X)$, $y_r\in \m^s \A $,  $L y\in \m^k \A$, then, with $\om
(j)= \om (j, k, m, s):= $ max $\{ 0, j+$ min $\{ k-m, s\}\}$, $y_r
^{(j)}\in \m^{\om (j)}\A$ for $0 \le j < m$ and $y ^{(j)}| \jj \in
\m ^{\om (j)}\A$ for $m \le j \le  p$.
\endproclaim
\demo{Proof} $s\Rightarrow s+1$: For fixed $h> 0$, define $z=M_h y
=(M_{h}\,  _1 y, \cdots , M_{h}\,  _n y) $; then $z\in
W^{1,m+1}_{loc} (\jj',X)\st  C^m(\jj',X)$, $z_r \in \m ^s \A$ and
$ Lz = M_h Ly \in \A\st \m \A$. The induction hypothesis yields
$z_r ^{(j)}\in \A$ for $0 \le j < m$. If $m>1$, then $z'_r \in
\A$. If $m=1$, then $z'| \jj \in \A$ with $Lz \in \A$; since $z\in
C^1 (\jj',X)$ with $z_r \in \A$, the $r'$-invariance of $\A$ gives
also $z'_r \in \A$. This means $ \Delta _h (y_r) \in \A$ for $h>
0$; $ (\Delta _1)$ for  $ \A$ gives $y_r - M_1 (y_r)\in \A$; since
$z_r \in  \A$, $h = 1$ gives $M_1(y_r)  =
                            (M_1(y))_r  \in  \A$,
 so $y_r
\in \A$. Theorem 10.1 gives then the conclusion of Corollary 10.3.

Corollary 10.4 follows inductively with  Corollary 10.3 and Lemma
2.3, we leave the details to the reader. \P
\enddemo

\proclaim {Examples 10.5} All the  $\A$  the following $\A$ are
admissible in Theorem 10.1 and its Corollaries: $C_0$, $L^p_w$ for
$1\le p\le \infty$, $AP$, $AAP$, $UAA$, $S^p$- $AP$, $EAP$,
various $PAP$ versions, $\T\E$, $\E$, $C_{ub}$, $C_u$, $O(w)$,
also all $\la$-classes [9, Definition 1.3.1], [13, p. 117]:
Examples 3.4/5/6/13, 4.7/9/10, 7.3/8/10/11/12,
                      8.9.
\endproclaim

\proclaim {Example 10.6} Without " positive-invariant" Theorem
10.1 becomes false already for $Ly = y''$:

$\A = \{ f \in C_{ub} ([0,\infty),\r): f(0)=0 \}$ is a Banach
space  with the supremum norm which satisfies $(\Delta)$ by
Proposition 4.2 (i). For the function $y(t) = (\sin\,t ) \cos
\frac{t^2}{t^2 +1}$ one has $y, y'' \in \A$ but $y'\not \in \A$.

Another such example would be $\A =C_0 + X\, t$.

There are also  such $r'$-invariant $\A$.
\endproclaim

\proclaim {Example 10.7}  Theorem 10.1 becomes also false  without
$r'$-invariance of $\A$: $Ly = y'' + y_ {-1}$, $\jj =[1,\infty)$,
$y'$ piecewise linear with $y'(0)=1$, $y'\equiv 0$ on
$[1/2,\infty)$, $\A = C^1 (\jj, \k)$ or $C_{ub} (\jj, \k)$.
\endproclaim

\proclaim {Remarks  10.8}

(a) $y_r ^{(j)}\in \A$  implies $y _{r_k} ^{(j)}| \jj\in \A$,
especially  $y ^{(j)}\in \A$, for $1 \le k \le q$, $0 \le j $.

 (b) If all $r_k \ge 0 $, i.e. $L$ is
only advanced or  even an ordinary differential operator, the
$r'$-invariance holds trivially.

 (c) For linear $\A$ " positive-invariant" and $\A\st \m\A$ are coupled:

(i) $\A$  positive-invariant  with $(\Delta)$ implies $\A\st
\m\A$.

(ii) If $\A$  is uniformly closed $\st C_u$, then $\A\st \m\A$ and
positive-invariant are equivalent.

(d) $\A$ linear, positive-invariant with $(\Delta)$  and $\A\st
\m\A$ does not imply $\A$ $r'$-invariant: Example 10.7.

(e) If $\A \st C(\jj,X)$, the a.e can be disregarded in the
conclusion of  Corollary 10.2, then $y | \jj \in C^m  (\jj,X)$.

(f) For $\A \st L^p_w$ see [51] and the references there, then
also neutral perturbations  can be treated with variable operator
valued $a_{j,k}$.

(g) In  Theorem 10.1 and Corollary 10.3 the assumption $Ly \in
\m\A$ cannot be weakened to $Ly \in \m^ {p(m)}\A$ with $p(m) >1$,
not even for $Ly = y^{(m)}$, $m=2, 3, 4 $, $\A = AP$.


(h) the case $\A = C_{ub}$ is also treated in [12,  Corollary
3.3a], [51], then variable $a_{j,k}\in C_{ub}$ are admissible ($y$
bounded implies $y\in \m C_{ub}$, our Corollary 10.3 applies).

(i) For  $\la$-classes  $\A$ an Esclangon-Landau result for
general neutral system has been obtained in [13, Corollary 2.7].

(j) Theorem 10.1 and  its corollaries hold also for $a_{j,k}$ with
components $a: X\to X$ bounded linear operators with $a (\A) \st
\A$.

(k) The results of this section hold also for $\jj=\r=\jj'$, with
the same proofs.
\endproclaim

For $\A$ which do not satisfy $(\Delta)$ we have

\proclaim {Theorem  10.9} Assume  $L$, $r$, $\jj$, $\jj'$, $X$ as
before Theorem 10.1, $\A$ linear, positive invariant, uniformly
closed , $r'$-invariant, $\st$ $L^1_{loc} (\jj,X)$, $s$ and $s^*
\in \N_0$.  Then $y \in W^{1,m}_{loc} (\jj',X)$ with $y_r \in
[\m^s\A] \cap\m^{s^*}  C_u$ and $L y\in \m C_u (\jj,X)$ imply
$y^{(j)}_r \in \A$, $0\le j< m$, $y^{(m)}| \jj \in \m\A$.
\endproclaim
\demo{Proof} The class $C_u (\jj,X)$ satisfies all the assumptions
of Theorem  10.1 by  Example 10.5, so Corollary 10.3 gives
$y^{(j)}_r \in C_u$ for  $0\le j< m$. $y_r \in \m^s \A$ gives
$y^{(j)}_r \in\m^{s+j} \A$ for  $0\le j\le m$ by Lemma 2.2, since
all $\m^{k} \A$ are also linear positive-invariant. $(\m^{s+j} \A)
\cap C_u \st \A $ of Proposition 5.6 gives $y^{(j)}_r, y^{(m-1)}|
\,\jj \in\A $ for $0\le j < m$,  Lemma 2.2 then $y^{(m)}| \,\jj
\in\m\A $. \P
\enddemo

\proclaim {Corollary 10.10} If $\A$,  $L$, $s$  are as  in Theorem
10.9 with  $\A\st \m C_u (\jj,X)$, then  $y\in W^{1,m}_{loc}
(\jj',X)$ with  $y_r\in \m^s \A $ and   $L y | \jj \in  \A$ imply
$y_r ^{(j)}\in \A$ for $0 \le j < m$, $y ^{(m)}| \jj \in \A$.
\endproclaim

The above can be generalized as in Corollary 10.4.

 \proclaim {Example 10.11} $\A= AP\cdot  e^{it^2}$  does not satisfy
$(\Delta_1)$ by  Example 4.21. This $\A$ is linear,
positive-invariant and uniformly closed, $\st C_b \st \m C_{ub}$,
Theorem 10.9 and Corollary 10.10 can still be applied.
\endproclaim

Similarly $\A= S^p AP$ is possible in Theorem 10.9 and Corollary
10.10 for $1\le p<\infty$, since $\m C_{ub}$ is closed with
respect to $S^p$-norm.

Also, all linear, positive-invariant, uniformly closed $\A \st
S^p_b$ with  $1 < p \le \infty$ are admissible,  since $S^p_b\st
\m C_{ub}$,

By Example 10.6  the "positive-invariant" in Theorem 10.9 and the
following are essential.

\head {\S 11 Open questions}\endhead

1- What in \S 2 can be extended to the half line?

2- Does there exist $\A$ with $(\Delta_1)$, but without
$(\Delta)$?

3-Do there exist $\A \st C_u$, with $\A$ positive-invariant
uniformly closed $\st \m \A$, but without $(\Delta)$?

4- Do all finite-dimensional $\A$ satisfy $(\Delta)$? (Yes for
dimension $\A = 1$)

5- Does $C_{wrc}$ have $(\Delta)$? Similarly for
$(L^{\infty}_{wrc})$, $L^{\infty}_{rc}:=$  the space of a.e.
bounded measurable functions with (weakly) relatively compact
range.

6- When does $(\Delta)$ for $\A$ imply $(\Delta)$ for
$\gamma_{\om}\A$, where $\om \in \r$? (See Remark 8.10).

7- When does $(\Delta)$ for $\m\A$ imply $(\Delta)$ for $\A$? For
example does $\m (AP\cdot e^{it^2})$ satisfy $(\Delta)$?

8-Does there exist a general theorem, which subsumes (most of) the
results of \S 4 with their ad hoc proofs?

9- Is Bohr-ap equivalent Bochner-ap, without  von Neumann's
countability axiom $(A_0)$ (see  Lemma 6.4 and Remark 6.3)?

10-What of Proposition 7.1 is true for $Av$, $Av_{ub}$, $Av_0$,
$\T Av$?

11- Do there exist $U$, $V$ with $(\Delta)$ and $U\cap V=\{0\}$
but $U +V$ does not satisfy  $(\Delta)$?

12- When does  $(\Delta)$ for $\A$ imply $(\Delta)$ for $X+\A$?
(See Corollary  7.7).

13-Does $Var (\jj,X)$ satisfy $(\Delta)$?

 14-For what $\A(\r,X)$ is $\m \A (\r,X)|\,\jj = \m \A(\jj,X)$
? (See Example 5.5).

15-  Does there exist a simple proof for $AAP\cdot e^{it^2}$ has
$(\Delta)$? (See Remark  4.23).

16- Is there an example of an $\A  \st  \m\A$ with  $\A +
\A'_{Loc}$ strictly
      $\st  \m\A$ (see Proposition 5.1)?

 For any  $\phi : \jj \to X$,  $\phi$ will be called $Maak\,\,
ergodic$
if to any positive $ \e$ there are $a  \in X$, $n  \in N$, $s_j
\in [0,\infty)$ such that

        $ || (1/n) \sum^n_{j=1}  \phi (s_j+t)  -  a || < \e $  for all  $ t  \in  \jj$;

\noindent the set of all these $\phi$  will be denoted by $\E
M(\jj,X)$ (see (see [61, p. 34, Mittelwertsatz])).

One can show $\E M(\jj,X)\cap L^{\infty}(\jj,X) \st \E (\jj,X)\cap
L^{\infty}(\jj,X)$.

 17- Is the inclusion $\E M(\jj,X)\cap L^{\infty}(\jj,X) \st \E
(\jj,X)\cap L^{\infty}(\jj,X)$ strict?

18- Does $\E M(\jj,X)\cap L^{\infty}(\jj,X)$ satisfy $(\Delta )$?

19- Is $ e^{it^2}$ Maak ergodic?

20-For what $\jj$, $X$ is $EAP(\r,X)|\,\jj \st EAP(\jj,X)$ strict?
( For  $\N$
 it is known  that
$ EAP(\z,\r)\, |\,\N \st\, EAP(\N,\r)$ is strict (see [20, p. 231,
 Example 5.1.15]).

21-For what $\jj$, $X$ is $EAP_0 \st  \E_n$ valid, or even $|f|
\in EAP_0(\jj,\r)$ whenever $f\in EAP_0 (\jj,X)$ (true for $X=\k$,
any $\jj$ by [20, p. 157, Theorem 4.3.13])?

\indent School of Math. Sci., P.O. Box No. 28M, Monash University,
 Vic. 3800.

\indent E-mail "bbasit\@vaxc.cc.monash.edu.au".

\indent Math. Seminar der  Univ. Kiel, Ludewig-Meyn-Str., 24098  Kiel, Deutschland.

\indent E-mail "guenzler\@math.uni-kiel.de".

\enddocument